\def\VersionDateTime{24/December/2025. Version $2.0$}
\documentclass[12pt, leqno]{amsart}
\usepackage{amscd, amsmath, amssymb, amsfonts}
\usepackage{latexsym}
\usepackage{mathrsfs}

\usepackage{fullpage}
\usepackage{microtype}
\usepackage{colonequals}

\usepackage{graphicx}
\usepackage{tikz}
\usepackage{tikz-cd}
\usetikzlibrary{calc, positioning}

\newtheorem{Theorem}{Theorem}[section]
\newtheorem{Proposition}[Theorem]{Proposition}
\newtheorem{Lemma}[Theorem]{Lemma}
\newtheorem{Corollary}[Theorem]{Corollary}

\theoremstyle{definition}
\newtheorem{Definition}[Theorem]{Definition}
\newtheorem{Remark}[Theorem]{Remark}
\newtheorem{Example}[Theorem]{Example}

\newcommand{\TT}{{\mathbb{T}}}
\newcommand{\GG}{{\mathbb{G}}}

\newcommand{\ZZ}{{\mathbb{Z}}}
\newcommand{\QQ}{{\mathbb{Q}}}
\newcommand{\RR}{{\mathbb{R}}}
\newcommand{\CC}{{\mathbb{C}}}
\newcommand{\PP}{{\mathbb{P}}}
\newcommand{\OO}{{\mathcal{O}}}

\newcommand{\Hom}{\operatorname{\operatorname{Hom}}}

\newcommand{\GL}{\operatorname{GL}}

\newcommand{\Ker}{{\operatorname{ker}}}

\newcommand{\Pic}{\operatorname{Pic}}

\newcommand{\rank}{\operatorname{rank}}

\newcommand{\rest}[2]{\left.{#1}\right\vert_{{#2}}}  % restriction of #1 to #2

\newcommand{\Spec}{{\operatorname{Spec}}}

\newcommand{\trop}{\operatorname{trop}}

\newcommand{\an}{{\operatorname{an}}}

\newcommand{\id}{\operatorname{id}}

\newcommand{\Supp}{\operatorname{Supp}}

\newcommand{\val}{\operatorname{val}}

\newcommand{\argmin}{\operatorname{arg\,min}}

\newcommand{\ndot}{\raisebox{.4ex}{.}}

\usepackage{mathrsfs}

\newcommand{\Proof}{{\sl Proof.}\quad}
\newcommand{\QED}{{\unskip\nobreak\hfil\penalty50\quad\null\nobreak\hfil
{$\Box$}\parfillskip0pt\finalhyphendemerits0\par\medskip}}
\begin{document}

%%%%%%%%%%%
%% Title                 %%
%%%%%%%%%%%

%\title{Ample divisors on tropical toric varieties}
\title{Effective faithful tropicalizations and embeddings for abelian varieties}
\author{Shu Kawaguchi}
\address{Department of Mathematics, Graduate School of Science, Kyoto University, Kyoto 606-8502, Japan}
\email{kawaguch@math.kyoto-u.ac.jp}
\author{Kazuhiko Yamaki}
\address{Mathematical Institute, Graduate School of Science, Tohoku University, Sendai 980-8578, Japan}
\email{kazuhiko.yamaki.d6@tohoku.ac.jp}
\date{\VersionDateTime}
\subjclass[2020]{14K25 (Primary); 14T25, 14G22 (Secondary)}
\keywords{abelian variety, line bundles, faithful embeddings, skeletons, faithful tropicalization, theta functions, Berkovich space, tropical geometry}

%% 14T05   	Tropical geometry
%% 14C20   	Divisors, linear systems, invertible sheaves
%% 14G22   	Rigid analytic geometry

\newcommand{\Proj}{\operatorname{\operatorname{Proj}}}
\newcommand{\Prin}{\operatorname{Prin}}
\newcommand{\Rat}{\operatorname{Rat}}
\newcommand{\zero}{\operatorname{div}}
\newcommand{\Func}{\operatorname{Func}}
\newcommand{\red}{\operatorname{red}}
\newcommand{\pr}{\operatorname{pr}}
\newcommand{\can}{\operatorname{can}}
\newcommand{\relin}{\operatorname{relin}}

%\definecolor{skyblue}{rgb}{0.0, 0.45, 0.73}
\def\red#1{\textcolor{red}{#1}}

\begin{abstract}
Let $A$ be an abelian variety over an algebraically closed field $k$ that is complete with respect to a nontrivial nonarchimedean absolute value. Let $A^{\an}$ denote the analytification of $A$ in the sense of Berkovich, and let $\Sigma$ be the canonical skeleton of $A^{\an}$.  In this paper, we obtain a faithful tropicalization of~$\Sigma$ by nonarchimedean theta functions, giving a tropical version of the classical theorem of Lefschetz on abelian varieties. Key ingredients of the proof are~(1) faithful embeddings of tropical abelian varieties by tropical theta functions and~(2) lifting of tropical theta functions to nonarchimedean theta functions, and they will be of independent interest. For~(1), we use some arguments similar to the case of complex abelian varieties as well as Voronoi cells of lattices. For~(2), we use Fourier expansions of nonarchimedean theta functions over the Raynaud extensions of abelian varieties. 
\end{abstract}

\maketitle

\section{Introduction}
A line bundle $L$ on a projective variety $X$ is said to be very ample if it admits a closed embedding into projective space; more precisely if $L$ has global sections $s_0 , \ldots , s_m$ such that the base loci of those sections have empty intersection and that the morphism $X \to \PP^m_k$ given by $x \mapsto (s_0 (x) : \cdots : s_m (x))$ is a closed embedding. There are some classically known sufficient conditions for a line bundle to be very ample. When $X$ is a smooth projective curve of genus $g$, then $L$ is very ample if $\deg (L) \geq 2g + 1$ (see e.g. \cite[Cor.~IV.3.2]{Ha}). When $X$ is an abelian variety and if $L$ is ample, then it is known as the Lefschetz theorem that $L^{\otimes d}$ is very ample for any $d \geq 3$ (see e.g. \cite[\S4.5]{BL}).

Recently, tropical geometry has seen much development. Through the tropicalization,  algebraic varieties are mapped to tropical varieties. Now, one can consider such kinds of embedding problems in tropical geometry. 
There are at least two versions: one is the faithful tropicalization problem of a skeleton of the analytification of an algebraic variety by using global sections of algebraic line bundles; the other is the faithful embedding problem of a tropical variety by regular global sections of a tropical line bundle.

The aim of this paper is to investigate those two problems for abelian varieties. We will establish tropical versions of the Lefschetz theorem for abelian varieties.

\subsection{Main results and ideas of proofs}

Let us briefly explain our main results and ideas of proofs.

\subsubsection{Faithful tropicalizations} \label{subsection:FT:intro}

Let $k$ be an algebraically closed field that is complete with respect to a nontrival nonarchimedean absolute value. Let $A$ be a projective variety over $k$, and let $A^{\an}$ be the analytification of $A$ in the sense of Berkovich. Then one often considers a skeleton $\Sigma$ in $A^{\an}$, which is a closed subset of $A^{\an}$ and is naturally a polyhedral set with an integral structure. In many cases, a skeleton is associated with a model of $A$ and hence is not unique, but in some cases, there is a canonical choice of a skeleton. For example, if $A$ is an abelian variety, one can define the canonical skeleton of $A^{\an}$.

Let $\tilde{\varphi}\colon A \to \PP^m_k$ be a morphism and let $\tilde{\varphi}^{\an}\colon A^{\an} \to (\PP^m_k)^{\an}$ be its analytification. Let $\TT\PP^m$ denote the $m$-dimensional tropical projective space. By the composition with the valuation map $\trop\colon (\PP^m_k)^{\an} \to \TT\PP^m$, we have a map $\trop \circ \tilde{\varphi}^{\an}\colon A^{\an} \to \TT\PP^m$, denoted by $\tilde{\varphi}_{\trop} : A^{\an} \to \TT\PP^m$. Let $\Sigma$ be a skeleton of $A^{\an}$. We say that {\em $\tilde{\varphi}$ faithfully tropicalizes $\Sigma$} if $\tilde{\varphi}_{\trop}$ restricts to a homeomorphism onto its image preserving the integral structures, in which case we say that  $\rest{\tilde{\varphi}_{\trop}}{\Sigma}\colon \Sigma \to \TT\PP^m$ is a {\em faithful embedding} of the skeleton $\Sigma$. 

Let $\tilde{L}$ be a line bundle on $A$. We say that {\em $\tilde{L}$ admits a faithful tropicalization for $\Sigma$} if there exist global sections $s_0 , \ldots , s_m$ of $\tilde{L}$ such that $\bigcap_{i=0}^m \Supp (\zero (s_i)) = \emptyset$ and $(\tilde{\varphi}_{s_0 , \ldots , s_m})_{\trop}\colon A^{\an} \to \TT\PP^m$ faithfully tropicalizes $\Sigma$, where $\Supp (\zero (s_i))$ denotes the support of the zero divisor $\zero (s_i)$ of the global section $s_i$ and $\tilde{\varphi}_{s_0 , \ldots , s_m}\colon A \to \PP^m_k$ is the morphism given by $x \mapsto (s_0(x): \cdots: s_m(x))$. 

When $A$ is a smooth projective curve, it is shown in \cite{KY2} that there is a concrete positive integer $C(g)$ depending only on the genus $g$ of $A$ such that if $\deg (\tilde{L}) \geq C(g)$, then for any skeleton $\Sigma$ of $A^{\an}$, $\tilde{L}$ admits a faithful tropicalization of $\Sigma$. This is a tropical version of the classical result for smooth projective curves that if $\deg (\tilde{L}) \geq 2g+1$, then $\tilde{L}$ admits a closed embedding, namely, $\tilde{L}$ is very ample.

Our first goal is to give the same type of result for abelian varieties.

\begin{Theorem} [Corollary~\ref{cor:main:FT1}~(2)] 
\label{thm:intro:main:FT}
Let $A$ be an abelian variety over $k$ and let $\tilde{L}$ be an ample line bundle on $A$. Let $d$ be an integer. Then if $d \geq \max \{ 3 , 2 (n-1)! \}$, then $\tilde{L}^{\otimes d}$ admits a faithful tropicalization for the canonical skeleton $\Sigma$ of $A^{\an}$.
\end{Theorem}

Theorem~\ref{thm:intro:main:FT} is seen as a tropical version of the Lefschetz theorem for abelian varieties, asserting that the $d$th power of an ample line bundle on an abelian variety is very ample if $d \geq 3$, while the condition on the power $d$ is not the same. 

In fact, we will prove a more precise theorem than Theorem~\ref{thm:intro:main:FT}: using the polarization type $(d_1, \ldots, d_n)$ of $\tilde{L}$ on the torus part,  we show that if $d_1 \geq \max \{ 3 , 2 (n-1)! \}$ and $\tilde{L}$ is basepoint free, then $\tilde{L}$ admits a faithful tropicalization for the canonical skeleton; see Theorem~\ref{thm:main:FT1}~(3). We note that if we are concerned only with whether or not $\tilde{L}^{\otimes d}$ associates a homeomorphism from the canonical skeleton onto its image (disregarding the integral structures), 
then we show in  Corollary~\ref{cor:main:FT1}~(1) that the condition $d \geq 3$ suffices. 

\smallskip
There are two  key ingredients of the proof of Theorem~\ref{thm:intro:main:FT}. One is 
faithful embeddings of tropical abelian varieties. The other is
lifting of theta functions.  

\subsubsection{Faithful embeddings}
\label{subsec:intro:faithful:embeddings}
Let $X$ be a rational polyhedral space and let $\varphi\colon X \to \TT\PP^m$ be a piecewise $\ZZ$-affine map. As we explain above, $\varphi$ is called a {\em faithful embedding} if it is a homeomorphism onto its image preserving the integral structures. On $X$, one can consider a tropical line bundles $L$ and its regular global sections. We say that {\em $L$ admits a faithful embedding} if there exist regular global sections $s_0 , \ldots , s_m$ of $L$ such that $\bigcap_{i=0}^m \Supp (\zero (s_i)) = \emptyset$ and the map $\varphi_{s_0 , \ldots , s_m}\colon X \to \TT\PP^m$ given by $x \mapsto  (s_0 (x) : \cdots : s_m (x))$ is a faithful embedding. 

Among a class of rational polyhedral spaces is a {\em tropical abelian variety}, which is a real torus with an integral structure admitting a polarization (see Definitions~\ref{def:real:torus} and~\ref{def:trop:av}.) Our second goal is to give faithful embedding results for tropical abelian varieties.  

\begin{Theorem}[Theorem~\ref{thm:main:1}~(3)] 
\label{thm:intro:main:FE}
Let $X$ be a tropical abelian variety of dimension $n$ and let $L$ be a positive line bundle on $X$. Let $d$ be an integer. Then if $d \geq \max \{ 3 , 2 (n-1)! \}$, then the $d$th multiple $dL$ of $L$ admits a faithful embedding.
\end{Theorem}

Theorem~\ref{thm:intro:main:FE} is also seen as a tropical version of the classical Lefschetz theorem for abelian varieties, while the condition for the multiple $d$ 
is not the same. 
%\red{the lower bound that appears in} the condition for the power $d$ is not $3$. 
In terms of the polarization type $(d_1, \ldots, d_n)$ of $L$, the condition $d \geq \max \{ 3 , 2 (n-1)! \}$ is replaced by $d_1 \geq \max \{ 3 , 2 (n-1)! \}$; see Remark~\ref{remark:twoexpressionofFE}.  
We also show in Theorem~\ref{thm:main:1}~(1) that if we disregard the integral structures, then for any $d \geq 3$, there are regular global sections of $dL$ that give a homeomorphism from the tropical abelian variety onto its image. 

\smallskip
Theorem~\ref{thm:intro:main:FE} implies the following, in particular. 
\begin{Corollary}[see Corollary~\ref{cor:trop:abel:proj}]
Let $X$ be a real torus with an integral structure. Then $X$ is a tropical abelian variety if and only if there exists a tropical line bundle on $X$ that admits a faithful embedding. 
\end{Corollary}

Let us briefly state some ideas of the proof of Theorem~\ref{thm:intro:main:FE}. A tropical abelian variety of dimension $n$ is written as $X \colonequals N_\RR/M^\prime$, where $N$ is a free $\ZZ$-module of rank $n$ and $M'$ is a lattice in $N_{\RR} := N \otimes \RR$. Let $\pi\colon N_\RR \to X$ be the natural projection. Let $L$ be a tropical line bundle on $X$. Then $\pi^{\ast} (L)$ is a trivial tropical line bundle on $N_{\RR}$. A nontrivial regular global section of $L$ is then identified with a concave piecewise $\ZZ$-affine function on $N_{\RR}$ that is quasi-periodic with respect to $M^\prime$. We call such a function (and the constant function $+ \infty$) a {\em tropical theta function} for $L$. Let $\TT = \RR \cup \{ + \infty \}$ be the tropical semifield with binary operations $\min$ and $+$. Then the set of tropical theta functions for $L$ becomes a $\TT$-semimodule. 

Assume that $L$ is positive. Then one defines the type $(d_1 , \ldots , d_n)$ of $L$, where $d_1 , \ldots , d_n$ are positive integers with $d_i \mid d_{i+1}$ for all $i=1 , \ldots , n-1$. Set $D:= d_1 \cdots d_n$. 
We will use a certain tropical basis $\vartheta_1 , \ldots , \vartheta_D$ 
of the $\TT$-semimodule of tropical theta functions for~$L$. Each $\vartheta_i$ equals, up to an additive real constant, the tropical theta function $\vartheta_i'$ that is explicitly constructed by Mikhalkin--Zharkov \cite{MZ} and Sumi \cite{Sumi}.When we regard $\vartheta_i$ and $\vartheta_i'$ as regular global sections of $L$, we write $\bar{\vartheta}_i$ and $\bar{\vartheta}_i'$, respectively. We will prove that $\varphi_{\bar{\vartheta}_1 , \ldots , \bar{\vartheta}_D}\colon X \to \TT\PP^{D-1}$ is a faithful embedding to obtain Theorem~\ref{thm:intro:main:FE}. Note that, since $\vartheta_i$ equals $\vartheta_i'$ up to an additive real constant, it suffices to show that $\varphi_{\bar{\vartheta}_1' , \ldots , \bar{\vartheta}_D'}$ is a faithful embedding.   (For the reason why we consider both the $\vartheta_i$ and the $\vartheta_i'$, see Remark~\ref{remark:intro:rational}.) 
Our proof of this assertion consists of two steps: (1) we show that if $d_1 \geq 3$, then the map  $\varphi_{\bar{\vartheta}_1' , \ldots , \bar{\vartheta}_D'}\colon X \to \TT\PP^{D-1}$ is injective; (2) we show that if $d_1 \geq 2 (n-1)!$, then $\varphi_{\bar{\vartheta}_1' , \ldots , \bar{\vartheta}_D'}$ is unimodular. It follows from those results that $L$ admits a faithful embedding if $d_1 \geq \max \{ 3 , 2(n-1)!\}$. The proof of the injectivity (1) is somewhat tricky, using some arguments similar to the case of complex abelian varieties. The proof of unimodularity (2), we use the Voronoi cells of lattices. 

\subsubsection{Lifting} \label{subsubsection:lifting}
\smallskip
Let $\Sigma$ be the canonical skeleton of an abelian variety $A$. Note that $\Sigma$ is  naturally a tropical abelian variety. Let $\tilde{L}$ be an ample line bundle on $A$, and let $L$ be the line bundle on the skeleton $\Sigma$ induced from $\tilde{L}$. It is known that we can identify the global sections of $\tilde{L}$ with the nonarchimedean theta functions on the Raynaud extension $E^{\an}$ of $A$; see Section~\ref{section:uniformization} for the Raynaud extensions. 
Foster--Rabinoff--Shokrieh--Soto \cite{FRSS} shows that a global section of $\tilde{L}$, which is regarded as a nonarchimedean theta function on  $E^{\an}$, corresponds by tropicalization to a tropical theta function for $L$. 

Let $\Gamma$ be the value group of $k$. Then we have a notion of $\Gamma$-rationality for tropical theta functions. Let $H^0(X, L)_\Gamma$ be the set of $\Gamma$-rational tropical theta functions for $L$. We set $\overline{\Gamma} := \Gamma \cup \{ + \infty \}$. Then $H^0(X, L)_\Gamma$ is a $\overline{\Gamma}$-semimodule. 
Further, we see that a tropical theta function that comes from a nonarchimedean theta function by tropicalization is a $\Gamma$-rational theta function (cf. Subsection~\ref{subsec:trop:na:theta:functions}). Thus the tropicalization by Foster--Rabinoff--Shokrieh--Soto induces a map 
\begin{equation}
\label{eqn:map:nonarch:to:trop:theta}
H^0(A, \tilde{L}) \to H^0(X, L)_\Gamma, \quad f \mapsto  f_{\rm trop}.
\end{equation}
To put it loosely, this asserts that {\em tropicalization of a nonarchimedean theta function is a tropical $\Gamma$-rational theta function.} 

Take tropical theta functions $\vartheta_1 , \ldots , \vartheta_D$ on $N_\RR$ which are used in the above faithful embedding of $X$. It follows from the definition that they are $\Gamma$-rational. Furthermore, it turns out that they lift to nonarchimedean theta functions, namely, the following theorem holds. 

\begin{Theorem}[see Corollary~\ref{cor:preciselift} for the precise statement]
\label{thm:intro:liftingtheorem} 
For $1 \leq i \leq D$, there exists a $\tilde{\vartheta}_i \in H^0(A, \tilde{L})$ such that $\tilde{\vartheta}_i$ is mapped to $\vartheta_i \in H^0(X, L)_\Gamma$ by \eqref{eqn:map:nonarch:to:trop:theta}. 
\end{Theorem}

To prove Theorem~\ref{thm:intro:liftingtheorem}, we use the Fourier expansions of nonarchimedean theta functions. Such a construction is efficient because $\vartheta_i$ is given in the form of ``tropical Fourier expansion''; see the last part of Subsection~\ref{subsection:construction:tropicaltheta}. 
We will find the Fourier coefficients of an expected nonarchimedean theta function $\tilde{\vartheta}_i$ in order that the tropicalization of the Fourier coefficients coincide with the tropical Fourier coefficients of $\vartheta_i$. Then we show that  $\vartheta_i$ lifts to a nonarchimedean theta function $\tilde{\vartheta}_i$. 

\smallskip
Once we have Theorems~\ref{thm:intro:main:FE} and~\ref{thm:intro:liftingtheorem}, 
it is not difficult to see from the definition of the canonical skeleton $\Sigma$ that $\varphi_{\tilde{\vartheta}_1,\ldots,\tilde{\vartheta}_D}$ faithfully tropicalizes $\Sigma$, which completes the proof of 
Theorem~\ref{thm:intro:main:FT}. 

\begin{Remark} \label{remark:intro:rational}
The tropical theta functions $\vartheta_i'$ constructed by Mikhalkin--Zharkov and Sumi are not $\Gamma$-rational, and hence we cannot lift them to nonarchimedean theta functions. That is the reason why we use not only $\vartheta_i'$ but also $\vartheta_i$ in the faithful embeddings.
\end{Remark}

%\smallskip
Regarding Theorem~\ref{thm:intro:liftingtheorem}, 
we prove, more strongly, the surjectivity of the tropicalization of theta functions.

\begin{Theorem}[see Theorem~\ref{thm:surjective:tropicalization:theta} for the precise statement]
\label{thm:intro:lift}
The map \eqref{eqn:map:nonarch:to:trop:theta} is surjective. 
\end{Theorem}

To put it loosely, Theorem~\ref{thm:intro:lift} asserts that 
{\em any tropical $\Gamma$-rational theta function lifts to 
a nonarchimedean theta function}. We note that lifting tropical objects to algebro-geometric objects has been a central theme in tropical geometry. In many cases, 
such a lifting is not possible. In this sense, Theorem~\ref{thm:intro:lift} exhibits a good property of abelian varieties. 

\subsection{Related works}
While the names ``tropical abelian variety'' and ``tropical theta function'' seem to have been used first in Mikhalkin--Zharkov \cite{MZ}, tropical abelian varieties have been implicit in the literature for years such as Mumford \cite[\S6]{Mu}, Faltings--Chai \cite[Chap.~VI, \S1]{FC}, and Alexeev--Nakamura \cite{AN}. As is written in \cite{MZ}, tropical theta functions are implicitly suggested in \cite{AN}. Mikhalkin--Zharkov \cite{MZ} considers 
tropical theta functions on tropical Jacobians of tropical curves. 
Subsequently,  the theory of tropical abelian varieties has been further developed in Foster--Rabinoff--Shokrieh--Soto \cite{FRSS}, Sumi \cite{Sumi}, and Gross--Shokrieh \cite{GS}. 

There is now a large body of literature on faithful tropicalization. To the authors' knowledge, the faithful tropicalization problem dates back to Mikhalkin's ICM talk \cite{Mikhalkin}, in which he considered a particular case of faithful tropicalizations for elliptic curves. Then faithful tropicalizations for elliptic curves 
are studied in detail in Katz--Markwig--Markwig \cite{KMM, KMM2}. 
Baker--Payne--Rabinoff \cite{BPR} studies skeletons and tropicalizations for curves in general, and they show that for a given skeleton of a curve, there exist rational functions that faithfully tropicalize the skeleton. Then Gubler--Rabinoff--Werner \cite{GRW} generalizes the result in \cite{BPR} to the higher dimensional case. In those works on faithful tropicalizations, they use rational functions and not global sections of a given line bundle, thus their arguments do not address whether or not a given line bundle admits a faithful tropicalization. The first work that discussed faithful tropicalizations by global sections of line bundles is \cite{KY2}, in which we give a concrete number $C(g)$ for any genus $g$ curve, as is mentioned before Theorem~\ref{thm:intro:main:FT}. In \cite{KY1}, the authors discuss faithful tropicalizations in higher dimensions by adjoint line bundles. For more literatures on faithful tropicalizations, see, for example, the references in \cite{KY2}.

For the faithful embedding problem, when $X$ is a tropical curve, it is known that $L$ admits a faithful embedding once $\deg (L) \geq 2g+1$; indeed, this follows from the combination of \cite[Lemma~41 and Theorem~45]{HMY} and \cite[Theorem~3.3.15]{Song}. This condition on the degree is the same as the classical result for a line bundle to be very ample. In \cite{KY3}, the authors discuss faithful embeddings of tropical toric varieties, where the condition of admitting a faithful embedding is {\em not} the same as that of being very ample. 

\subsection{Organization of this paper}
This paper consists of ten sections including this section. Sections from \ref{section:TVandTL} to \ref{section:Pf:unimodular} are devoted to faithful embeddings of tropical abelian varieties, and Sections from \ref{section:uniformization}  to \ref{section:surjective:trop:theta} are devoted to faithful tropicalizations of the canonical skeletons of abelian varieties. In Section~\ref{section:TVandTL}, we recall the basic definitions in tropical geometry. In Section~\ref{section:trop:theta}, we recall tropical theta functions and the tropical Appell--Humbert theory for tropical line bundles on real tori. In Section~\ref{section:FE:trop:abel}, we state the main results on faithful embedding of tropical abelian varieties. In Section~\ref{section:injectivity}, we prove the injectivity part of the faithful embedding, and in Section~\ref{section:Pf:unimodular}, we prove the unimodularity part. In Section~\ref{section:uniformization}, we recall the Raynaud extensions of abelian varieties and the nonarchimedean Appell--Humbert theory by Bosch--L\"utkebohmert. In Section~\ref{sec:FT:skeleton}, we state the main results on faithful tropicalization of canonical skeletons of abelian varieties, and in Section~\ref{section:construction:NA:theta}, we prove those results.
In Section~\ref{section:surjective:trop:theta}, we prove the surjectivity of the tropicalization map for theta functions. 

\bigskip\noindent
{\sl Acknowledgment.}\quad
The first named author partially supported by KAKENHI 23K03041 and 20H00111,  and the second named author partially supported by KAKENHI 23K03046.

\setcounter{equation}{0}
\section{Faithful embeddings of rational polyhedral spaces} \label{section:TVandTL}

In this section, we define the notion of faithful embeddings.
Before that, we briefly review a general theory of tropical line bundles on rational polyhedral spaces.  

\subsection{Polyhedral spaces}
In this section, we recall the notion of rational polyhedral spaces that we will need. There are various incarnations of rational polyhedral spaces such as the tropical varieties from \cite{Mikhalkin}. Here we follow \cite[Subsection~2.A]{JRS}. 

A {\em polyhedron} (resp. \emph{rational polyhedron}) $\sigma$ in $\RR^n$ is a subset defined by a finite system of inequalities $\langle x, v\rangle \geq b$ in $\RR^n$, where $v \in \RR^n$ (resp. \emph{$v \in \ZZ$}) and $b \in \RR$. 
A {\em face} of $\sigma$ is $\sigma$ itself or a subset defined by changing some of the defining inequalities into equalities. Note that a face of a (rational) polyhedron is again a (rational) polyhedron. A {\em polyhedral complex} in $\RR^n$ is a locally finite collection $\mathcal{C}$ of polyhedra in $\RR^n$ with the following properties: 
\begin{enumerate}
\item[(i)]
If $\sigma$ belongs to $\mathcal{C}$, then any face of $\sigma$ belongs to $\mathcal{C}$; 
\item[(ii)]
If $\sigma$ and $\sigma'$ belong to $\mathcal{C}$ and $\sigma\cap \sigma' \neq \emptyset$ , then 
$\sigma\cap\sigma'$ is a face of both $\sigma$ and $\sigma'$. 
\end{enumerate}
A polyhedral complex $\mathcal{C}$ is said to be \emph{rational} if any $\sigma \in \mathcal{C}$ is a rational polyhedron.
A polyhedral complex is said to be \emph{finite} if it is a finite set. 
The underlying set of a polyhedral complex $\mathcal{C}$ is defined to be $|\mathcal{C}|  \colonequals \bigcup_{\sigma \in \mathcal{C}} \sigma \subseteq \RR^n$. 

A subset $X$ of $\RR^n$ is called a \emph{rational polyhedral set} if it is the underlying set of some rational polyhedral complex $\mathcal{C}$; then we call $\mathcal{C}$ a \emph{rational polyhedral decomposition} of $X$.
An {\em open rational polyhedral set} in $\RR^n$ is an open subset $U$ of $|\mathcal{C}|$ for some rational polyhedral complex $\mathcal{C}$ in~$\RR^n$. Here, we always endow $\RR^n$ with the Euclidean topology. 

A map $F\colon \RR^n \to \RR^m$ is said to be {\em $\ZZ$-affine} if there exist 
an $m \times n$ matrix $A$ with entries in $\ZZ$ and a $b \in \RR^m$ such that $F(x) = A x + b$ for all $x \in \RR^n$, where $x$ is regarded as a column vector. A map $\varphi\colon U \to V$, where $U \subseteq \RR^n$ and $V \subseteq \RR^m$ are open rational polyhedral sets, is called a {\em tropical morphism} if for any $p \in U$, there are an open neighborhood $U_p$ of $p$ in $U$ and a $\ZZ$-affine map $F_p\colon \RR^n \to \RR^m$ such that $\rest{\varphi}{U_p} = \rest{F_p}{U_p}$. Note that the composite of two tropical morphisms is a tropical morphism. A bijective tropical morphism whose inverse is a tropical morphism is called an \emph{isomorphism}.

\begin{Definition}[Rational polyhedral space]
\label{def:polyhedral:space}
Let $X$ be a Hausdorff topological space. A {\em rational polyhedral atlas} for $X$ is a collection of 
tuples $\{(X_i, \alpha_i, U_i)\}_{i \in \mathcal{I}}$ with the following properties: 
\begin{enumerate}
\item[(i)]
$\{X_i\}_{i \in \mathcal{I}}$ is an open covering of $X$; 
\item[(ii)]
for each $i \in \mathcal{I}$, $U_i$ is an open rational polyhedral
set in $\RR^{n_i}$ for some $n_i$, and $\alpha_i\colon X_i \to U_i$ is a homeomorphism; 
\item[(iii)]
for each pair $i, j \in \mathcal{I}$ with $X_i \cap X_j \neq \emptyset$, the transition map $\rest{\alpha_i \circ \alpha_j^{-1}}{\alpha_j(X_i \cap X_j)}: \alpha_j(X_i \cap X_j) \to \alpha_i(X_i \cap X_j)$ 
is a tropical morphism. 
\end{enumerate}
Two rational polyhedral atlases on $X$ are {\em equivalent} if their union forms a rational polyhedral atlas. 
A {\em rational polyhedral space} is a Hausdorff topological space $X$ together with 
the choice of an equivalence class of rational atlases on $X$. 
\end{Definition}

\begin{Remark}
A rational polyhedral space in Definition~\ref{def:polyhedral:space} is a special type of a rational polyhedral space defined in \cite[Definition~2.1]{JRS}. There, they allow the rational polyhedral spaces to have infinity boundaries. In this sense, our rational polyhedral spaces should be called ``boundaryless'' rational polyhedral space, but we use this convention because we do not need to mention rational polyhedral spaces with boundaries.
\end{Remark}

\subsection{Tropical line bundles and regular global sections}
\label{subsection:tropical:line:bundle}

Let $\TT = \RR \cup\{+ \infty\}$ denote the tropical semifield, as in the introduction. Note that in our convention, the tropical addition is $\min$ and the tropical multiplication is $+$. We endow $\TT$ with the Euclidean topology.

\begin{Definition}[tropical line bundle]
Let $X$ be a rational polyhedral space. 
We consider a pair $(L, \pi)$ of a topological space $L$ and a continuous map $\pi\colon L \to X$. Assume that there exist a rational polyhedral atlas $\{(X_i, \alpha_i, U_i)\}_{i \in \mathcal{I}}$ of $X$ and homeomorphisms $\{\beta_i\colon \pi^{-1}(X_i) \to X_i \times \TT\}_{i \in \mathcal{I}}$ such that: 
\begin{enumerate}
\item[(i)]
the diagram
\[
\begin{tikzcd}
 \pi^{-1}(X_i)  \arrow[dr, "\pi"] \arrow[rr, "\beta_i"] &   &   X_i\times \TT \arrow[dl, "{\rm pr}_1"] \\
& X_i & 
\end{tikzcd}
\]
commutes, where ${\rm pr}_1$ denotes the projection to the first factor; 
\item[(ii)]
for each pair $i, j \in \mathcal{I}$ with $X_i \cap X_j \neq \emptyset$, 
there is a tropical morphism $t_{ij}\colon \alpha_j(X_i \cap X_j) \to \RR$, called the \emph{transition function}, 
such that  
\[
\rest{(\alpha_i \times {\rm id_{\TT}}) \circ \beta_i \circ \beta_j^{-1} \circ (\alpha_j \times {\rm id_{\TT}})^{-1}}{\alpha_j(X_i \cap X_j) \times \TT}: \alpha_j(X_i \cap X_j) \times \TT 
\to \alpha_i(X_i \cap X_j) \times \TT
\]
is given by $(\alpha_j(x), a) \mapsto (\alpha_i(x), a + t_{ij}(x))$ for any $x \in X_i \cap X_j$ and $a \in \TT$. 
\end{enumerate}
We call such $\{(X_i, \alpha_i, U_i), \beta_i\}_{i \in \mathcal{I}}$ an \emph{atlas} for $(L, \pi)$. Two atlases $\{(X_i, \alpha_i, U_i), \beta_i\}_{i \in \mathcal{I}}$ and $\{(X_j^\prime, \alpha_j^\prime, U_j^\prime), \beta_j^\prime\}_{j \in \mathcal{J}}$ for $(L, \pi)$ are said to be {\em equivalent} if their union is again an atlas. 
A {\em tropical line bundle} over $X$ is a pair $(L, \pi)$ together with an equivalence class of atlases $\{(X_i, \alpha_i, U_i), \beta_i\}_{i \in \mathcal{I}}$ for $(L,\pi)$. We simply write $L$ for this.
\end{Definition}

Let $L^{(1)} = (L^{(1)}, \pi^{(1)})$ and $L^{(2)} = (L^{(2)}, \pi^{(2)})$ be tropical line bundles over $X$. We call a homeomorphism $f\colon L^{(1)} \to L^{(2)}$ satisfying $\pi^{(1)} =\pi^{(2)} \circ f$  an {\em isomorphism of tropical line bundles} if there exist atlases $\{(X_i, \alpha_i, U_i), \beta_i^{(1)}\}_{i \in \mathcal{I}}$ and $\{(X_i, \alpha_i, U_i), \beta_i^{(2)}\}_{i \in \mathcal{I}}$ of $L^{(1)}$ and $L^{(2)}$, respectively,  such that 
for each $i \in \mathcal{I}$, 
there is a tropical morphism $t_i\colon U_i \to \RR$ 
such that  
\[
(\alpha_i \times {\rm id_{\TT}}) \circ \beta_i \circ f \circ \beta_i^{-1} \circ (\alpha_i \times {\rm id_{\TT}})^{-1}: U_i \times \TT 
\to U_i \times \TT
\]
is given by $(\alpha_i(x), a) \mapsto (\alpha_i(x), a + t_{i}(x))$ for any $x \in X_i \cap X_j$ and $a \in \RR$. 

We can naturally define the pullback of a tropical line bundle by a tropical morphism. Indeed, the pullback as a topological space  is defined by fiber product, and since the composite of the tropical morphism is a tropical morphism, the pullbacks of the transition functions for $L$ become transition functions over $Y$.

Let $\Pic (X)$ be the set of isomorphism classes of tropical line bundles over $X$. Then $\Pic (X)$ has a natural structure of abelian group. Indeed, we take two tropical line bundles $L$ and $L' \in \Pic (X)$. We may take the same atlas $\{(X_i, \alpha_i, U_i)\}_{i \in \mathcal{I}}$ for $L$ and $L^\prime$, and let $\{ t_{ij} \}$ and $\{ t_{ij}' \}$ denote the family of transition functions of 
$L$ and $L^\prime$ for this atlas, respectively. 
Then we can construct a tropical line bundle whose transition functions are given by $\{ t_{ij} + t_{ij}'\}$ that for the atlas $\{(X_i, \alpha_i, U_i)\}_{i \in \mathcal{I}}$. Further, such a tropical line bundle is unique up to isomorphism, which we write for $L + L'$. We always regard $\Pic (X)$ as an abelian group with this group structure and call it the \emph{tropical Picard group of $X$}. A tropical morphism between rational polyhedral spaces induces a group homomorphism between the tropical Picard groups by pullback.

We are going to define {\em regular global sections} of a tropical line bundle. 
A map $F\colon \RR^n \to \RR^m$ is said to be {\em piecewise $\ZZ$-affine} if there exists a rational polyhedral decomposition $\mathcal{C}^\prime$ of $\RR^n$ such that for each $\sigma^\prime \in \mathcal{C}^\prime$, there exists a $\ZZ$-affine map $F_{\sigma^\prime}\colon \RR^n \to \RR^m$ such that $\rest{F}{{\sigma^\prime}} = \rest{F_{\sigma^\prime}}{\sigma^\prime}$. Note that then $F$ is continuous. 
Let $U \subseteq \RR^n$ be an open rational polyhedral set. A map $\varphi\colon U \to \RR^m$ is said to be {\em 
piecewise $\ZZ$-affine} if for any $p \in U$, there are an open neighborhood $U_p$ of $p$ in $U$ 
and a piecewise $\ZZ$-affine map $F_p\colon \RR^n \to \RR^m$ such that $\rest{\varphi}{U_p} = \rest{F_p}{U_p}$. Note that for a map $F : \RR^n \to \RR^m$, there are two definitions of piecewise $\ZZ$-affine maps, but they coincide with each other.
When $m=1$, we call such an $F$ a piecewise $\ZZ$-affine {\em function}. 
 
Let $V \subseteq \RR^n$ be a convex set. A function $f\colon V \to \RR$ is {\em concave} if for any $p_1, p_2 \in V$ and $t \in [0, 1]$, one has 
$f((1-t) p_1 + t p_2) \geq t f(p_1) + (1-t) f(p_2)$. 
A function $\varphi\colon U \to \RR$, where $U \subseteq \RR^n$ is an open rational polyhedral set, is said to be {\em concave} 
if for any $p \in U$, there are an open neighborhood $U_p$ of $p$ in $U$, 
a convex subset $V$ of $\RR^n$ containing $U_p$, and a function $F_p\colon \RR^n \to \RR$ such that $\rest{\varphi}{U_p} = \rest{F_p}{U_p}$ and $\rest{F_p}{V}\colon V \to \RR$ is concave. 

A \emph{global section} of a tropical line bundle $\pi : L \to X$ is a map $s : X \to L$ such that $\pi \circ s = \id_X$.

\begin{Definition}[Regular section of a tropical line bundle]
Let $L$ be a tropical line bundle on a rational polyhedral space $X$. 
We take an atlas $\left((L, \pi), \{(X_i, \alpha_i, U_i), \beta_i\}_{i \in \mathcal{I}}\right)$ of $L$. 
A global section~$s$ of the tropical line bundle $L$ is said to be \emph{regular} if one of the following holds;
\begin{enumerate}
\item[(i)]
for each $i \in \mathcal{I}$, 
\begin{equation}
\label{eqn:def:global:section}
 {\rm pr}_2 \circ \beta_i \circ s \circ \alpha_i^{-1}\colon U_i \to \TT
\end{equation}
has its image in $\{ + \infty \} \subseteq \TT$; 
\item[(ii)]
for each $i \in \mathcal{I}$, the map in (\ref{eqn:def:global:section}) has image in $\RR \subseteq \TT$ and it is a concave locally piecewise $\ZZ$-affine map.
\end{enumerate}
There exists a unique regular global section as in (i). We call it the \emph{trivial} section, denoted by $+ \infty$. It is easy to see that the notion of regular global sections does not depend on the choice of an atlas of $L$. A regular global section of the trivial tropical line bundle is called a \emph{regular function}. 
\end{Definition}

\begin{Remark} \label{remark:regular-locallypiecewiseaffine}
Let $U$ be an open rational polyhedral set and let $L = U \times \TT$ be the trivial tropical line bundle on $U$. Let $s$ be a nontrivial section of $L$. Then it is regular if and only if it is an $\RR$-valued concave piecewise $\ZZ$-affine function. 
\end{Remark}

We write $H^0 (X ,L)$ for the set of regular global sections of $L$. Then $H^0 (X ,L)$ has a natural structure of $\TT$-semimodule. Indeed, for $s_1, s_2 \in H^0(X, L)$, one defines 
$s_1 \oplus s_2 \colonequals \min\{s_1, s_2\} \in H^0(X, L)$ by taking the minimum for the maps (induced by $s_1$ and $s_2$) in \eqref{eqn:def:global:section} for each $i \in \mathcal{I}$. Further, one defines the action $\TT \times H^0 (X ,L) \to H^0 (X ,L)$ by $(a , s) \mapsto a + s$. Then these makes $H^0 (X,L)$ a $\TT$-semimodule.

\subsection{Unimodularity and faithful embeddings}
In this subsection, we define the notion concerning faithful embeddings. 
First, we recall the notion of unimodularity.
A $\ZZ$-affine map $F(x) = Ax + b$ with  an $m \times n$ matrix $A$ with entries in $\ZZ$ and $b \in \RR^n$ is said to be {\em unimodular} if $A$ is unimodular, i.e., $\rank (A) = n$ and $\ZZ^m/A\ZZ^n$ is torsion-free. Let  $U$ be an open rational polyhedral set in $\RR^n$. A map $\varphi\colon U \to \RR^m$ is said to be {\em unimodular} if there is a rational polyhedral complex $\mathcal{C}$ in $\RR^n$ such that $U$ is an open subset of $|\mathcal{C}|$ with the property that for each $\sigma \in \mathcal{C}$, there is a unimodular $\ZZ$-affine map $F\colon \RR^n \to \RR^m$ such that $\rest{\varphi}{\sigma \cap U} = \rest{F}{\sigma \cap U}$. 
Let $X$ be a rational polyhedral space. We say that a map $\varphi : X \to \RR^m$ is \emph{unimodular} if for some atlas $\{ (X_i , \alpha_i , U_i) \}_i$ of $X$ and for any $i \in \mathcal{I}$, the map $\varphi \circ \alpha_i^{-1} : U_i \to \RR^m$ is unimodular. Note that then the same holds for any atlas of $X$.

The $m$-dimensional {\em tropical projective space} is defined to be 
\[
  \TT\PP^m
  \colonequals 
  \left(\TT^{m+1}\setminus\{(+\infty, \ldots, +\infty)\}\right)/\!\sim,  
\]
where $(a_0, \ldots, a_m) \sim (b_0, \ldots, b_m)$ if 
there exists a $c \in \RR$ such that $b_0 = c + a_0, \ldots, b_m = c + a_m$. 
The subspace $\RR^{m+1} /\!\sim = \{ (a_0 : \cdots : a_m) \in \TT\PP^m \mid a_0 , \ldots , a_m \in \RR \}$ has a canonical structure of rational polyhedral space such that the homomorphism $\RR^m \to \TT\PP^m$ given by $(a_1 , \ldots , a_m) \mapsto (0 : a_1 : \cdots : a_m)$ is an isomorphism of rational polyhedral spaces. Note that via this isomorphism, $\RR^{m+1} /\!\sim$ has a $\ZZ$-structure.
We call $\RR^{m+1} /\!\sim$ the \emph{maximal tropical torus} of $\TT\PP^m$.
Note that for any $i =0, \ldots , m$, the map $\RR^m \to \RR^{m+1} /\!\sim$ given by $(a_1 , \ldots , a_m) \mapsto ( a_1 : \cdots : a_{i} : 0 : a_{i+1} : \cdots : a_m)$ is also an isomorphism of rational polyhedral spaces. 

\begin{Definition}[unimodularity, faithful embedding] \label{def:unimodular:FE}
Let $X$ be a rational polyhedral space. Let $\varphi : X \to \TT\PP^m$ be a map. Assume that $\varphi (X)$ is contained in the maximal tropical torus of $\TT\PP^m$. We say that $\varphi$ is \emph{unimodular} if it is unimodular as a map from $X$ to $\RR^m$. We call $\varphi$ a \emph{faithful embedding} if it is unimodular and a homeomorphism onto its image.
\end{Definition}

We remark that we can define the notions of unimodularity and faithful embeddings even when $\varphi (X)$ is not contained in the maximal tropical torus of $\TT\PP^m$, but we omit those definitions because we will not need them in the sequel.

Let $L$ be a tropical line bundle on $X$ and let $s_0, \ldots, s_m \in H^0(X, L) \setminus \{+\infty\}$. Then we define a map $\varphi_{s_0 , \ldots , s_m} : X \to \TT\PP^m$ by $x \mapsto (s_0 (x) : \cdots : s_m (x))$, and the image is in the maximal tropical torus.

\begin{Definition} \label{def:FE}
Let $L$ be a tropical line bundle on a rational polyhedral space $X$. We say that $s_0 , \ldots , s_m \in H^0 (X ,L) \setminus \{ + \infty \}$ give a \emph{unimodular} map (resp. a \emph{faithful embedding}) for $X$ if the map $\varphi_{s_0 , \ldots , s_m} : X \to \TT\PP^m$ is a unimodular map (resp. a faithful embedding). 
%We say that $s_0 , \ldots , s_m \in H^0 (X ,L) \setminus \{ + \infty \}$ \emph{give a faithful embedding} of $X$ if $\varphi_{s_0 , \ldots , s_m} : X \to \TT\PP^m$ is a faithful embedding.
We say that $L$ \emph{admits a faithful embedding} if there exist $s_0 , \ldots , s_m \in H^0 (X ,L) \setminus \{ + \infty \}$ that give a faithful embedding for $X$.
\end{Definition}

%\begin{Remark} \label{remark:FE:uptoconstant}
%With the notation in Definition~\ref, if $s_0 , \ldots , s_m \in H^0 (X,L) \setminus \{ + \infty \}$ give a faithful embedding of $X$, then for any $c_0 , \ldots , c_m \in \RR$, so do $s_0 + c_0 , \ldots , s_m + c_m$.
%\end{Remark}

\begin{Remark} \label{remark:unimdular:part}
With the notation in Definition~\ref{def:FE}, suppose that a part of $s_0 , \ldots , s_m \in H^0 (X ,L) \setminus \{ + \infty \}$ give a unimodular map (resp. faithful embedding). Then the following hold: 
\begin{enumerate}
\item
$s_0 + c_0 , \ldots , s_m + c_m$ give a unimodular map (resp. faithful embedding);
\item
$s_0 , \ldots , s_m$ give a unimodular map (resp. faithful embedding). 
\end{enumerate}
\end{Remark}

\setcounter{equation}{0}
\section{Tropical theta functions and tropical line bundles on real tori}
\label{section:trop:theta}

In this section, we recall known facts on tropical theta functions and tropical line bundles on real tori with integral structures. The basic references in this section are \cite{MZ, FRSS, GS, Sumi}.

\subsection{Notation} \label{subsection:notation1}

We fix the notation. Let $M$ be a free $\ZZ$-module of finite rank $n$ and set $N := \mathrm{Hom} (M , \ZZ)$. We have $N_{\RR} := N \otimes \RR = \mathrm{Hom} (M , \RR)$. Let $\langle \ndot, \ndot \rangle : M \times N_{\RR} \to \RR$ be the canonical pairing.
Let $M'$ be another free $\ZZ$-module of rank $n$. Assume that we are given an embedding $M' \hookrightarrow N_{\RR}$ of lattice, i.e., an injective homomorphism that induces an isomorphism $M'_{\RR} \to N_{\RR}$. We usually regard $M'$ as a lattice in $N_{\RR}$ via this embedding.

\subsection{Tropical theta functions}
First, we recall the definition of tropical theta functions. Our basic 
references are \cite{MZ, FRSS, Sumi}. 

To explain the tropical theta functions, we basically follow the convention in \cite{FRSS}, because this is useful when we consider lifting 
of tropical theta functions to nonarchimedean theta functions in Sections~\ref{section:construction:NA:theta} and \ref{section:surjective:trop:theta}.  
We consider a pair $(\lambda , \gamma )$ consisting of a homomorphism $\lambda : M' \to M$ and a map $\gamma : M' \to \RR$ such that for any $u_1' , u_2' \in M^\prime$,
\begin{align} \label{align:tropicaldescentcondition}
\gamma (u_1' + u_2') - \gamma (u_1') - \gamma (u_2') = \langle \lambda (u_2') , u_1' \rangle
.
\end{align}
We call such a pair $(\lambda , \gamma)$ a \emph{tropical descent datum}. Note that since the left-hand side on equality (\ref{align:tropicaldescentcondition}) is symmetric on $(u_1' , u_2')$, the bilinear form $(u_1' ,u_2') \mapsto \langle \lambda (u_2') , u_1' \rangle$ is symmetric.

Let $(\lambda , \gamma)$ be a tropical descent datum. A function $\vartheta : N_{\RR} \to \TT$ is said to be \emph{quasi-periodic} with respect to $(\lambda , \gamma)$ if for any $u' \in M'$ and $x \in N_{\RR}$, one has
\[
\vartheta (x+u') = \vartheta  (x) - \langle \lambda (u') , x \rangle - \gamma (u').
\]
A \emph{tropical theta function} with respect to $(\lambda , \gamma)$ is a regular function on $N_{\RR}$ (cf. Remark~\ref{remark:regular-locallypiecewiseaffine}) that is quasi-periodic with respect to $(\lambda,\gamma)$. 

There is another convention to characterize the quasi-periodicity, which is given in terms of a symmetric bilinear form on $N_{\RR}$ and a homomorphism $M' \to \RR$.
We say that a symmetric $\RR$-bilinear form $Q\colon N_\RR \times N_\RR \to \RR$ is \emph{integral} over $M' \times N$, or integral, simply, if $Q(M^\prime \times N) \subseteq \ZZ$. The following lemma gives a correspondence between an integral symmetric 
$\RR$-bilinear form   $Q\colon N_\RR \times N_\RR \to \RR$ and a group homomorphisms $\lambda : M' \to M$.

\begin{Lemma}
\label{lem:polarization:equiv}
Let $M$, $N$, $N_{\RR}$, and $M'$ be as above.
Then the following hold.
\begin{enumerate}
\item
There is a natural bijection between 
\begin{enumerate}
\item[(i)]
the set of symmetric $\RR$-bilinear forms $Q\colon N_\RR \times N_\RR \to \RR$ with $Q(M^\prime \times N) \subseteq \ZZ$; and 
\item[(ii)]
the set of $\ZZ$-linear maps $\lambda\colon M^\prime \to M$ such that 
$\langle \lambda (u_1^\prime), u_2^\prime \rangle = \langle \lambda (u_2^\prime), u_1^\prime \rangle$ for any 
$u_1^\prime, u_2^\prime \in M^\prime$. 
\end{enumerate}
Further, this bijection is given by 
\begin{align} \label{align:Q-lambda}
\rest{Q}{M^\prime \times M^\prime} = \langle \lambda (\cdot), \cdot \rangle.
\end{align} 
\item
Let $\lambda : M' \to M$ be a homomorphism and let $Q$ be an integral symmetric bilinear form on $N_{\RR}$ corresponding to $\lambda$ by (1) above. Let $\gamma : M' \to \RR$ and $\ell : M' \to \RR$ be maps. Assume that 
\begin{align} \label{align:gamma-ell}
\gamma (u') = \frac{1}{2} Q(u',u') - \ell (u')
.
\end{align}
Then $\gamma$ satisfies (\ref{align:tropicaldescentcondition}) if and only if $\ell$ is a homomorphism. Further, when $\gamma$ satisfies (\ref{align:tropicaldescentcondition}) (and hence $\ell$ is a homomorphism), a function $\vartheta : N_{\RR} \to \RR$ is quasi-periodic if and only if
\begin{equation}
\label{eqn:quasi:periodicity:0}
\vartheta (x +  u^\prime)
= \vartheta (x) -  Q(x,  u^\prime) - \frac{1}{2} Q(u^\prime, u^\prime)
  + \ell(u^\prime)
\end{equation}
for all $x \in N_{\RR}$ and $u^\prime \in M^\prime$.
\end{enumerate}
\end{Lemma}

\Proof
We prove (1). Suppose that we are given a symmetric $\RR$-bilinear form $Q\colon N_\RR \times N_\RR \to \RR$ with $Q(M^\prime \times N) \subseteq \ZZ$. 
For any $u^\prime \in M^\prime \subseteq N_{\RR}$,  
noting that $Q(M^\prime \times N) \subseteq \ZZ$, 
we define $\lambda(u^\prime) := Q(u^\prime, \cdot) \in \Hom_\ZZ(N, \ZZ) = M$. Then 
$\lambda\colon M^\prime \to M$ is $\ZZ$-linear, and $\langle \lambda(u_1^\prime), u_2^\prime \rangle = Q(u_1^\prime, u_2^\prime)$ for any 
$u_1^\prime, u_2^\prime \in M^\prime$. 
Since $Q$ is symmetric, we have 
$\langle \lambda (u_1^\prime), u_2^\prime \rangle = \langle \lambda(u_2^\prime), u_1^\prime \rangle$. 

Next, suppose that we are given a $\ZZ$-linear map $\langle \lambda (u_1^\prime), u_2^\prime \rangle = \langle \lambda(u_2^\prime), u_1^\prime \rangle$. 
Then the map $\langle \lambda (\cdot) , \cdot \rangle \colon M^\prime \times M^\prime \to \RR$ is a $\ZZ$-bilinear map. Since $M^\prime_\RR = N_\RR$, we extend $\langle \lambda (\cdot) , \cdot \rangle$ to $Q\colon N_\RR \times N_\RR \to \RR$ by linearity. Then $Q$ is symmetric. 
Further, for any $u^\prime \in M^\prime$ and $v \in N = \Hom_{\ZZ}(M, \ZZ)$. we have $Q(v, u^\prime) = \langle \lambda(u^\prime),v \rangle \in \ZZ$, and thus $Q(N \times M^\prime) \subseteq \ZZ$. By symmetry  $Q(M^\prime \times N) \subseteq \ZZ$. This proves (1)

Assertion (2) is shown by a straightforward computation using (\ref{align:tropicaldescentcondition}).
\QED

Let $(\lambda , \gamma)$ be a tropical descent datum. Then by Lemma~\ref{lem:polarization:equiv}, we have a unique pair $(Q,\ell)$ of a symmetric $\RR$-bilinear form  $Q : N_{\RR} \times N_{\RR} \to \RR$ that is integral over $M' \times N$ given by (\ref{align:Q-lambda}) and a homomorphism $\ell : M' \to \RR$ given by (\ref{align:gamma-ell}), which we call the the pair corresponding to $(\lambda , \gamma)$. 

\subsection{Tropical Appell--Humbert theory}

In this subsection, we recall the tropical Appell--Humbert theory. This theory works over the real tori with integral structures, which is defined below. 

\begin{Definition}[real torus with an integral structure]
\label{def:real:torus}
Let $M$, $N_{\RR}$, $M'$ be as in Subsection~\ref{subsection:notation1}
The quotient $X = N_\RR/M^\prime$ is called a {\em real torus with an integral structure}. Here, the integral structure means the lattice $N$ in $N_\RR$. 
\end{Definition}

In our convention, \emph{a real torus always means a real torus with an integral structure, and we often omit the words ``with an integral structure.''}

The topological space $N_{\RR}$ has a canonical structure of polyhedral space. Since the canonical projection $\pi : N_{\RR} \to X$ is a local isomorphism, we endow $X$ with a structure of rational polyhedral space from $N_{\RR}$ by $\pi$. Then $\pi$ is a tropical morphism, obviously.

We recall the tropical Appell--Humbert theory. Let $X = N_\RR/M^\prime$ be a real torus. 
Let $(\lambda , \gamma)$ be a tropical descent datum.
We define an action of $M^\prime$ on $N_\RR \times \TT$ by 
\[
  u^\prime \cdot (x, t) 
  \colonequals 
  \left(x + u^\prime, \; t - \langle  \lambda (u^\prime) , x \rangle 
 - \gamma (u') \right). 
\]
Then the quotient $L(\lambda, \gamma) \colonequals (N_\RR \times \TT)/M^\prime$ is a tropical line bundle on $X = N_\RR/M^\prime$. 

\begin{Remark} \label{remark:translateQell}
Let $(Q, \ell)$ be a pair of a symmetric bilinear $\RR$-from $Q : N_{\RR} \times N_{\RR} \to \RR$ such that $Q(M' \times N) \subseteq \ZZ$ and a homomorphism $\ell : M' \to \RR$ corresponding to the tropical descent datum $(\lambda , \gamma)$ by Lemma~\ref{lem:polarization:equiv}.
The above action of $M'$ on $N_{\RR} \times \TT$ is described as
\[
  u^\prime \cdot (x, t) 
  \colonequals 
  \left(x + u^\prime, \; t - Q(u^\prime, x) - \frac{1}{2} Q(u^\prime, u^\prime)
  + \ell(u^\prime)\right). 
\]
\end{Remark}

For a tropical line bundle $L$ on $X$, we call the isomorphism between the fiber of $L \to X$ over $0 \in X$ and $\TT$ a \emph{rigidification}, and we say that a tropical line bundle is \emph{rigidified} if a rigidification is specified. Any tropical line bundle on $X$ has a rigidification. Further, for any tropical descent datum $(\lambda , \gamma)$, $L(\lambda, \gamma)$ has a canonical rigidification coming from the canonical rigidification of the trivial tropical line bundle on $N_{\RR}$ through the quotient construction. 

Note that an isomorphism between tropical line bundles on $X$ that respects given rigidifications is unique, because there are only constant functions that are regular on $X$. Thus for rigidified tropical line bundles $L_1$ and $L_2$ on $X$ that are isomorphic to each other, it makes sense to write $L_1 = L_2$. 

\begin{Remark} \label{remark:additive}
It is straightforward to see that $L(\lambda_1 , \gamma_1) + L(\lambda_2 , \gamma_2) = L(\lambda_1 + \lambda_2 , \gamma_1 + \gamma_2)$ as rigidified line bundles for any tropical descent data $(\lambda_1 , \gamma_1)$ and $(\lambda_2 , \gamma_2)$.
\end{Remark}

The following theorem shows that any (rigidified) tropical line bundle on a real torus is of form $L(\lambda, \gamma)$ for some tropical descent datum $(\lambda , \gamma)$.
Thus we need to consider tropical line bundles of this form only.

\begin{Theorem}[Tropical Appell--Humbert Theorem {\cite[Proposition~28]{Sumi}}, {\cite[Theorem~7.2]{GS}}]
\label{thm:trop:Appell:Humbert} 
Let $L$ be a rigidified tropical line bundle on a real torus $X = N_\RR/M^\prime$. Then there exists a tropical descent datum $(\lambda , \gamma)$ such that
$L = L(\lambda, \gamma)$. Further, $L (\lambda, \gamma) = L(\lambda^\prime, \gamma^\prime)$ if and only if $\lambda^\prime = \lambda$ and $(\gamma^\prime - \gamma)(N) \subseteq \ZZ$. 
\end{Theorem}

\Proof
Noting Remark~\ref{remark:translateQell}, this is nothing but \cite[Proposition~28]{Sumi} and \cite[Theorem~7.2]{GS}.
\QED

It follows from the construction of the rigidified line bundle $L(\lambda, \gamma)$ that the canonical projection $\pi : N_{\RR} \to X$ induces a natural map $\pi^{\ast} : H^{0} (X, L(\lambda ,\gamma)) \to H^0 (N_{\RR} , N_{\RR} \times \TT)$ by pullback, where $N_{\RR} \times \TT$ is a trivial tropical line bundle on $N_{\RR}$. In fact, one sees that $\pi^{\ast}$ induces an isomorphism from $H^{0} (X, L(\lambda ,\gamma))$ onto the $\TT$-subsemimodule consisting of the tropical theta functions with respect to $(\lambda , \gamma)$.

\begin{Remark} \label{remark:translate:u}
We take any $u \in M$ and naturally regard it as a function $N_{\RR} \to \RR$. Let $(\lambda, \gamma)$ be a tropical descent datum. Then $(\lambda , \gamma - \rest{u}{M'})$ is also a tropical descent datum. Further, for any tropical theta function $\vartheta$ with respect to $(\lambda, \gamma)$, $\vartheta + u $ is a tropical theta function with respect to $(\lambda, \gamma -  \rest{u}{M'} )$, and this translation by $u$ gives an isomorphism $H^0 (X , L(\lambda , \gamma)) \to H^0 (X, L (\lambda , \gamma -  \rest{u}{M'} ))$.
\end{Remark}

Keeping the above identification in mind, we have the following theorem due to Sumi \cite{Sumi}. Here, a map $f : X \to Y$ between real tori called a \emph{homomorphism} if it is a tropical morphism and is a homomorphism of groups.

\begin{Proposition} \label{prop:Globalsections:Sumi}
Let $X$ be as above. Let $(\lambda , \gamma)$ be a descent datum. Let 
$Q$ be the symmetric $\RR$-bilinear form corresponding to $\lambda$ by Lemma~\ref{lem:polarization:equiv}~(1). For the tropical line bundle $L(\lambda,\gamma)$, the following hold.
\begin{enumerate}
\item
If $Q$ is not positive-semidefinite, then $H^0 ( X , L(\lambda , \gamma)) = \{ + \infty \}$.
\item
Assume that $Q$ is positive-semidefinite and that $L(\lambda,\gamma)$ has a nontrivial regular global section. Then there exist a real torus $Y$ with an integral structure, a line bundle $L(\lambda' , \gamma')$ on $Y$, and a homomorphism $f : X \to Y$ such that the symmetric $\RR$-bilinear form $Q'$ corresponding to $\lambda'$ by Lemma~\ref{lem:polarization:equiv}~(1) is positive-definite and the induced map $f^{\ast} : H^0 (Y , L (\lambda' , \gamma')) \to H^0 (X , L (\lambda , \gamma))$ is an isomorphism of $\TT$-modules. Further, $Q$ is positive-definite if and only if $f$ is an isomorphism.
\end{enumerate}
\end{Proposition}

\Proof
(1) is due to \cite[Remark~29]{Sumi}, and (2) follows from the argument in \cite[Section~3.5]{Sumi}; see also in the alternative description of tropical theta functions in Subsection~\ref{subsec:Alternative}.
\QED

The above projection indicates that it is important to understand the tropical theta functions when $Q$ is positive-definite. In this case, we have concrete tropical theta functions, as is seen in the next subsection.

\subsection{Construction of tropical theta functions} \label{subsection:construction:tropicaltheta}

In this subsection, we recall concrete tropical theta functions, which will form a system of generators of the $\TT$-semimodule of tropical theta functions with respect to a give descent datum. This is a modification of the construction of Mikhalkin--Zharkov~ \cite{MZ} and Sumi~\cite{Sumi}; see (\ref{eqn:def:theta:b2}) and Remark~\ref{remark:construction:trop:theta} in the sequel.

We keep the notation in Subsection~\ref{subsection:notation1}, such as $M$, $N$, $N_{\RR}$, and $M'$. Set $X = N_{\RR} / M'$ and let $\pi : N_{\RR} \to X$ denote the canonical surjection.  Let $(\lambda , \gamma)$ be a tropical descent datum. Let $Q:N_{\RR} \times N_{\RR} \to \RR$ be the $\RR$-bilinear form corresponding to $\lambda$ by Lemma~\ref{lem:polarization:equiv}.  

\begin{Definition}[polarization, tropical abelian varieties]
\label{def:trop:av}
\begin{enumerate}
\item[(1)]
We call $\lambda$ a \emph{polarization} on $N_{\RR} / M'$ if $Q$ is positive-definite. By abuse of words, we sometimes call $Q$ a polarization. 
\item[(2)]
A real torus $X$ is called a {\em tropical abelian variety} if it has a polarization.
\end{enumerate}
\end{Definition}

The existence of a polarization in Definition~\ref{def:trop:av} is a tropical version of  Riemann's period relations for complex abelian varieties. 

\begin{Remark} \label{remark:polarization:injective}
If $\lambda : M' \to M$ is a polarization, then $Q$ is nondegenerate. It follows that the $\RR$-linear extension $\lambda_{\RR} : M'_{\RR} \to N_{\RR}$ is an isomorphism, and thus $\lambda$ is injective.
\end{Remark}

\begin{Lemma} \label{lemma:towardSumiTheta_new1}
We fix any $b \in M$. With the above notation, assume that $\lambda$ is a polarization. Then 
for any $x \in N_{\RR}$, the subset
$\left\{
\langle b + \lambda (u') , x \rangle + \gamma (u') + \langle b, u' \rangle \mid u' \in M'
\right\}
$
of $\RR$ has a minimum.
\end{Lemma}

\Proof
Fix any $x \in N_{\RR}$. For any $u' \in M'$, a straightforward computation shows us 
\begin{align*}
\langle b + \lambda (u') , x \rangle + \gamma (u') + \langle b ,u' \rangle 
=  \langle b + \lambda (u') , x \rangle  - \ell (u') + \langle b ,u' \rangle  +
\frac{1}{2} Q (u'  , u' ) 
.
\end{align*}
Note that the first three terms on the right-hand side are linear on $u'$. Since $Q(u',u')$ a positive-definite quadratic form on $u'$, it follows that for a fixed $x \in N_{\RR}$,
\[
\left\{\left.   \langle b + \lambda (u') , x \rangle  - \ell (u') + \langle b ,u' \rangle  +
\frac{1}{2} Q (u'  , u' )  \;\right|\; u' \in M'\right\}
\]
has a minimum. This proves the lemma.
\QED

Assume that $\lambda$ is a polarization on $N_{\RR}/M'$. We define a function $\vartheta_b^{(\lambda ,\gamma)} : N_{\RR} \to \RR$ by
\begin{align}
\label{eqn:def:theta:b0}
\vartheta_b^{(\lambda ,\gamma)} (x) :=
\min_{u' \in M'} (
\langle b + \lambda (u') , x \rangle + \gamma (u') + \langle b, u' \rangle
) ,
\end{align}
which is well-defined by Lemma~\ref{lemma:towardSumiTheta_new1}. This is a nontrivial regular function and is quasi-periodic. Thus $\vartheta_b^{(\lambda , \gamma)}$ is a tropical theta function with respect to $(\lambda , \gamma)$.

%%%%%%%%%%%%%
% Faithful embeddings %%
%%%%%%%%%%%%%a
\setcounter{equation}{0}
\section{Faithful embeddings of tropical abelian varieties}
\label{section:FE:trop:abel}

In this section, we state the main results on faithful embeddings and show some examples. We keep the notation in Subsection~\ref{subsection:notation1}, such as $M$, $N$, $N_{\RR}$, and $M'$. Set $X : = N_{\RR} / M'$, which is a real torus of dimension $n$, and let $\pi : N_{\RR} \to X$ denote the canonical surjection. 

\subsection{Polarization types}

Suppose that $X = N_{\RR} / M'$ is a tropical abelian variety, and 
let $\lambda : M' \to M$ be a polarization on $X$. By Remark~\ref{remark:polarization:injective}, $\lambda$ is injective. It follows that, 
by choosing suitable free $\ZZ$-bases of $M^\prime$ and $M$, there exists a unique $n$-tuple of positive integers $(d_1, \ldots, d_n)$ such that $d_i \mid d_{i+1}$ for $1 \leq i \leq n-1$ and $\lambda$ can be represented by the diagonal matrix $\mathrm{diag}(d_1, \ldots , d_n)$. We call $(d_1 , \ldots , d_n)$ the \emph{type} of $\lambda$ or the corresponding $\RR$-bilinear form $Q$. Note that $M/\lambda(M^\prime) \cong \ZZ/(d_1) \oplus \cdots \oplus \ZZ/(d_n)$. In particular, if $\mathfrak{B}$ is a complete system of representatives of $M/\lambda(M^\prime)$, then $|\mathfrak{B}| = d_1 \cdots d_n$.

\begin{Remark}
Let $\lambda : M' \to M$ be a polarization of type $(d_1 , \ldots , d_n)$ and let $Q$ be the corresponding bilinear $\RR$-form. Then, by a suitable choice of $\ZZ$-basis of $M'$ and $N$, $Q$ is represented by the diagonal matrix ${\rm diag}\,(d_1, \ldots, d_n)$.
\end{Remark}

\begin{Definition}[positive, polarization type]
Let $X = N_{\RR} /M'$ be a real torus and 
$L$ be a tropical line bundle. We take  tropical descent datum $(\lambda , \gamma)$ with $L \cong L(\lambda, \gamma)$. 
We say that $L$ is \emph{positive} if $\lambda$ is a polarization; by Theorem~\ref{thm:trop:Appell:Humbert}, this is well-defined. 
Further, we define the \emph{polarization type} of $L$ to be the type of $\lambda$, which is also well-defined by Theorem~\ref{thm:trop:Appell:Humbert}.
\end{Definition}

\subsection{Main results}

In the theory of complex tori, the Lefschetz theorem gives a sufficient condition in terms of the type for a line bundle to be very ample. Our main results will give an analogous sufficient condition for a tropical line bundle on a real torus with an integral structure to admit a faithful embedding. 

The following theorem is one of the main results in this paper. It is a tropical version of the classical Lefschetz theorem.

\begin{Theorem}
\label{thm:main:1-1}
Let $X = N_{\RR}/M'$ be a tropical abelian variety of dimension $n$ 
and let $\pi : N_{\RR} \to X$ denote the canonical surjection. Let $L$ be a positive tropical line bundle on $X$ of polarization type $(d_1 , \ldots , d_n)$, and we write $L = L(\lambda , \gamma)$ with a tropical descent datum $(\lambda , \gamma)$. Set $D := d_1 \cdots d_n$ and let $\{ b_1 , \ldots , b_D\}$ be a complete system of representatives of $M / \lambda (M')$. For each $i=1 , \ldots , D$, let $\bar{\vartheta}_b^{(\lambda , \gamma)} \in H^{0} (X,L)$ be a regular global section such that $\pi^{\ast} (\bar{\vartheta}_{b_i}^{(\lambda , \gamma)}) = \vartheta_{b_i}^{(\lambda ,\gamma)}$. Let $\varphi : X \to \TT\PP^{D-1}$ be the map defined by $x \mapsto \left(\bar{\vartheta}_{b_1}^{(\lambda , \gamma)} (x) : \cdots : \bar{\vartheta}_{b_D}^{(\lambda , \gamma)} (x) \right)$. Then the following hold.
\begin{enumerate}
\item
If $d_1 \geq 3$, then 
$\varphi$ is a homeomorphism onto its image. 
\item
If $d_1 \geq 2 (n-1)!$, then $\varphi$ is unimodular. 
\item
If $d_1 \geq \max \{ 3 , 2 (n-1)! \}$, then 
$\varphi$ is a faithful embedding. 
\end{enumerate}
\end{Theorem}

Since (3) is an immediate consequence of (1) and (2), the essential parts of the theorem are (1) and (2).

If $L$ is a positive tropical line bundle of polarization type 
$(d_1, \ldots, d_n)$, then for any integer $d \geq 1$, $d L$ has polarization type $(d d_1, \ldots, d d_n)$; see Remark~\ref{remark:additive}. By 
Theorem~\ref{thm:main:1-1}~(3), we obtain the following corollary. 

\begin{Corollary}
\label{cor:main:1:2}
Let $L$ be a positive tropical line bundle on a 
tropical abelian variety of dimension $n$. Then 
for any $d \geq \max \{ 3 , 2 (n-1)! \}$, $dL$ admits a faithful embedding of $X$.
\end{Corollary}

\begin{Remark} \label{remark:twoexpressionofFE}
We note that Corollary~\ref{cor:main:1:2} implies Theorem~\ref{thm:main:1-1}~(3), and thus they are equivalent to each other. To see that, let $L = L (\lambda , \gamma)$ be a positive tropical line bundle with type $(d_1 , \ldots , d_n)$ and suppose that $d_1 \geq 3$. Then $\frac{1}{d_1} \lambda$ is a homomorphism from $M'$ to $M$, and in fact, it is a polarization. Further, $\left( \frac{1}{d_1} \lambda , \frac{1}{d_1} \gamma \right)$ is a tropical descent datum, and $L = L (\lambda , \gamma) = d_1 L \left(\frac{1}{d_1} \lambda , \frac{1}{d_1} \gamma \right)$ by Remark~\ref{remark:additive}. It follows that Corollary~\ref{cor:main:1:2} implies Theorem~\ref{thm:main:1-1}~(3). Note that this kind of equivalence also holds for line bundles on abelian varieties over $\CC$.
\end{Remark}

Recall that in the classical Lefschetz theorem over $\CC$, for a positive line bundle $L$, $d \geq 3$ is sufficient for $L^{\otimes d}$ to be very ample. (As we mentioned in Remark~\ref{remark:twoexpressionofFE}, this is equivalent that for a positive line bundle on $L$ of type $(d_1 , \ldots , d_n)$, if $d_1 \geq 3$, then $L$ is very ample.)
It would be natural to ask the following question: for a positive tropical line bundle $L$ on a real torus $X$, if $d \geq 3$, then does $dL$ admit a faithful embedding of $X$?
When $n \leq 2$, since $\max \{ 3 , 2 (n-1)! \} = 3$, we know that the answer is affirmative by Corollary~\ref{cor:main:1:2}, but we do not know if this is affirmative or not when $n \geq 3$. Even if the answer to the question is negative, it would be an interesting problem to find a sharp bound for $d$ such that $dL$ admits a faithful embedding.

Next, we consider tropical line bundles that are not necessarily positive. 

\begin{Proposition} \label{prop:notadmitFE}
Let $L$ be a tropical line bundle on a real torus $X$. 
Suppose that $L$ is not positive. Then for any $d \in \ZZ_{\geq 1}$ and for any $s_0 , \ldots , s_d \in H^0 (X , L) \setminus \{ + \infty \}$, there exists a real subtorus of positive dimension $Y \subseteq X$ such that $\varphi_{s_0 , \ldots , s_d} (Y)$ is a singleton. 
\end{Proposition}

\Proof
Since $L$ is not positive, we have $\dim (X) \geq 1$. 
We may assume that $L$ has a nontrivial regular global section. By Proposition~\ref{prop:Globalsections:Sumi}, there exist a real torus $X'$ of dimension less than $\dim (X)$, a surjective homomorphism $f : X \to X'$, and a line bundle on $L'$ on $Y$ such that $L \cong f^{\ast} (L')$ and $f^{\ast} : H^0 (Y,L') \to H^0 (X,L)$ is an isomorphism of $\TT$-modules.

We take any $s_0 , \ldots , s_d \in H^0(X,L) \setminus  \{ + \infty \}$. Then there exist $t_0 , \ldots , t_d \in H^0 (Y,L')$ such that $s_i = f^{\ast} (t_i)$ for $i=0 , \ldots , d$. It tuns out that $\varphi_{s_0 , \ldots , s_d} (\Ker (f))$ is a singleton. Let $Y$ be the connected component of $\Ker (f)$ with $0$. Since $\dim (Y) \geq 1$ and $\varphi_{s_0 , \ldots , s_d} (Y) \subseteq \varphi_{s_0 , \ldots , s_d} (\Ker (f))$, this $Y$ suffices.
\QED

Here is an immediate application of Theorem~\ref{thm:main:1-1}. 
In the classical theory of complex tori, we know that a complex torus has an integral positive-definite hermitian form on the universal cover if and only if it can be embedded in a projective space by a line bundle. Since we define a tropical abelian variety as a real torus that satisfies a tropical analogue of Riemann's period relations (see Definition~\ref{def:trop:av}), the following equivalence between (i) and (ii) is regarded as a tropical analogue of that classical fact.

\begin{Corollary} \label{cor:trop:abel:proj}
Let $X$ be a real torus. Then the following are equivalent to each other:
\begin{enumerate}
\item[(i)]
$X$ is a tropical abelian variety;
\item[(ii)]
there exists a tropical line bundle on $X$ that admits a faithful embedding;
\item[(iii)]
there exist a tropical line bundle $L$ on $X$ and regular global sections $s_0 , \ldots , s_m \in H^0 (X, L) \setminus \{ + \infty \}$ such that the map $\varphi_{s_0,\ldots , s_m} : X \to \TT\PP^m$ is injective.
\end{enumerate}
\end{Corollary}

\Proof
The implication from (i) to (ii) follows from Theorem~\ref{thm:main:1-1}~(3). It is obvious that (ii) implies (iii). The implication from (iii) to (i) is a consequence of Proposition~\ref{prop:notadmitFE}.
\QED

\subsection{Alternative description of tropical theta functions} \label{subsec:Alternative}

In order to prove Theorem~\ref{thm:main:1-1}, 
in this subsection, we describe the tropical theta functions in terms of the pair corresponding to a tropical descent datum by Lemma~\ref{lem:polarization:equiv}.

Let $(Q,\ell)$ be a pair of a symmetric $\RR$-bilinear form $Q$ on $N_{\RR} \times N_{\RR}$ such that $Q (M' \times N) \subseteq \ZZ$ and a homomorphism $\ell : M' \to \RR$,  with corresponding tropical descent datum $(\lambda , \gamma)$. Let $T(Q,\ell)$ denote the $\TT$-semimodule of tropical theta functions with respect to $(\lambda ,\gamma)$; by Lemma~\ref{lem:polarization:equiv}~(2), we have
\begin{equation} \label{eqn:TQell}
T(Q,\ell) := \{ \vartheta \in H^0 (N_{\RR} , N_{\RR} \times \TT) \mid \text{$\vartheta$ satisfies (\ref{eqn:quasi:periodicity:0})}\}
.
\end{equation}
Note that for any $\vartheta \in T(Q , \ell)$ and $\vartheta' \in T(Q', \ell')$, $\vartheta + \vartheta' \in T(Q+Q' , \ell + \ell')$.
Recall that the pullback homomorphism $\pi^{\ast} : H^0 (X,L(\lambda , \gamma)) \to H^0 (N_{\RR} , N_{\RR} \times \TT)$ induces an isomorphism $\pi^{\ast} : H^0 (X,L(\lambda , \gamma)) \to T(Q,\ell)$; see the paragraph after Theorem~\ref{thm:trop:Appell:Humbert}.

Assume that $\lambda$ is a polarization, i.e., $Q$ is positive-definite. For any $b \in M$, we describe the theta function $\vartheta_b^{(\lambda , \gamma)}$ in terms of $(Q,\ell)$. Since $Q$ is nondegenerate,  there exists a unique $r \in N_\RR$ such that $\ell(\cdot) = Q(\cdot, r)$. By Remark~\ref{remark:polarization:injective},  the $\RR$-linear extension $\lambda_\RR\colon  M^\prime_\RR \to M_\RR$ is an isomorphism. Note that $\lambda_{\RR}^{-1} (b) \in M'_{\RR} = N_{\RR}$. We set
\begin{equation}
\label{eqn:def:theta:b}
\vartheta_b^{(Q,\ell)} (x) =
\min_{u^\prime \in M^\prime} 
\left(
Q\!\left(x,\, \lambda_\RR^{-1}(b) + u^{\prime} \right) + \frac{1}{2} Q\!\left(u^\prime + \lambda_\RR^{-1}(b) - r,\, u^\prime + \lambda_\RR^{-1}(b) - r\right)
\right) , 
\end{equation}
which 
is nothing but (up to sign convention) the tropical theta function given in \cite[Subsection~3.4]{Sumi}.
A direct computation shows that
\begin{equation} \label{eqn:def:theta:b2}
\vartheta_b^{(Q,\ell)} = \vartheta_b^{(\lambda , \gamma)} + \frac{1}{2} Q (  \lambda_{\RR}^{-1} (b) - r ,  \lambda_{\RR}^{-1} (b) - r )
.
\end{equation}
Note in particular that $\vartheta_b^{(Q,\ell)} - \vartheta_b^{(\lambda , \gamma)} \in \RR$.

\begin{Remark} \label{remark:construction:trop:theta}
Mikhalkin--Zharkov \cite{MZ} defines $\vartheta_0^{(Q,\ell)}$ for tropical Jacobians. They also mentioned the idea how to construct $\vartheta_b^{(Q,\ell)}$'s for any real torus, and Sumi \cite{Sumi} gives an explicit construction for them.
\end{Remark}

The following theorem is proved in \cite[\S3.4, Theorem~36]{Sumi}.

\begin{Theorem}[{\cite[Remark~5.5]{MZ}}, {\cite[\S3.4, Theorem~36]{Sumi}}]
\label{thm:generators0}
Let $(Q,\ell)$ and $(\lambda , \gamma)$ be as above. Assume that $\lambda$ is a polarization.
Let $\mathfrak{B} \subseteq M$ be a complete system of representatives of $M/\lambda (M')$. Then $\left\{\vartheta_b^{(Q , \ell)}\right\}_{b \in \mathfrak{B}}$ 
generates $T(Q,\ell)$ as a $\TT$-semimodule, 
i.e., for any tropical theta function $\vartheta \in T(Q ,\ell)$,
there exist $c_b \in \TT$ for all $b \in \mathfrak{B}$ such that $\vartheta = \min_{b \in \mathfrak{B}} \left(\vartheta_b^{(Q , \ell)} + c_b \right)$. 
\end{Theorem}

\begin{Remark} \label{remark:theta:transf:inv}
For any  $b \in M$ and $\widetilde{u}^\prime \in M^\prime$, we have
\begin{align*}
\begin{split}
&\vartheta_{b+\lambda(\widetilde{u}^\prime)}^{(Q,\ell)}(x) 
\\
%& 
& = 
\min_{u^\prime \in M^\prime} 
\left(
Q\!\left(x,\,  u^\prime +  \widetilde{u}^\prime+ \lambda_\RR^{-1}(b) \right) + \frac{1}{2} Q\!\left(u^\prime + \widetilde{u}^\prime + \lambda_\RR^{-1}(b) - r,\, u^\prime + \widetilde{u}^\prime + \lambda_\RR^{-1}(b) - r\right) 
\right)
\\
&
=  \vartheta_{b}^{(Q,\ell)}(x) 
\end{split}
\end{align*}
as functions on $x \in N_{\RR}$. It follows that $\vartheta_b^{(Q,\ell)}$ depends only on the class of $b$ in $M/\lambda (M')$, and thus $\left\{ \vartheta_b^{(Q , \ell)} \right\}_{b \in \mathfrak{B}}$ does not depend on the choice of the complete system of representatives $\mathfrak{B}$. On the other hand, $\left\{ \vartheta_b^{(\lambda , \gamma)} \right\}_{b \in \mathfrak{B}}$ depends on the choice of $\mathfrak{B}$ (cf. (\ref{eqn:def:theta:b2})). 
\end{Remark}

Let $(d_1, \ldots, d_n)$ be the type of the polarization $\lambda$ on $X$, where $n = \dim (X) = \dim (N_{\RR})$. We set $D := |M' / \lambda (M)| = d_1 \cdots d_n$. Let $\mathfrak{B} \subseteq M$ be a complete system of representatives of $M / \lambda (M')$. We write $\mathfrak{B} = \{b_1, \ldots, b_{D}\}$. For all $i=1 , \ldots , D$, we consider the tropical theta function $\vartheta_{b_i}^{(Q,\ell)}$ given by (\ref{eqn:def:theta:b}). For each $i=1 , \ldots , D$, let $\bar{\vartheta}_{b_{i}}^{(Q,\ell)}$ denote the regular global section of $L(\lambda , \gamma)$ such that $\pi^{\ast} (\bar{\vartheta}_{b_{i}}^{(Q,\ell)}) = \vartheta_{b_i}^{(Q,\ell)}$. Then we have a map $\varphi_{\bar{\vartheta}_{b_1}^{(Q,\ell)}, \ldots ,\bar{\vartheta}_{b_{D}}^{(Q,\ell)}}\colon X \to \TT\PP^{D-1}$ defined just before Definition~\ref{def:FE} for $\bar{\vartheta}_{b_1}^{(Q,\ell)} , \ldots , \bar{\vartheta}_{b_D}^{(Q,\ell)} \in H^0 (X , L(\lambda ,\gamma)) \setminus  \{ + \infty \}$.

Note that if we change the numbering of the elements of $\mathfrak{B}$, then only the ordering of the homogeneous coordinates changes. As in Remark~\ref{remark:theta:transf:inv}, we have $\left\{ \vartheta_b^{(Q,\ell)} \right\}_{b \in \mathfrak{B}} = \left\{ \vartheta_b^{(Q,\ell)} \right\}_{b \in M}$, so 
$\left\{ \vartheta_b^{(Q,\ell)} \right\}_{b \in \mathfrak{B}}$ 
 is independent of the choice of $\mathfrak{B}$. It follows that the map $\varphi_{\bar{\vartheta}_{b_1}^{(Q,\ell)} , \ldots , \bar{\vartheta}_{b_D}^{(Q,\ell)}}\colon X \to \TT\PP^{D-1}$ modulo the ordering of homogeneous coordinates of $\TT\PP^{D-1}$ does not depend on the choice of $\mathfrak{B}$ or the numbering of the elements of $\mathfrak{B}$. Thus when we consider this map modulo the ordering of homogeneous coordinates of $\TT\PP^{D-1}$, we write $\varphi^{(Q,\ell)}$ for it. Since the properties that a map $X \to \TT\PP^{D-1}$ is injective, unimodular, and a faithful embedding, respectively, does not depend on the choice of the ordering of homogeneous coordinates, it makes sense to ask whether or not $\varphi^{(Q,\ell)}$ has such a property.

Since $\vartheta_{b}^{(Q,\ell)} - \vartheta_b^{(\lambda , \gamma)} \in \RR$ for any $b \in \mathfrak{B}$, it follows Remark~\ref{remark:unimdular:part}~(1) that Theorem~\ref{thm:main:1-1} is equivalent to the following theorem.

\begin{Theorem}
\label{thm:main:1}
Let $X =N_{\RR} /M'$ be a tropical abelian variety, and let 
$(Q,\ell)$ be a pair of a polarization $Q$ on $X$ and a homomorphism $\ell : M' \to \RR$.  Then the following hold. 
\begin{enumerate}
\item
If $d_1 \geq 3$, then 
$\varphi^{(Q,\ell)}$ is a homeomorphism onto its image. 
\item
If $d_1 \geq 2 (n-1)!$, then the map
$\varphi^{(Q,\ell)}$ is unimodular. 
\item
If $d_1 \geq \max \{ 3 , 2 (n-1)! \}$, then 
$\varphi^{(Q,\ell)}$ is a faithful embedding. 
\end{enumerate}
\end{Theorem}

In Theorem~\ref{thm:main:1}, since (3) is an immediate consequence of (1) and (2), we only have to prove (1) and (2). We will prove (1) in Section~\ref{section:injectivity} and (2) in Section~\ref{section:Pf:unimodular}.

\begin{Remark} \label{remark:suffice:injectivity}
Since $X$ is a compact Hausdorff space, it suffices for (1) above to show that $\varphi^{(Q,\ell)}$ is injective if $d_1 \geq 3$.
\end{Remark}

\begin{Remark} \label{remrak:lift:UM:FE}
The canonical surjection $\pi : N_{\RR} \to X$ is a local isomorphism of rational polyhedral spaces, and we have 
$\varphi_{\bar{\vartheta}_{b_1}^{(Q,\ell)}, \ldots ,\bar{\vartheta}_{b_{D}}^{(Q,\ell)}} \circ \pi = \varphi_{\vartheta_{b_1}^{(Q,\ell)}, \ldots ,\vartheta_{b_{D}}^{(Q,\ell)}}$. This means that $\varphi_{\bar{\vartheta}_{b_1}^{(Q,\ell)}, \ldots ,\bar{\vartheta}_{b_{D}}^{(Q,\ell)}}$ is unimodular if and only if so is $\varphi_{\vartheta_{b_1}^{(Q,\ell)}, \ldots ,\vartheta_{b_{D}}^{(Q,\ell)}}$. Further $\varphi_{\bar{\vartheta}_{b_1}^{(Q,\ell)}, \ldots ,\bar{\vartheta}_{b_{D}}^{(Q,\ell)}}$ is injective if and only if for any $y,y' \in N_{\RR}$ with $\varphi_{\vartheta_{b_1}^{(Q,\ell)}, \ldots ,\vartheta_{b_{D}}^{(Q,\ell)}} (y) = \varphi_{\vartheta_{b_1}^{(Q,\ell)}, \ldots ,\vartheta_{b_{D}}^{(Q,\ell)}} (y')$, we have $y - y ' \in M'$.
Thus in order to study $\varphi_{\bar{\vartheta}_{b_1}^{(Q,\ell)}, \ldots ,\bar{\vartheta}_{b_{D}}^{(Q,\ell)}}$, we will often focus on $\varphi_{\vartheta_{b_1}^{(Q,\ell)} , \ldots , \vartheta_{b_D}^{(Q,\ell)}}$.
\end{Remark}

Let us explain the reason why we use both $\vartheta_b^{(\lambda ,\gamma)}$ and $\vartheta_b^{(Q,\ell)}$. The tropical theta functions $\vartheta_b^{(Q,\ell)}$ have an advantage that $\left\{ \vartheta_b^{(Q,\ell)} \right\}_{b \in \mathfrak{B}}$ does  not depend on the choice of a complete system $\mathfrak{B}$ of representatives of $M /\lambda (M')$, while $\left\{ \vartheta_b^{(\lambda,\gamma)} \right\}_{b \in \mathfrak{B}}$ does. On the other hand, as we will see in Subsection~\ref{subsection:construct:lift} and Section~\ref{section:surjective:trop:theta}, the tropical theta functions $\vartheta_b^{(\lambda , \gamma)}$ have an advantage when we consider the lifting of tropical theta functions to nonarchimdean theta functions.

\subsection{Examples}

We illustrate Theorem~\ref{thm:main:1} 
by giving examples
in a simple case. 
%%%%%%%
% Example %
%%%%%%%
\begin{Example}
\label{eg:trop:ell:faithful}
We fix $\varpi > 0$, and we take $M^\prime =\varpi \ZZ \subseteq \RR$ and $N = \ZZ$, so that $X = \RR/\varpi \ZZ$ with integral structure $\ZZ$. 
Namely, we consider a tropical abelian variety of dimension one (i.e., a tropical elliptic curve). 

Let $Q\colon \RR \times \RR \to \RR$ be an $\RR$-bilinear form with $Q(\varpi \ZZ \times \ZZ) \subseteq \ZZ$.  Then there is an integer~$d$ such that 
$Q = d \varpi^{-1}$. This is a polarization on $X$ if and only if $d \geq 1$, and then it has type~$(d)$. The corresponding $\ZZ$-linear map $\lambda\colon \varpi \ZZ \to \ZZ$ (see Lemma~\ref{lem:polarization:equiv}) is given by $\varpi k \mapsto d\, k$. Theorem~\ref{thm:main:1} 
in particular asserts that $\varphi^{(d \varpi^{-1}, 0)}$ is unimodular when $d = 2$ and is a 
faithful embedding whenever $d \geq 3$. Let us see this by explicit computation 
when $d = 2, 3$, and $4$. 

\begin{enumerate}
\item[(1)] 
First,  let $d = 2$.  
We show that $\varphi^{(2 \varpi^{-1}, 0)}$ is unimodular. 
A subset $\mathfrak{B} = \{0, 1\} \subseteq \ZZ$ is a complete system of representatives of $\ZZ / \lambda (\varpi \ZZ)$. For $b \in \mathfrak{B}$, the tropical theta function $\vartheta_b\colon \RR \to \RR$ with respect to $(2 \varpi^{-1}, 0)$ becomes  
\[
\vartheta_b(x) 
= \min_{m^\prime \in \varpi \ZZ} 
\left(
(2  \varpi^{-1} m^\prime + b) x  + \frac{\varpi}{4} (2  \varpi^{-1} m^\prime + b)^2 \right)
= \min_{k \in \ZZ} 
\left(
(2 k + b) x  + \frac{\varpi}{4} (2 k + b)^2 \right). 
\] 
By a straightforward computation, for $x \in [0, \varpi)$, we obtain
{\allowdisplaybreaks
\begin{align*}
& 
\vartheta_0(x)
= 
\begin{cases}
0 & (x \in \left[0, \frac{1}{2}\varpi\right]) \\
- 2 x + \varpi & (x \in \left[\frac{1}{2}\varpi, 1\right)) 
\end{cases}, 
\quad 
\vartheta_1(x)
= -x + \frac{\varpi}{4}. 
\end{align*}
}
By adjusting the first coordinate of $\TT\PP^{1}$ to $0$, the map $\widetilde{\varphi}_{(\vartheta_0 , \vartheta_1)}\colon \RR \to \RR$ is given by  
\[
\widetilde{\varphi}_{\vartheta_0 , \vartheta_1}(x) 
= (\vartheta_1(x) - \vartheta_0(x)) 
= 
\begin{cases}
-x + \frac{\varpi}{4} & (x \in \left[0, \frac{1}{2}\varpi\right]) \\
x - \frac{3\varpi}{4} & (x \in \left[\frac{1}{2}\varpi, 1\right)) 
\end{cases}. 
\]
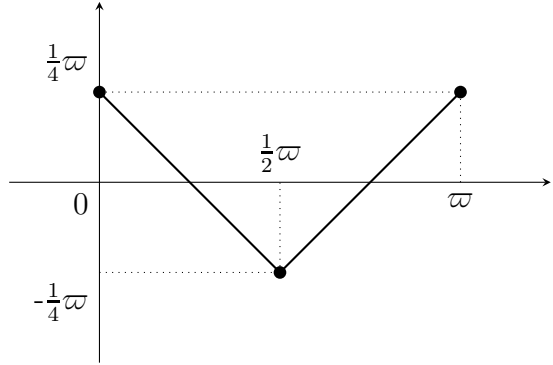
\[
\begin{tikzpicture}[scale = 1.2]
  \begin{scope}
  \draw[->, >=stealth] (-1, 0) -- (5, 0);
  \draw[->, >=stealth] (0, -2) -- (0, 2);
  \draw[thick] (0,1) -- (2, -1) -- (4, 1);
  \fill (0, 1) circle [radius=0.07] 
  node[above left]{$\frac{1}{4}\varpi$}; 
  \fill (2, -1) circle [radius=0.07]; 
  \fill (4, 1) circle [radius=0.07]; 
   \draw[dotted] (2,-1) -- (2, 0);
  \draw[dotted] (4,0) -- (4, 1);
  \draw[dotted] (0,1) -- (4, 1);
  \draw[dotted] (0,-1) -- (2, -1);
  \draw (0, -1) node[below left]{-$\frac{1}{4}\varpi$}; 
  \draw (2, 0) node[above]{$\frac{1}{2}\varpi$}; 
  \draw (4, 0) node[below]{$\varpi$}; 
  \draw (0, 0) node[below left]{$0$}; 
  \end{scope}
\end{tikzpicture}
\]
Since the slope of $\widetilde{\varphi}_{\vartheta_0 , \vartheta_1}$ is 
$-1$ on $\left[0, \frac{1}{2}\varpi\right]$ and $1$ 
on $\left[\frac{1}{2}\varpi, \varpi\right)$, 
we see that this map and hence $\varphi^{(2\varpi^{-1}, 0)}$ is unimodular.  We remark that $\widetilde{\varphi}_{\vartheta_0 , \vartheta_1}$ is not injective, so $\varphi^{(2\varpi^{-1}, 0)}$ is not a faithful embedding.
\item[(2)]
Let $d = 3$. 
We prove that $\varphi^{(3 \varpi^{-1}, 0)}$ is a faithful embedding.
The subset $\mathfrak{B} = \{0, 1, 2\} \subseteq \ZZ$ is a complete system of representatives of $\ZZ / \lambda (\varpi \ZZ)$. For $b \in \mathfrak{B}$, the tropical theta function $\vartheta_b\colon \RR \to \RR$ with respect to $(3 \varpi^{-1}, 0)$ becomes  
\[
\vartheta_b(x) 
= \min_{k \in \ZZ} 
\left(
(3 k + b) x  + \frac{\varpi}{6} (3 k + b)^2 \right). 
\] 
By a straightforward computation, for $x \in [0, \varpi)$, we obtain
{\allowdisplaybreaks
\begin{align*}
& 
\vartheta_0(x)
= 
\begin{cases}
0 & (x \in \left[0, \frac{1}{2}\varpi\right]) \\
- 3 x + \frac{3}{2} \varpi & (x \in \left[\frac{1}{2}\varpi, \varpi\right)) 
\end{cases}, 
\quad 
\vartheta_1(x)
= \begin{cases}
x + \frac{1}{6}\varpi & (x \in \left[0, \frac{1}{6}\varpi\right]) \\
-2 x  + \frac{2}{3}\varpi & (x \in \left[\frac{1}{6}\varpi, \varpi\right)) 
\end{cases}, 
\\
& 
\vartheta_2(x)
= 
\begin{cases}
-x + \frac{1}{6}\varpi & (x \in \left[0, \frac{5}{6}\varpi\right]) \\
-4x + \frac{8}{3}\varpi & (x \in \left[\frac{5}{6}\varpi, \varpi\right)) 
\end{cases}. 
\end{align*}
}

By adjusting the first coordinate of $\TT\PP^{2}$ to $0$, the map $\widetilde{\varphi}_{\vartheta_0 , \vartheta_1 , \vartheta_2}\colon \RR \to \RR^{2}$ is given by  
$
\widetilde{\varphi}_{\vartheta_0 , \vartheta_1 , \vartheta_2}(x) = (\vartheta_1(x) - \vartheta_0(x), \vartheta_{2}(x) - \vartheta_0(x)) 
$. 
It follows that, for $x \in [0, \varpi)$,  
\[
\widetilde{\varphi}_{\vartheta_0 , \vartheta_1 , \vartheta_2}(x) 
= \begin{cases}
\left(
x +\frac{1}{6}\varpi, \, -x +\frac{1}{6}\varpi 
\right) & (x \in \left[0, \frac{1}{6}\varpi\right]),  
\\
\left(
-2x + \frac{2}{3}\varpi, \,-x +\frac{1}{6}\varpi
\right) & (x \in \left[\frac{1}{6}\varpi, \frac{1}{2}\varpi\right]), 
\\
\left(
x -\frac{5}{6}\varpi, \,2x - \frac{4}{3}\varpi
\right) & (x \in \left[\frac{1}{2}\varpi, \frac{5}{6}\varpi\right]), 
\\
\left(
x -\frac{5}{6}\varpi, \,-x +\frac{7}{6}\varpi 
\right) & (x \in \left[\frac{5}{6}\varpi, \varpi\right)). 
\end{cases}
\] 
The image of $\rest{\widetilde{\varphi}_{\vartheta_0 , \vartheta_1 , \vartheta_2}}{[0, \varpi)}$ is a triangle and each edge has lattice length $\frac{\varpi}{3}$, and indeed the map $\varphi_{\bar{\vartheta}_0 , \bar{\vartheta}_1, \bar{\vartheta}_2} : \ZZ / \varpi \ZZ \to \TT\PP^2$ is a faithful embedding, and thus $\varphi^{(3\varpi^{-1}, 0)}$ is a faithful embedding.
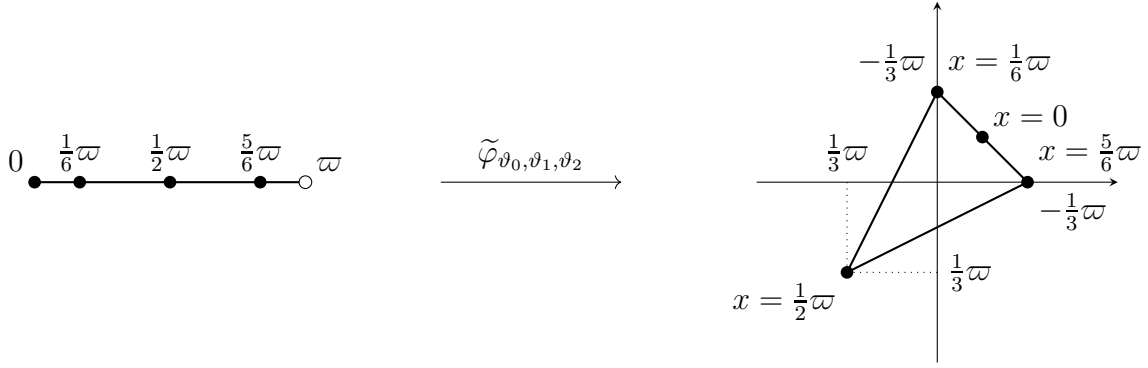
\[
\begin{tikzpicture}[scale = 1.2]
  \draw[thick] (0,0) -- (2.93, 0);
  \fill (0, 0) circle [radius=0.07]  node[above left]{$0$}; 
  \fill (0.5, 0) circle [radius=0.07] node[above]{$\frac{1}{6}\varpi$}; 
  \fill (1.5, 0) circle [radius=0.07] node[above]{$\frac{1}{2}\varpi$}; 
  \fill (2.5, 0) circle [radius=0.07] node[above]{$\frac{5}{6}\varpi$}; 
  \draw (3, 0) circle [radius=0.07]  node[above right]{$\varpi$};  
  \draw[->] (4.5, 0) -- (6.5, 0);
  \draw (5.5, 0) node[above]{$\widetilde{\varphi}_{\vartheta_0 , \vartheta_1 , \vartheta_2}$}; 
  \begin{scope}[shift={(10, 0)}]
  \draw[->, >=stealth] (-2, 0) -- (2, 0);
  \draw[->, >=stealth] (0, -2) -- (0, 2);
  \draw[thick] (1,0) -- (0, 1) -- (-1, -1) -- (1, 0) -- cycle;
  \fill (1, 0) circle [radius=0.07] node[below right]{$-\frac{1}{3}\varpi$}
  node[above right]{$x = \frac{5}{6}\varpi$}; 
  \fill (0, 1) circle [radius=0.07] node[above left]{$-\frac{1}{3}\varpi$} 
  node[above right]{$x = \frac{1}{6}\varpi$}; 
  \fill (-1, -1) circle [radius=0.07] node[below left]{$x = \frac{1}{2}\varpi$}; 
  \fill (0.5, 0.5) circle [radius=0.07] node[above right]{$x = 0$}; 
  \draw[dotted] (-1,0) -- (-1, -1);
  \draw[dotted] (0,-1) -- (-1, -1);
  \draw (-1, 0) node[above]{$\frac{1}{3}\varpi$}; 
  \draw (0, -1) node[right]{$\frac{1}{3}\varpi$}; 
  \end{scope}
\end{tikzpicture}
\]
\item[(3)]
Let $d = 4$. 
The subset $\mathfrak{B} = \{0, 1, 2, 3 \} \subseteq \ZZ$ is a complete system of representatives of $\ZZ / \lambda (\varpi \ZZ)$. 
For $b \in \mathfrak{B}$, the tropical theta function $\vartheta_b\colon \RR \to \RR$ with respect to $(4 \varpi^{-1}, 0)$ becomes  
\[
\vartheta_b(x) 
= \min_{k \in \ZZ} 
\left(
(4 k + b) x  + \frac{\varpi}{8} (4 k + b)^2 \right). 
\] 
By a straightforward computation, for $x \in [0, \varpi)$, we obtain
\begin{align*}
& 
\vartheta_0(x)
= 
\begin{cases}
0 & (x \in \left[0, \frac{1}{2}\varpi\right]) \\
- 4 x + 2 \varpi & (x \in \left[\frac{1}{2}\varpi, \varpi\right)) 
\end{cases}, 
&
\vartheta_1(x)
= \begin{cases}
x + \frac{1}{8}\varpi & (x \in \left[0, \frac{1}{4}\varpi\right]) \\
-3 x  + \frac{9}{8}\varpi & (x \in \left[\frac{1}{4}\varpi, \varpi\right)) 
\end{cases}, 
\\
& 
\vartheta_2(x)
= 
-2x + \frac{1}{2}\varpi  \quad (x \in \left[0, \varpi\right]) 
& 
\vartheta_3(x)
= 
\begin{cases}
-x + \frac{1}{8}\varpi & (x \in \left[0, \frac{3}{4}\varpi\right]) \\
-5x + \frac{25}{8}\varpi & (x \in \left[\frac{3}{4}\varpi, \varpi\right)) 
\end{cases}. 
\end{align*}
Thus the map 
\[
\widetilde{\varphi}_{\vartheta_0 , \vartheta_1 , \vartheta_2 ,\vartheta_3}\colon \RR \to \RR^{3}, 
\quad 
 x \mapsto \left(\vartheta_1(x) - \vartheta_0(x), \vartheta_{2}(x) - \vartheta_0(x), \vartheta_3(x) - \vartheta_0(x)\right) 
\]
is given by 
\[
\widetilde{\varphi}_{\vartheta_0 , \vartheta_1 , \vartheta_2 ,\vartheta_3}
= \begin{cases}
\left(
x -\frac{1}{8}\varpi, \, 
-2x + \frac{1}{2}\varpi, \, 
-x +\frac{1}{8}\varpi 
\right) & (x \in \left[0, \frac{1}{4}\varpi\right]),  
\\
\left(
-3x +\frac{9}{8}\varpi, \, 
-2x + \frac{1}{2}\varpi, \, 
-x +\frac{1}{8}\varpi 
\right) & (x \in \left[\frac{1}{4}\varpi, \frac{1}{2}\varpi\right]), 
\\
\left(
x -\frac{7}{8}\varpi, \, 
2x - \frac{3}{2}\varpi, \, 
3 x -\frac{15}{8}\varpi 
\right) & (x \in \left[\frac{1}{2}\varpi, \frac{3}{4}\varpi\right]), 
\\
\left(
x -\frac{7}{8}\varpi, \, 
2x - \frac{3}{2}\varpi, \, 
-x + \frac{9}{8}\varpi 
\right) & (x \in \left[\frac{3}{4}\varpi, \varpi\right)). 
\end{cases}
\] 
The image of $\rest{\widetilde{\varphi}_{\vartheta_0 , \vartheta_1 , \vartheta_2 ,\vartheta_3}}{[0, \varpi)}$ is a quadrangle and each edge has lattice length $\frac{\varpi}{4}$, and indeed $\varphi_{\bar{\vartheta_0} , \bar{\vartheta}_1 , \bar{\vartheta}_2 ,\bar{\vartheta}_3}\colon \RR/\varpi\ZZ \to \RR^3$ is a faithful embedding. Thus $\varphi^{(4\varpi , 0)}$ is a faithful embedding.
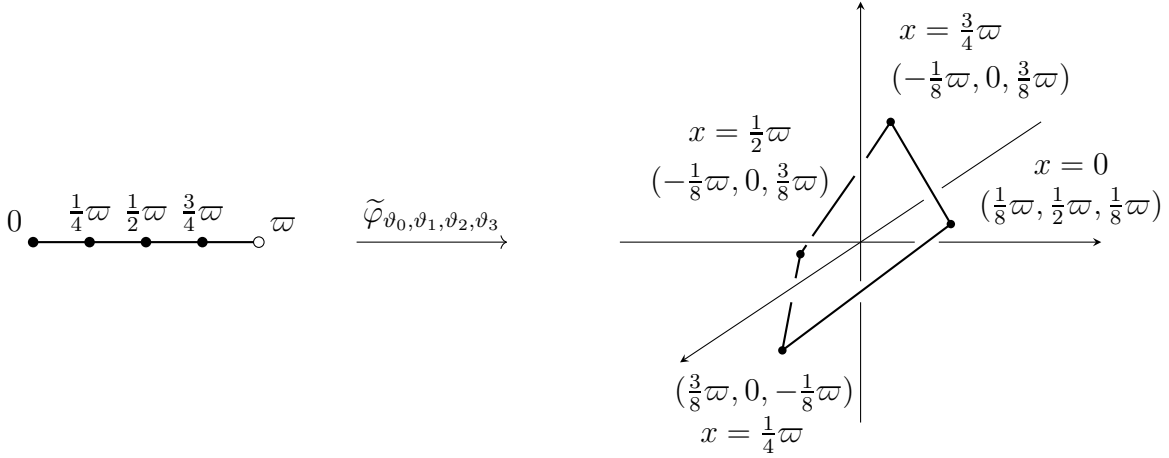
\[
\begin{tikzpicture}
  \draw[thick] (0,0) -- (2.93, 0);
  \fill (0, 0) circle [radius=0.07]  node[above left]{$0$}; 
  \fill (0.75, 0) circle [radius=0.07] node[above]{$\frac{1}{4}\varpi$}; 
  \fill (1.5, 0) circle [radius=0.07] node[above]{$\frac{1}{2}\varpi$}; 
  \fill (2.25, 0) circle [radius=0.07] node[above]{$\frac{3}{4}\varpi$}; 
  \draw (3, 0) circle [radius=0.07]  node[above right]{$\varpi$};  
  \draw[->] (4.3, 0) -- (6.3, 0);
  \draw (5.3, 0) node[above]{$\widetilde{\varphi}_{\vartheta_0 , \vartheta_1 , \vartheta_2 ,\vartheta_3}$}; 
  \begin{scope}[shift={(11, 0)}, scale = 0.8]
  \draw (-4, 0) -- (0.9, 0);
  \draw[->, >=stealth] (1.3, 0) -- (4, 0);
  \draw (0, -3) -- (0, -1);
  \draw[->, >=stealth] (0, -0.65) -- (0, 4);
   \draw (3, 2) -- (1.35, 0.9);
  \draw[->, >=stealth] (1.05, 0.7) -- (-3, -2);  
  \fill (-1.3, -1.8) circle [radius=0.07]; 
  \fill (1.5, 0.3) circle [radius=0.07]; 
  \fill (0.5, 2) circle [radius=0.07]; 
  \fill (-1, -0.2) circle [radius=0.07];  
  \draw[thick] (-1.3, -1.8) -- (1.5, 0.3);
  \draw[thick] (1.5, 0.3) -- (0.5, 2);
  \draw[thick] (0.5, 2) -- (0.1, 1.4);
   \draw[thick] (-0.9, -0.05) -- (-1, -0.2);
  \draw[thick] (-0.8, 0.1) -- (-0.1, 1.1);
    \draw[thick] (-1, -0.2) -- (-1.08, -0.6);
   \draw[thick] (-1.15, -1) -- (-1.3, -1.8);
   \draw (-1.6, -2.5) node{$(\frac{3}{8}\varpi, 0, -\frac{1}{8}\varpi)$}; 
   \draw (-1.8, -3.2) node{$x = \frac{1}{4}\varpi$}; 
   \draw (3.5, 0.6) node{$(\frac{1}{8}\varpi, \frac{1}{2}\varpi, \frac{1}{8}\varpi)$}; 
   \draw (3.5, 1.3) node{$x = 0$}; 
   \draw (2, 2.7) node{$(-\frac{1}{8}\varpi, 0, \frac{3}{8}\varpi)$}; 
   \draw (1.5, 3.5) node{$x = \frac{3}{4}\varpi$}; 
   \draw (-2, 1) node{$(-\frac{1}{8}\varpi, 0, \frac{3}{8}\varpi)$}; 
   \draw (-2, 1.8) node{$x = \frac{1}{2}\varpi$};   
  \end{scope}
\end{tikzpicture}
\]
\end{enumerate}
\end{Example}

%%%%%%%%%%%%%
% Reduction to r= 0. %%
%%%%%%%%%%%%%

\subsection{Reduction to the case $\ell = 0$}
Let $(Q , \ell)$ be a pair of a symmetric $\RR$-bilinear form $Q : N_{\RR} \times N_{\RR} \to \RR$ such that $Q(M' \times N) \subseteq \ZZ$ and a homomorphism $\ell : M' \to \RR$. In this subsection, assume that $Q$ is a polarization, i.e., positive-definite, and we will see that our main results Theorem~\ref{thm:main:1}~(1), (2) are reduced to the results for $(Q,\ell)$ with $\ell = 0$. 

We keep the notation in Subsection~\ref{subsec:Alternative}; $X = N_{\RR}/M'$ is a real torus and $\pi : N_{\RR} \to X$ is the canonical surjection. Recall that $T(Q,\ell)$ be the $\TT$-semimodule of tropical theta functions with respect to the tropical descent datum corresponding to the pair $(Q,\ell)$.

First, we show how a translation changes a tropical theta function. For any $v \in N_{\RR}$, let $t_v\colon N_\RR \to N_\RR$ denote the translation by $v$ and let $Q(\cdot , v)$ denote the linear function $N_{\RR} \to \RR$ given by $x \mapsto Q(x,v)$.

\begin{Lemma}
\label{lem:theta:translation}
The following hold.
\begin{enumerate}
\item
For any $b \in M$,
we have $t_{v}^* \left( \vartheta_b^{(Q, \ell)} \right) 
= \vartheta_b^{(Q, \ell - Q (\cdot , v))} - \frac{1}{2} Q(v, v)$. 
\item
For any $\eta \in T(Q , \ell)$, $t_v^* (\eta) \in T(Q, \ell - Q (\cdot , v))$. 
\end{enumerate}
\end{Lemma}

\Proof
For $x \in N_\RR$, 
we compute
{\allowdisplaybreaks
\begin{align*}
& t_{v}^* \left( \vartheta_b^{(Q, \ell)} \right) (x) 
 = \vartheta_b^{(Q, \ell)}(x + v)\\
& \quad  = \min_{u^\prime \in M^\prime} 
\left(
Q\left(x + v, u^\prime + \lambda_{\RR}^{-1}(b)\right) + \frac{1}{2} Q\left(u^\prime + \lambda_{\RR}^{-1}(b) - r, u^\prime + \lambda_{\RR}^{-1}(b) - r \right)
\right) \\
& \quad  = \min_{u^\prime \in M^\prime} 
\left(
Q\left(x, u^\prime + \lambda_{\RR}^{-1}(b)\right) + \frac{1}{2} Q\left(u^\prime + \lambda_{\RR}^{-1}(b) -r + v, u^\prime + \lambda_{\RR}^{-1}(b)-r+ v\right)\right) \\
& \qquad \qquad - \frac{1}{2} Q(v, v)
\\
& \quad  = \vartheta_b^{(Q, \ell - Q (\cdot , v))}(x) - \frac{1}{2} Q(v, v). 
\end{align*}
}
We obtain (1). 

By (1), $t_v^{\ast} \left( \vartheta_b^{(Q ,\ell)} \right) \in T (Q , \ell - Q (\cdot , v))$.  
Since $\left\{ \vartheta_b^{(Q ,\ell)} \right\}_{b \in M}$ generates $T (Q,\ell)$ as a $\TT$-semimodule (Theorem~\ref{thm:generators0}), (2) follows from (1). Thus the proof is complete.
\QED

The following lemma indicates that if one wants to show that the classes $\varphi^{(Q,\ell)}$ of maps $X \to \TT\PP^{D-1}$ is injective (resp. unimodular), it suffices to show that it is injective (resp. unimodular) under the assumption that $\ell = 0$ (by taking $\ell' := - \ell$). 

\begin{Proposition}
\label{prop:reduction:to:ell:zero}
Let $\ell^\prime\colon M^\prime \to \RR$ be a $\ZZ$-linear map. Then 
if $\varphi^{(Q,\ell)}$ is injective (resp. unimodular), then 
so is $\varphi^{(Q,\ell + \ell^\prime)}$. 
\end{Proposition}

\Proof
Since $Q$ is nondegenerate, there exist $r,v \in N_\RR$ such that $\ell(\cdot) = Q(\cdot, r)$ and $\ell^\prime = - Q(\cdot, v)$.  By Lemma~\ref{lem:theta:translation}~(1), 
$t_{v}^* \left( \vartheta_b^{(Q, \ell)} \right)
= \vartheta_b^{(Q, \ell +\ell^\prime)} - \frac{1}{2} Q(v, v)$. 
Let $\mathfrak{B} \subseteq M$ be a complete system of representatives of $M / \lambda (M')$. We number the elements of $\mathfrak{B}$ to write $\mathfrak{B} = \{ b_1 , \ldots , b_D\}$, where $D = |M / \lambda (M')|$.
Then
\begin{align*}
 \left(\vartheta_{b_1}^{(Q, \ell +\ell^\prime)}, \ldots , \vartheta_{b_D}^{(Q, \ell +\ell^\prime)}\right) 
&= \left(t_{v}^*\left( \vartheta_{b_1}^{(Q, \ell)} \right) + \frac{1}{2} Q(v, v), \ldots , t_{v}^* \left( \vartheta_{b_D}^{(Q, \ell )} \right) + \frac{1}{2} Q(v, v)\right) 
\\
 &= \left(t_{v}^* \left( \vartheta_{b_1}^{(Q, \ell )} \right), \ldots , t_{v}^* \left( \vartheta_{b_D}^{(Q, \ell)} \right)\right) 
 + 
 \frac{1}{2}Q(v,v) (1, \ldots, 1)
 . 
\end{align*}
It follows that
\[
\varphi_{\bar{\vartheta}_{b_1}^{(Q, \ell +\ell^\prime)} , \ldots , \bar{\vartheta}_{b_D}^{(Q, \ell +\ell^\prime)} } 
= \varphi_{\bar{\vartheta}_{b_1}^{(Q, \ell)} , \ldots , \bar{\vartheta}_{b_D}^{(Q, \ell)} } \circ t_{\bar{v}}, 
\]
where $t_{\bar{v}} : X \to X$ is the translation by the image $\bar{v} \in X$ of $v$.
Since $t_{\bar{v}}\colon X \to X$ is unimodular and bijective, we obtain that, if $ \varphi_{\bar{\vartheta}_{b_1}^{(Q, \ell)} , \ldots , \bar{\vartheta}_{b_D}^{(Q, \ell)} }$ is  injective (resp. unimodular), 
then so is $\varphi_{\bar{\vartheta}_{b_1}^{(Q, \ell +\ell^\prime)} , \ldots , \bar{\vartheta}_{b_D}^{(Q, \ell +\ell^\prime)} } $. 
\QED

%%%%%%%%%%%%%
% Preparatory results. %%
%%%%%%%%%%%%%
\setcounter{equation}{0}
\section{Proof of the injectivity}
\label{section:injectivity}

In this section, we prove Theorem~\ref{thm:main:1}~(1). As is noted in Remark~\ref{remark:suffice:injectivity}, it suffices to show that, under the setting of Theorem~\ref{thm:main:1}, $\varphi^{(Q,\ell)}$ is injective if $d_1 \geq 3$. 

We keep the notation in Subsection~\ref{subsection:notation1}. Let $Q$ be  a polarization. Let $\vartheta_b := \vartheta_b^{(Q,0)}$ be the tropical theta function give n by (\ref{eqn:def:theta:b}) for $(Q,0)$.

\subsection{Key property of tropical theta functions}
Our goal in this subsection is to prove Proposition~\ref{thm:important:for:main}, 
which will be the key for the injectivity.

\begin{Proposition}
\label{thm:important:for:main}
Let $Q$ be a polarization on $X = N_{\RR} /M'$ of type
$(d_1, \ldots, d_n)$ and let $z \in N_{\RR}$.
Assume that $d_1 \geq 2$. 
and that for any $b \in M$,
\begin{equation}
\label{eqn:thm:important:for:main:1}
 \vartheta_b(x + z) - \vartheta_0(x + z)
 =  \vartheta_b(x) - \vartheta_0(x)
\end{equation}
holds as functions on $x \in N_\RR$.
Then $z \in M^\prime$. 
\end{Proposition}

Note that the condition for $d_1$ in Proposition~\ref{thm:important:for:main} 
is $d_1 \geq 2$ and not
$d_1 \geq 3$ as in Theorem~\ref{thm:main:1}~(1).

\begin{Lemma}
\label{lem:for:thm:important:for:main:1}
Under setting and the assumption in Proposition~\ref{thm:important:for:main},  
for any $m \in \ZZ_{\geq 1}$,  
\begin{equation}
\label{eqn:lem:for:thm:important:for:main:1}
 \vartheta_0\left(d_1^m(x + z)\right) 
 - \vartheta_0\left(d_1^m x\right)
 =  d_1^{2m} \left(\vartheta_0(x+z) - \vartheta_0(x)\right)
\end{equation}
holds as a function on $x \in N_\RR$.
\end{Lemma}

\Proof
It follows from the definition of $d_1, \ldots ,  d_n$ that $\lambda (M') \subseteq d_1 M$, and thus $\frac{1}{d_1} \lambda (M') \subseteq M$. There exists a complete system of representatives $\mathfrak{B}_1$ of $\frac{1}{d_1} \lambda (M') / \lambda (M')$. Then $d_1 \mathfrak{B}_1$ is a complete system of representatives of $\lambda (M')/ d_1 \lambda (M')$, and hence $ \lambda_{\RR}^{-1} (d_1\mathfrak{B}_1) \subseteq M'$ is a complete system of representatives of $M'/dM'$. It follows that the map
\begin{align} \label{align:bijection}
M' \times \mathfrak{B}_1 \to M' ; \quad (u' , b) \mapsto d_1 u' +  \lambda_{\RR}^{-1} (d_1b)
\end{align}
is bijective.

We set 
$\eta (x) \colonequals \min_{b \in \mathfrak{B}_1} \vartheta_b(x)$. 
Then for any $x \in N_\RR$, noting the bijection in (\ref{align:bijection}), we compute
{\allowdisplaybreaks
\begin{align}
\label{eqn:lem:for:thm:important:for:main:3}
\begin{split}
d_1^2 \eta(x) 
& = \min_{b \in \mathfrak{B}_1} 
\min_{u^\prime \in M^\prime} 
\left(
d_1^2Q\!\left( x,\, u^\prime + \lambda_\RR^{-1}(b)\right) + \frac{1}{2} d_1^2 Q\!\left(u^\prime + \lambda_\RR^{-1}(b),\, u^\prime + \lambda_\RR^{-1}(b)\right)
\right) \\
&=
\min_{(u' , b) \in M' \times \mathfrak{B}_1} 
\left(
Q\!\left(d_1 x,\, d_1u^\prime + \lambda_{\RR}^{-1}(d_1b)\right) + \frac{1}{2} Q\!\left(d_1u^\prime + \lambda_{\RR}^{-1}(d_1b),\, d_1u^\prime + \lambda_{\RR}^{-1}(d_1b)\right)
\right) \\
& =
\min_{\tilde{u}^\prime \in M^\prime} 
\left(
Q\!\left(d_1 x,\, \tilde{u}^\prime\right) + \frac{1}{2} Q\!\left(\tilde{u}^\prime, \, 
\tilde{u}^\prime\right)
\right) 
=  \vartheta_0(d_1 x). 
\end{split}
\end{align}}
By \eqref{eqn:thm:important:for:main:1}, 
we have, for any $x \in N_\RR$,  
\begin{align}
\label{eqn:lem:for:thm:important:for:main:2}
\begin{split}
\eta (x+z) - \vartheta_0(x+z)  
%\\
%\notag
& = \min_{b \in \mathfrak{B}_1} \vartheta_b(x+z)  - \vartheta_0(x+z)
= \min_{b \in \mathfrak{B}_1} (\vartheta_b(x+z)  - \vartheta_0(x+z)) \\
%\notag
& = \min_{b \in \mathfrak{B}_1} (\vartheta_b(x)  - \vartheta_0(x))
= \min_{b \in \mathfrak{B}_1} \vartheta_b(x)
  - \vartheta_0(x)
= \eta(x) - \vartheta_0(x).   
\end{split}
\end{align}
Combining 
\eqref{eqn:lem:for:thm:important:for:main:2} and 
\eqref{eqn:lem:for:thm:important:for:main:3}, we get for any $x \in N_\RR$, 
\[
\frac{1}{d_1^2} \vartheta_0(d_1 (x+z)) - \vartheta_0(x + z)
= 
\frac{1}{d_1^2} \vartheta_0(d_1 x)- \vartheta_0(x) 
\]
for any $x \in N_\RR$, which is equivalent to the $m = 1$ case of \eqref{eqn:lem:for:thm:important:for:main:1}.

Next, we assume \eqref{eqn:lem:for:thm:important:for:main:1} for $m \geq 1$, and we will prove it for $m+1$. A telescoping argument gives 
{\allowdisplaybreaks
\begin{align*}
& \vartheta_0\left(d_1^{m+1}(x + z)\right) 
 - \vartheta_0\left(d_1^{m+1} x\right) \\
& \quad = 
\sum_{k=1}^{d_1} \left(
\vartheta_0\left(d_1^m\left((d_1 x + (k-1) z) + z\right)\right)
-  \vartheta_0\left(d_1^m\left((d_1 x + (k-2) z) + z\right)\right)
\right)
\\
& \quad = 
d_1^{2m}  
\sum_{k=1}^{d_1} \left(
\vartheta_0\left((d_1 x + (k-1) z) + z\right)
-  \vartheta_0\left((d_1 x + (k-2) z) + z\right)
\right)
\\
& \quad = 
 d_1^{2m} \left(\vartheta_0(d_1(x+z)) - \vartheta_0(d_1x)\right)
=  
d_1^{2m+2} \left(\vartheta_0(x+z) - \vartheta_0(x)\right).
\end{align*}
}
This completes the proof. 
\QED

\begin{Lemma}
\label{lem:for:thm:important:for:main:2}
We set  
\[
\mathfrak{F}
\colonequals \left\{
x \in N_\RR \;\left\vert\; 
\text{
$Q(x + u^\prime, x + u^\prime) \geq Q(x, x)$ 
for any $u^\prime \in M^\prime$
}
\right. 
\right\}. 
\]
Then the following hold: 
\begin{enumerate}
\item[(1)]
$0 \in \mathfrak{F}$;
\item[(2)]
$N_{\RR} = \mathfrak{F} + M^\prime$;
\item[(3)]
$\vartheta_0(x) = 0$ for any $x \in \mathfrak{F}$, where $\vartheta_0$ is the one defined by (\ref{eqn:def:theta:b}) for $(Q,0)$. 
\end{enumerate}
\end{Lemma}

\Proof
Since $Q$ is positive-definite, we have 
$Q(u^\prime, u^\prime) \geq 0 = Q(0, 0)$ for any $u^\prime \in M^\prime$. 
Thus (1) holds.

To prove
(2), we take any $y \in N_\RR$. Since $Q$ is positive-definite, 
there is an element $u^\prime_0 \in M^\prime$ such that 
$Q(y + u^\prime_0, y + u^\prime_0) 
=  
  \min_{u^\prime \in M^\prime} 
  Q(y + u^\prime, y + u^\prime)$.  
We set $x = y + u^\prime_0$. Then $x \in \mathfrak{F}$, and 
we obtain $y = x - u^\prime_0 \in \mathfrak{F} + M^\prime$. 

We prove
(3). For $x \in \mathfrak{F}$, we have 
\begin{equation}
\label{eqn:vartheta:0}
\vartheta_0(x) = 
\min_{u^\prime \in M^\prime} 
\left(
Q(x, u^\prime) + \frac{1}{2} Q(u^\prime, u^\prime)
\right)
= 
\min_{u^\prime \in M^\prime} 
%\left(
\frac{1}{2} Q(x + u^\prime, x + u^\prime)
%\right\}
- \frac{1}{2}  Q(x,  x). 
\end{equation}
By the definition of $\mathfrak{F}$, the minimum of the right-hand side 
is attained at $u^\prime = 0$, and we get $\vartheta_0(x) = 
\frac{1}{2}  Q(x,  x) - \frac{1}{2}  Q(x,  x) = 0$. 
\QED

{\sl Proof of Proposition~\ref{thm:important:for:main}.} \quad 
By Lemma~\ref{lem:for:thm:important:for:main:2}~(2), 
there is $z^\prime \in \mathfrak{F}$ such that $z = z^\prime + u^\prime$ for 
some $u^\prime \in M^\prime$. By quasi-periodicity \eqref{eqn:quasi:periodicity:0}, we have 
$
\vartheta_b(x + z^\prime) - \vartheta_0(x + z^\prime)
 =  \vartheta_b(x+z) - \vartheta_0(x+z). 
$
Thus, replacing $z$ by $z - u^\prime$, we may assume that 
$z \in \mathfrak{F}$. 

We claim that $z = 0$. 
Indeed, by Lemma~\ref{lem:for:thm:important:for:main:2}, 
we have $0 \in  \mathfrak{F}$, $\vartheta_0(0) = 0$, and 
$\vartheta_0(z) = 0$. 
Substituting $x = 0$ in \eqref{eqn:lem:for:thm:important:for:main:1}, 
we get 
$
 \vartheta_0\left(d_1^m z \right)  = 0 
$
for all $m \geq 1$ by Lemma~\ref{lem:for:thm:important:for:main:1}. We endow $N_\RR$ with an Euclidean norm $\Vert \ndot \Vert$.  If $z \neq 0$, then as $m$ tends to $+\infty$ we have $\Vert d_1^m z  \Vert \to +\infty$. Since $N_\RR/M^\prime$ is compact,  
$\min_{u^\prime \in M^\prime} 
\frac{1}{2} Q(x + u^\prime, x+u^\prime)$ is bounded as above as a function on $x \in N_\RR$. 
Then, by \eqref{eqn:vartheta:0}, we see that 
$\vartheta_0\left(d_1^m z \right) \to -\infty$. 
This is a contradiction. Thus $z = 0$, and we are done. 
\QED

%%%%%%%%%%%%%%%%%%%%%%%%%%%%%%
% Injectivity and proof of Theorem~\ref{thm:main:1} %%
%%%%%%%%%%%%%%%%%%%%%%%%%%%%%%
\subsection{Proof of Theorem~\ref{thm:main:1}~(1)} 

We prove Theorem~\ref{thm:main:1}~(1).
We keep the setting that we fixed at the beginning of this section. 
By Proposition~\ref{prop:reduction:to:ell:zero} and Remark~\ref{remark:suffice:injectivity}, it suffices to show that $\varphi^{(Q, 0)}$ is injective.
We take a complete system of representatives $\mathfrak{B} \subseteq M$ with $0 \in \mathfrak{B}$, and 
we write $\mathfrak{B} = \{b_1, \ldots, b_D\}$. We take $b_1 = 0$, so $\vartheta_{b_1} = \vartheta_0$. 

To prove the injectivity of $\varphi^{(Q,0)}$,
suppose that
\[
\left(\vartheta_{b_1}(y): \cdots: \vartheta_{b_D}(y)\right)
= \left(\vartheta_{b_1}(y^\prime): \cdots: \vartheta_{b_D}(y^\prime)\right). 
\]
Then
\begin{equation}
\label{eqn:equality}
\vartheta_b(y) - \vartheta_0(y) 
= \vartheta_b(y^\prime) - \vartheta_0(y^\prime) 
\end{equation}
for any $b \in M'$. We want to show that $\bar{y} = \bar{y^\prime}$ in $X = N_{\RR}/M'$, i.e.,  
$y-y^\prime \in M^\prime$. 
Note that for any $\vartheta, \vartheta^\prime \in T(Q,0)$, we have 
\begin{equation}
\label{eqn:equality:2}
\vartheta(y) - \vartheta^\prime(y) 
= \vartheta(y^\prime) - \vartheta^\prime(y^\prime);  
\end{equation}
see (\ref{eqn:TQell}) for the notation $T(Q,0)$.
Indeed, since $\{\vartheta_b\}_{b \in \mathfrak{B}}$ generates 
$T(Q,0)$ (Theorem~\ref{thm:generators0}), it follows from \eqref{eqn:equality} that 
$\vartheta(y) - \vartheta_0(y) 
= \vartheta(y^\prime) - \vartheta_0(y^\prime)$  for any $\vartheta \in T(Q,0)$. Thus equality (\ref{eqn:equality:2}) holds.

For $1 \leq i \leq n$, we write $d_i = \delta_i d_1$ with $\delta_i \in \ZZ_{\geq 1}$. We set $Q_1 = \frac{1}{d_1} Q$. Then 
$Q_1 \colon N_\RR \times N_\RR \to \RR$ satisfies $Q_1 (M^\prime\times N) \subseteq \ZZ$ and is a polarization on $X$ of type $(1, \delta_2, \ldots, \delta_n)$. Since $d_1 \geq 3$ by assumption, we have $d_1 - 2 \geq 1$ and 
$2 Q_1 + (d_1 -2) Q_1 = d_1 Q_1 = Q$. 

For any $\eta, \eta^\prime \in T( 2Q_1,0)$, we prove that 
\begin{align} \label{align:translate:eta}
t_{v}^*\eta(y) - t_{v}^*\eta^\prime(y) 
= 
t_{v}^*\eta(y^\prime) - t_{v}^*\eta^\prime(y^\prime) 
\end{align}
for any $v \in N_\RR$. Take any $\eta, \eta^\prime \in T(2 Q_1 , 0)$ and any $v \in N_{\RR}$. Since $d_1 - 2 \geq 1$, the $\RR$-bilinear form $(d_1 -2) Q_1$ is a polarization, and hence there exists a nontrivial 
$\zeta \in T( (d_1 -2) Q_1 , 0)$. Then by Lemma~\ref{lem:theta:translation}~(2), we have 
\begin{align*}
&t_{v}^*(\eta) , t_{v}^*(\eta') \in T \left( 2Q_1 , -2Q_1 (\cdot , v)  \right)
,
\\
&t_{-2v/(d_1 -2)}^*(\zeta) \in T \left( (d_1-2)Q_1 , - (d_1-2) Q_1 \left( \cdot , \frac{-2v}{d_1 -2} \right) \right) = T ((d_1 - 2) Q_1 , 2 Q_1 (\cdot , v))
,
\end{align*}
and thus 
$t_{v}^*(\eta) +   t_{-2v/(d_1 -2)}^*(\zeta), \; 
t_{v}^*(\eta^\prime) +   t_{-2v/(d_1 -2)}^*(\zeta) \in T (Q,0)$. 
By \eqref{eqn:equality:2}, it follows that 
\begin{multline*}
\left(t_{v}^*\eta(y) +   t_{-2v/(d_1 -2)}^*\zeta(y)\right) - 
\left(t_{v}^*\eta^\prime(y) +   t_{-2v/(d_1 -2)}^*\zeta(y)\right)  \\
 = 
\left(t_{v}^*\eta(y^\prime) +   t_{-2v/(d_1 -2)}^*\zeta(y^\prime)\right) - 
\left(t_{v}^*\eta^\prime(y^\prime) +   t_{-2v/(d_1 -2)}^*\zeta(y^\prime)\right),  
\end{multline*}
which shows (\ref{align:translate:eta}).

We set $z \colonequals y - y^\prime \in N_\RR$. 
Then for any $\eta, \eta^\prime \in  T(2 Q_1 , 0)$ and for any $x \in N_\RR$, since (\ref{align:translate:eta}) holds, we have 
\begin{align*}
\eta(x+z) -  \eta^\prime(x+z)
& = \eta(y + x-y^\prime) -  \eta^\prime(y + x-y^\prime)
= t_{x-y^\prime}^* \eta(y) -  t_{x-y^\prime}^* \eta^\prime(y) \\
& = t_{x-y^\prime}^* \eta(y^\prime) -  t_{x-y^\prime}^* \eta^\prime(y^\prime) 
= \eta(x) -  \eta^\prime(x). 
\end{align*}
Since $2Q_1$ has polarization type $(2, 2 \delta_2, \ldots, 2 \delta_n)$, 
we apply Proposition~\ref{thm:important:for:main} to $2Q_1$ to conclude that $y - y' = z \in M^\prime$.
Thus the proof is complete.

%%%%%%%%%%%%%
% Unimodularity %%
%%%%%%%%%%%%%
\setcounter{equation}{0}
\section{Proof of the unimodularity} \label{section:Pf:unimodular}
In this section, we prove Theorem~\ref{thm:main:1}~(2). We keep the notation in Subsection~\ref{subsection:notation1}. Let $Q$ be  a polarization on $N_{\RR} / M'$. Let $(d_1 , \ldots , d_n)$ be the type of the polarization $Q$. For any $b \in M$, let $\vartheta_b := \vartheta_b^{(Q,0)}$ be the tropical theta function on $N_{\RR}$ give by (\ref{eqn:def:theta:b}) for $(Q,0)$. Recall that, by Proposition~\ref{prop:reduction:to:ell:zero}, it suffices to show that $\varphi^{(Q,0)} $ is unimodular.

\subsection{Main proposition}

In this subsection, we state the main proposition, Proposition~\ref{prop:main:unimodular:new}, of this section, which will imply  Theorem~\ref{thm:main:1}~(2). 
A \emph{polyhedral tiling} on a polyhedron $P$ of dimension $n$ means a locally finite family $\Sigma$ of polyhedra of dimension $n$ in $P$ with the following properties: 
%set $\Sigma^{(n)} := \{ \sigma' \in \Sigma \mid  \dim (\sigma') = n \}$; then 
for any $\sigma_1 , \sigma_2 \in \Sigma^{(n)}$ with $\sigma_1 \neq \sigma_2$, $\relin (\sigma_1) \cap \relin (\sigma_2) = \emptyset$ holds, and 
$\bigcup_{\sigma \in \Sigma^{(n)}} \relin (\sigma)$ is a dense open subset of $P$. We remark that for a polyhedral tiling $\Sigma$ on a polyhedron and for any $\sigma \in \Sigma$, $\relin (\sigma)$ is an open subset of the polyhedron.

For a function $g$ on a set $A$, let $\argmin_{a \in A} g(a)$ denote  
the subset of $A$ consisting of those $a  \in A$ that attain $\min_{a \in A} g(a)$.

\begin{Proposition} \label{prop:main:unimodular:new}
Assume that $d_1 \geq 2 (n-1)!$. Then there exists a polyhedral tiling $\Sigma$ of $N_{\RR}$ such that for any $\sigma \in \Sigma$, there exist $b_0 , b_1, \ldots , b_n \in M$ with the following properties: 
\begin{enumerate}
\item[(i)]
$b_1 - b_0 , \ldots , b_n - b_0$ is a basis of $M$;
\item[(ii)]
there exists a $\widetilde{u}' \in M'$ such that for any $x \in \sigma$ and for any $j=0,\ldots , n$, we have 
%$\tilde{u}' \in \argmin \vartheta_{b_{j}}(z)$.
\[
\tilde{u}' \in
\underset{u' \in M'}{\argmin} 
\left( 
Q (x, \lambda_{\RR}^{-1} (b_j) + u') +
\frac{1}{2} Q (u' + \lambda_{\RR}^{-1} (b_j) , u' + \lambda_{\RR}^{-1} (b_j))
\right).
\]
\end{enumerate}
\end{Proposition}

In the rest of this subsection, we prove that Proposition~\ref{prop:main:unimodular:new} implies Theorem~\ref{thm:main:1}~(2). 

\begin{Lemma} \label{lemma:unimodular:dense:open}
Let $\varphi : N_{\RR} \to \RR^m$ be a piecewise $\ZZ$-affine map. Let $U$ be a dense open subset of $N_\RR$. Suppose that for any $p \in U$, there exists an open neighborhood $V \subseteq U$ of $p$ such that $\rest{\varphi}{V} : V \to \RR^m$ is a unimodular $\ZZ$-affine map. Then $\varphi$ is unimodular.
\end{Lemma}

\Proof
Since $\varphi$ is piecewise $\ZZ$-affine, there exists a rational polyhedral decomposition $\mathcal{C}$ of $N_{\RR}$ such that $\rest{\varphi}{\sigma}$ is $\ZZ$-affine for any $\sigma \in \mathcal{C}$. For any $\sigma \in \mathcal{C}$, it follows from the definition of unimodularity that if $\rest{\varphi}{\relin (\sigma)}$ is unimodular, then $\rest{\varphi}{\sigma}$ is unimodular and hence $\rest{\varphi}{\tau}$ is unimodular for any face $\tau$ of $\sigma$. Thus it suffices to show that $\rest{\varphi}{\relin (\sigma)}$ is unimodular for any $\sigma \in \mathcal{C}^{(n)}$, where $\mathcal{C}^{(n)} := \{ \sigma \in \mathcal{C} \mid \dim (\sigma) = n \}$.

We take any $\sigma \in \mathcal{C}^{(n)}$. We take $U$ as in the lemma. Since $\relin (\sigma)$ is open in $N_{\RR}$ and $U$ is dense in $N_{\RR}$, there exists a point $p \in \relin (\sigma) \cap U$. By the assumption, there exists an open neighborhood $V$ of $p$ in $\relin (\sigma) \cap U$ such that $\rest{\varphi}{V}$ is a unimodular $\ZZ$-affine map. Since $\rest{\varphi}{\relin (\sigma)}$ is $\ZZ$-affine, it follows that $\rest{\varphi}{\relin (\sigma)}$ is unimodular. This completes the proof of the lemma.
\QED

Assume that Proposition~\ref{prop:main:unimodular:new} holds. Noting Remark~\ref{remrak:lift:UM:FE}, we only have to show that the map $\varphi^{(Q,0)} \circ \pi : N_{\RR} \to \TT\PP^{D-1}$ is unimodular, where $\pi : N_{\RR} \to N_{\RR} / M'$ is the canonical surjection and $D:= d_1 \cdots d_n$.
We take a polyhedral tiling $\Sigma$ as in Proposition~\ref{prop:main:unimodular:new}. By Lemma~\ref{lemma:unimodular:dense:open}, it suffices to show that the map $\psi$ is unimodular over $\bigcup_{\sigma \in \Sigma} \relin (\sigma)$.
We take any $\sigma \in \Sigma^{(n)}$, and we take $b_0 , b_1, \ldots , b_n$ as in the proposition. It suffices to show that the map $N_{\RR} \to \RR^n$ given by
\[
x \mapsto ( \vartheta_{b_1} (x) - \vartheta_{b_0} (x) , \ldots , \vartheta_{b_n} (x) - \vartheta_{b_0} (x) )
\]
is already unimodular on $\relin (\sigma)$; see Remarks~\ref{remark:unimdular:part}~(2).

We take $\widetilde{u}' \in M'$ as in Proposition~\ref{prop:main:unimodular:new}~(ii).
For any $b \in M$, we note that, as functions on $x \in N_{\RR}$,
\[
Q (x, \lambda_{\RR}^{-1} (b) + u') +
\frac{1}{2} Q (u' + \lambda_{\RR}^{-1} (b) , u' + \lambda_{\RR}^{-1} (b))
= \langle b + \lambda (u') , x \rangle + C_{u',b},
\]
where $C_{u' b}$ is a constant on $x$, and thus
\[
\vartheta_{b_j} (x) - \vartheta_{b_0} (x) = \langle b_j - b_0 , x \rangle + C_{j}' 
\]
for some constant $C_j'$. Since, by property (i) in the proposition, $b_1 - b_0 , \ldots , b_n-b_0$ is a basis of $M$, this shows that the map \[
( \vartheta_{b_1} (x) - \vartheta_{b_0} (x) , \ldots , \vartheta_{b_n} (x) - \vartheta_{b_0} (x) )
\]
is unimodular as a map on $x \in \relin (\sigma)$. Thus Proposition~\ref{prop:main:unimodular:new} implies Theorem~\ref{thm:main:1}~(2).

\subsection{Paraphrase of Proposition~\ref{prop:main:unimodular:new}}

In this subsection, we paraphrase Proposition~\ref{prop:main:unimodular:new} by taking bases of $M$ and $M'$ and using the canonical inner product on $\RR^n$.
% Proposition~\ref{prop:main:unimodular:new} to another proposition.

First, note that 
\[
\vartheta_b (x) =
\frac{1}{2} \min_{u' \in M'}
%\left\{
Q(x+ \lambda_{\RR}^{-1} (b) + u', x+ \lambda_{\RR}^{-1} (b) + u')
%\right\}
- \frac{1}{2} Q(x,x).
\]
Indeed, for any $b \in M$ and $u' \in M'$, we have, as functions on $x \in N_{\RR}$,
\[
Q (x, \lambda_{\RR}^{-1} (b) + u') +
\frac{1}{2} Q (u' + \lambda_{\RR}^{-1} (b) , u' + \lambda_{\RR}^{-1} (b))
= \frac{1}{2} Q (x+ \lambda_{\RR}^{-1} (b) + u', x+ \lambda_{\RR}^{-1} (b) + u') - \frac{1}{2} Q(x,x).
\]
Thus
\begin{multline} \label{multline:for:prop:main:unimodular:new}
\underset{u' \in M'}{\argmin} \  
\left( 
Q (x, \lambda_{\RR}^{-1} (b) + u') +
\frac{1}{2} Q (u' + \lambda_{\RR}^{-1} (b) , u' + \lambda_{\RR}^{-1} (b))
\right)
\\
=
\underset{u' \in M'}{\argmin} \  
%\argmin_{u' \in M'}
Q(x+ \lambda_{\RR}^{-1} (b) + u', x+ \lambda_{\RR}^{-1} (b) + u').
\end{multline}

In view of property (ii) in Proposition~\ref{prop:main:unimodular:new}, 
for fixed $b \in M$ and $x \in N_{\RR}$, we then are interested in
\begin{align*} 
%\label{align:argmin:Q}
\underset{u' \in M'}{\argmin} \  
%\argmin_{u' \in M'}
Q(x+ \lambda_{\RR}^{-1} (b) + u', x+ \lambda_{\RR}^{-1} (b) + u').
\end{align*}
% minimum on the right-hand side. 
To look into this,
%(\ref{align:argmin:Q})
we take suitable bases of $M$ and $M'$.
% to describe $Q(x+ \lambda_{\RR}^{-1} (b) + u', x+ \lambda_{\RR}^{-1} (b) + u')$.  
%Let $(d_1 , \ldots , d_n)$ be the type of the polarization $Q$. 
Let $\lambda : M' \to M$ be the homomorphism corresponding to $Q$ by Lemma~\ref{lem:polarization:equiv}~(1). There exist free  $\ZZ$-bases $\mathsf{f}^\prime_1, \ldots, \mathsf{f}^\prime_n$ and $\mathsf{e}_1, \ldots, \mathsf{e}_n$ of
 $M^\prime$ and $M$, respectively, such that $\lambda(\mathsf{f}^\prime_1) = d_1 \mathsf{e}_1, \ldots, \lambda(\mathsf{f}^\prime_n) = d_n \mathsf{e}_n$. Since $Q$ is symmetric and positive-definite, the $n \times n$ matrix
%Let $Q^{\rm mtrx}$ be the $n \times n$ positive symmetric real matrix that represents 
%$Q$ with respect to $\mathsf{f}^\prime_1, \ldots, \mathsf{f}^\prime_n$, i.e., 
%\[
%Q^{\rm mtrx} := 
$
\begin{pmatrix}
Q(\mathsf{f}^\prime_i, \mathsf{f}^\prime_j)
\end{pmatrix}_{i, j}
$
%= ( [\mathsf{f}^\prime_i, d \mathsf{e}_j])_{1 \leq i, j \leq n}.
%\]
is positive and symmetric.
%Recall that $d_1 , \ldots , d_n \in \ZZ_{\geq 1}$ and that $d_i \mid d_{i+1}$ for all $i=1 , \ldots , n-1$.
%Since $Q^{\rm mtrx}$ is positive and symmetric, 
It follows that there exists a $P = \begin{pmatrix}
 \mathsf{p}_1 & \ldots & \mathsf{p}_n
\end{pmatrix}
  \in \GL_n(\RR)$ 
such that 
\begin{align} \label{eqn:def:pi}
\begin{pmatrix}
Q(\mathsf{f}^\prime_i, \mathsf{f}^\prime_j)
\end{pmatrix}_{i, j} = 
%\begin{pmatrix}
%\mathsf{p}_1 & \ldots & \mathsf{p}_n
%\end{pmatrix}^{{\sf T}}
%\begin{pmatrix}
%\mathsf{p}_1 & \ldots & \mathsf{p}_n
%\end{pmatrix}
P^{{\sf T}} P,
\end{align}
where $\mathsf{p}_j$ is the $j$th column vector of $P$ and $P^{\mathsf{T}}$ is the transpose of $P$. 

% and $P^{\mathsf{T}}$ is the transpose of $P$.
%for a matrix $A$, $A^{\mathsf{T}}$ denotes the transpose of $A$.
%We set 
%\begin{align}
%\label{eqn:def:pi}
%\mathsf{p}_i \colonequals 
%P\, (0, \ldots, 0, 1, 0 \ldots, 0)^{{\sf T}} \in \RR^n \qquad (i = 1, 2, \ldots, n), 
%\end{align}
%where $(0, \ldots, 0, 1, 0 \ldots, 0)^{{\sf T}}$ is the $i$-th standard column vector. 
%Since $P \in \GL_n(\RR)$, the vectors $\mathsf{p}_1, \ldots , \mathsf{p}_n$ form an $\RR$-basis of $\RR^n$.

Let $\phi : \ZZ^n \to M$ and $\psi : \ZZ^n \to M'$ be isomorphisms
arising from the bases $\mathsf{e}_1 , \ldots , \mathsf{e}_n$ of $M$ and $\mathsf{f}^{\prime}_1 , \ldots , \mathsf{f}^{\prime}_n$ of $M'$, repectively, and let $\psi_{\RR} : \RR^n \to N_{\RR}$ is the $\RR$-linear extension of $\psi$. 
Let us fix any $x \in N_{\RR}$ and $b \in M'$. Then there exist unique
$\ell = (\ell_1 , \ldots , \ell_n) \in \ZZ^n$ and $y = (y_1 , \ldots , y_n) \in \RR^n$ such that
%.
%For any $\ell = (\ell_1 , \ldots , \ell_n) \in \ZZ^n$, 
%We write
\[
b := \phi (\ell) = \sum_{i=1}^n \ell_i \mathsf{e}_i , \quad
x := \psi_{\RR} (y) = \sum_{i=1}^n y_i \mathsf{f}_i'.
\]
%For $u' \in M'$, we write
Let $[\ndot, \ndot]$ denote the standard inner product on $\RR^n$. We naturally regard $\mathsf{p}_1, \ldots , \mathsf{p}_n \in \RR^n$. Then
we have
$Q (x,x) = \left[ \sum_{i=1}^n y_i \mathsf{p}_i , \sum_{i=1}^n y_i \mathsf{p}_i \right]$. Further, 
if $u' = \sum_{i=1}^n a_i \mathsf{f}_i^\prime \in M'$ with $a = (a_1 , \ldots , a_n) \in \ZZ^{n}$,  
then 
%with the notation so far, we compute
%$Q (x,x) = \left[ \sum_{i=1}^n y_i \mathsf{p}_i , \sum_{j=1}^n y_j \mathsf{p}_j \right]$ and
\[
Q (x+ u' + \lambda_{\RR}^{-1} (b_{\ell}) , x+u' + \lambda_{\RR}^{-1} (b_{\ell}))
=
\left[
\sum_{i=1}^n \left(a_i + \frac{\ell_i}{d_i} + y_i\right) \mathsf{p}_i, \, 
\sum_{i=1}^n \left(a_i + \frac{\ell_i}{d_i} + y_i\right) \mathsf{p}_i
\right]
.
\]
Thus we have
%Then, from the expression \eqref{eqn:theta:b:bis}, we have an expression 
%{\allowdisplaybreaks
\begin{align}
\label{eqn:theta:b:bis:contd}
& \vartheta_{b_\ell}(x) 
%\\
%\notag
%& \qquad = 
%\min_{a \in \ZZ^n} 
%\left\{
%(x_1, \ldots, x_n) Q^{\rm mtrx} 
%\begin{pmatrix} a_1 + \frac{\ell_1}{d_1} \\ \vdots \\  a_n + \frac{\ell_n}{d_n}
%\end{pmatrix}
% + 
% \frac{1}{2} 
% \left(a_1 + \frac{\ell_1}{d_1}, \ldots, a_n + \frac{\ell_n}{d_n} \right) 
% Q^{\rm mtrx} 
%\begin{pmatrix} a_1 + \frac{\ell_1}{d_1} \\ \vdots \\  a_n + \frac{\ell_n}{d_n}
%\end{pmatrix}
%\right\}
%\\
%\notag
%& \qquad  =  
% \min_{a \in \ZZ^n} 
%\left\{
%\left[
%y, \, 
%\sum_{i=1}^n \left(a_i + \frac{\ell_i}{d_i}\right) \mathsf{p}_i
%\right]
%+  
%\frac{1}{2} 
%\left[
%\sum_{i=1}^n \left(a_i + \frac{\ell_i}{d_i}\right) \mathsf{p}_i, \, 
%\sum_{i=1}^n \left(a_i + \frac{\ell_i}{d_i}\right) \mathsf{p}_i
%\right]
%\right\}
%\\
%\notag
%& 
%\qquad  
=  
\frac{1}{2}
 \min_{a \in \ZZ^n} 
%\left\{
\left[
\sum_{i=1}^n \left(a_i + \frac{\ell_i}{d_i} + y_i\right) \mathsf{p}_i, \, 
\sum_{i=1}^n \left(a_i + \frac{\ell_i}{d_i} + y_i\right) \mathsf{p}_i
\right]
%\right\}
- \frac{1}{2} Q(x,x),
\end{align}
%}
and the isomorphism $\psi : \ZZ^n \to M'$ restricts to a bijection
\begin{multline} \label{multline:bijection}
\underset{a \in \ZZ^n} {\argmin} \ 
%\argmin_{a \in \ZZ^n} 
%\left\{
%\frac{1}{2} 
\left[
\sum_{i=1}^n \left(a_i + \frac{\ell_i}{d_i} + y_i\right) \mathsf{p}_i, \, 
\sum_{i=1}^n \left(a_i + \frac{\ell_i}{d_i} + y_i\right) \mathsf{p}_i
\right]
%\right\}
\\
\longrightarrow
\underset{u' \in M'}{\argmin} \  
%\argmin_{u' \in M'}
Q(x+ \lambda_{\RR}^{-1} (b) + u', x+ \lambda_{\RR}^{-1} (b) + u').
\end{multline}
It follows that Proposition~\ref{prop:main:unimodular:new} amounts to the following proposition.

\begin{Proposition}
\label{cor:voronoi:variant:2} 
Assume that $d_1 \geq 2 (n-1)!$. Then there exists a polyhedral tiling $\widetilde{\Sigma}$ of $\RR^n$ such that for any $\widetilde{\sigma} \in \widetilde{\Sigma}$, there exist $\ell^{(0)} , \ell^{(1)}, \ldots, \ell^{(n)} \in \ZZ^n$ that satisfy the following conditions:
\begin{enumerate}
\item[(i)]
$\ell^{(1)} - \ell^{(0)}, \ldots, \ell^{(n)} - \ell^{(0)}$ form a basis of the free $\ZZ$-module $\ZZ^n$; 
\item[(ii)]
%for any $z \in \widetilde{\sigma}$ and 
%for any $j = 0, 1, \ldots, n$, 
there exists an $\widetilde{a} \in  \ZZ^n$ such that for any $(y_1 , \ldots , y_n) \in \RR^n$ with $\sum_{i=1}^n y_i \mathsf{p}_i \in \widetilde{\sigma}$ and for any $j = 0, 1, \ldots, n$, we have 
\begin{align} \label{align:cor:voronoi:variant:2}
\widetilde{a} \in 
\underset{a \in \ZZ^n} {\argmin} \ 
%\argmin_{a \in \ZZ^n} 
%\left\{
%\frac{1}{2} 
\left[
\sum_{i=1}^n \left(a_i + \frac{\ell_i^{(j)}}{d_i} + y_i\right) \mathsf{p}_i, \, 
\sum_{i=1}^n \left(a_i + \frac{\ell_i^{(j)}}{d_i} + y_i\right) \mathsf{p}_i
\right]
,
\end{align}
where $\ell^{(j)} = (\ell_1^{(j)} , \ldots , \ell_n^{(j)})$.
%$\widetilde{a} \in 
%\argmin \xi_{\ell}(y)$  
%for any $z \in \widetilde{\sigma}$ and 
%for any $j = 0, 1, \ldots, n$. 
\end{enumerate}
\end{Proposition}

Indeed, 
we take a polyhedral titling $\widetilde{\Sigma}$ of $\RR^n$ as in Proposition~\ref{cor:voronoi:variant:2}. Let $\rho : \RR^n \to \RR^n$ be the isomorphism given by $(y_1 , \ldots , y_n) \mapsto \sum_{i=1}^n y_i \mathsf{p}_i$. Set    $\Sigma := \{ \psi_{\RR} (\rho^{-1}(\widetilde{\sigma})) \mid \widetilde{\sigma} \in \widetilde{\Sigma} \}$. Fix any $\sigma \in \Sigma$. For $\widetilde{\sigma} := \rho (\psi_{\RR}^{-1} (\sigma)) \in \widetilde{\Sigma}$, we take $\ell^{(0)} , \ldots , \ell^{(n)}$ as in Proposition~\ref{cor:voronoi:variant:2} and set $b_i := \phi (\ell^{(i)})$ for $i = 1, \ldots, n$. Then by Proposition~\ref{cor:voronoi:variant:2}~(i), $b_0 , \ldots , b_n$ satisfy Proposition~\ref{prop:main:unimodular:new}~(i). We take $\widetilde{a} \in \ZZ^n$ as in Proposition~\ref{cor:voronoi:variant:2}~(ii) and set $\widetilde{u}' := \psi (\widetilde{a})$. We take any $x \in N_{\RR}$ with $x \in \sigma$ and set $(y_1 , \ldots , y_n) := \psi_{\RR}^{-1} (x) \in \RR^n$. Then $\sum_{i=1} y_i \mathsf{p}_i = \rho (\psi_{\RR}^{-1} (x)) \in \widetilde{\sigma}$, and hence (\ref{align:cor:voronoi:variant:2}) holds for any $j=0,1,\ldots , n$. Via the bijection in (\ref{multline:bijection}) and the identity in  (\ref{multline:for:prop:main:unimodular:new}), we have
\[
\tilde{u}' \in 
\underset{u' \in M'}{\argmin} \  
\left( 
Q (x, \lambda_{\RR}^{-1} (b_j) + u') +
\frac{1}{2} Q (u' + \lambda_{\RR}^{-1} (b_j) , u' + \lambda_{\RR}^{-1} (b_j))
\right)
\]
for any $j=0,1,\ldots , n$, which is the requirement in Proposition~\ref{prop:main:unimodular:new}~(ii). Thus Proposition~\ref{prop:main:unimodular:new} is deduced from Proposition~\ref{cor:voronoi:variant:2}.

%\subsection{Matrix representations of the map given by the theta functions}\label{subsec:unimodular}

\subsection{Polyhedral tiling of a Voronoi cell}
\label{subsec:continued}

In this subsection, we construct a certain polyhedral tiling of the Voronoi cell for a lattice in $\RR^n$ with suitable properties (Proposition~\ref{prop:gooddecomposition}). This polyhedral tiling will behave well with respect to the lattice to prove
Proposition~\ref{cor:voronoi:variant:2}. 

We fix the notation in this subsection. Let $\mathsf{p}_1 , \ldots , \mathsf{p}_n \in \RR^n$ be an $\RR$-linearly independent sequence of vectors. We set
\begin{equation}
\label{eqn:sfP}
\mathsf{P}  = \ZZ \mathsf{p}_1 + \cdots + 
\ZZ  \mathsf{p}_n,
\end{equation}
which is a lattice in $\RR^n$. Here, a lattice means a full lattice in the vector space. 
Let $\mathsf{V}$ be the {\em Voronoi cell} of $\mathsf{P}$ around the origin, 
i.e., 
\[
\mathsf{V}
= \{
y \in \RR^n \mid [y, y] \leq [y + \mathsf{p}, y + \mathsf{p}] \quad \text{for any $\mathsf{p} \in \mathsf{P}$} 
\},
\]
where we recall that $[\cdot, \cdot]$ is the standard inner product on $\RR^n$. The translate $\mathsf{p} + \mathsf{V}$ of $\mathsf{V}$ by $\mathsf{p}$ is called the Voronoi cell of $\mathsf{P}$ around $\mathsf{p} \in \mathsf{P}$. We have a polyhedral tiling $\RR^n = \bigcup_{\mathsf{p} \in \mathsf{P}} (\mathsf{p} + \mathsf{V})$ of~$\RR^n$. 

We start with several lemmas concerning the Voronoi cell and lattices. 

\begin{Lemma} \label{lemma:exist-polydecomp}
Let $\mathsf{S} \subseteq \RR^n$ be a finite subset. Then there exists a polyhedral tiling $\Sigma$ of $\mathsf{V}$ such that for any $\mathsf{s} \in \mathsf{S}$ and for any $\sigma \in \Sigma^{(n)}$, the following hold:
\begin{enumerate}
\item[(i)]
there exists a unique $\mathsf{p} \in \mathsf{P}$ such that $\mathsf{s} + \relin (\sigma) \subseteq \mathsf{p} + \relin (\mathsf{V})$;
\item[(ii)]
for any $\mathsf{q} \in \mathsf{P} \setminus \{ \mathsf{p} \}$, $(\mathsf{s} + \relin (\sigma)) \cap (\mathsf{q} + \mathsf{V}) = \emptyset$.
\end{enumerate}
\end{Lemma}

\Proof
Note that for any $x_0 \in \RR^n$, $\bigsqcup_{\mathsf{p} \in \mathsf{P}} (- x_0 + \relin (V) + \mathsf{p})$ is a dense open subspace of $\RR^n$, and for each $\mathsf{p} \in \mathsf{P}$, $- x_0 + \relin (\mathsf{V}) + \mathsf{p}$ is a connected component of this subspace. Note also that this is the union of the relative interiors of the $n$-dimensional polyhedra in $\{ - x_0 + \mathsf{V} + \mathsf{p} \}_{\mathsf{p} \in \mathsf{P}}$.
% maximal faces of the polyhedral decomposition $\RR^n = \bigcup_{\mathsf{p} \in \mathsf{P}}  (- x_0 + \mathsf{V} + \mathsf{p})$, the translate by $-x_0$ of the polyhedral decomposition $\RR^n = \bigcup_{\mathsf{p} \in \mathsf{P}} (\mathsf{p} + \mathsf{V})$.
Then, since $\mathsf{V}$ is a polyhedron and $\mathsf{S}$ is finite, we see that 
\[
U := \bigcap_{\mathsf{s} \in \mathsf{S}}  \left( \bigsqcup_{\mathsf{p} \in \mathsf{P}} (- \mathsf{s} + \relin (\mathsf{V}) + \mathsf{p}) \right) \cap \relin (\mathsf{V})
\]
is a finite disjoint union of the relative interiors of polyhedra of dimension $n$ contained in $\relin (\mathsf{V})$ and is dense in $\mathsf{V}$. 

Let $\Sigma$ be the set of polyhedra that are the closures of the connected components of $U$. We take any $\mathsf{s} \in \mathsf{S}$ and $\sigma \in \Sigma^{(n)}$. By the above construction, $\relin (\sigma)$ is a connected component of $U$. In particular,
$\relin (\sigma) \subseteq \bigsqcup_{\mathsf{p} \in \mathsf{P}} (- \mathsf{s} + \relin (\mathsf{V}) + \mathsf{p})$. Since $\relin (\sigma)$ is connected, there exists a unique $\mathsf{p} \in \mathsf{P}$ such that $\relin (\sigma) \subseteq - \mathsf{s} + \relin (\mathsf{V}) + \mathsf{p}$. This completes the proof of (i).

Note that $\relin (\sigma) \cap (- \mathsf{s} + \relin (\mathsf{V}) + \mathsf{q}) = \emptyset$ for any $\mathsf{q} \in \mathsf{P} \setminus \{  \mathsf{p} \}$. Since $\relin (\sigma)$ is open in $\RR^n$ and $\relin (\mathsf{V})$ is dense in $\mathsf{V}$, this means that $\relin (\sigma) \cap (- \mathsf{s} + \mathsf{V} + \mathsf{q}) = \emptyset$ for any $\mathsf{q} \in \mathsf{P} \setminus \{  \mathsf{p} \}$. This completes the proof of (ii).
\QED

\begin{Lemma} \label{lemma:forcompatibledecomposition}
For any $x \in \mathsf{V}$, there exists a sequence $\mathsf{q}_1 , \ldots , \mathsf{q}_n \in \frac{1}{2} \mathsf{P} \subseteq \RR^n$ such that
$\mathsf{q}_1 , \ldots , \mathsf{q}_n$ is $\RR$-linearly independent in $\RR^n$ and that $\mathsf{q}_i + x \in \mathsf{V}$ for any $i=1, \ldots , n$.
\end{Lemma}

\Proof
Let $\pi : \RR^n \to \RR^n / \mathsf{P}$ be the quotient. 
The subgroup $\frac{1}{2}\mathsf{P} / \mathsf{P} \subseteq \RR^n / \mathsf{P}$ is the subgroup of $2$-torsion points of $\RR^n / \mathsf{P}$ and has order $2^n$. 

We take any $x \in \mathsf{V}$. 
For any $\mathsf{t} \in \frac{1}{2}\mathsf{P}/\mathsf{P}$, since $\pi (\mathsf{V}) = \RR^n / \mathsf{P}$, there exists a $y_{\mathsf{t}} \in \mathsf{V}$ such that $\pi (y_{\mathsf{t}}) = \pi (x) + \mathsf{t}$. 
We set $\mathsf{q}_{\mathsf{t}} := y_{\mathsf{t}} - x$ for each $\mathsf{t} \in \frac{1}{2}\mathsf{P}/\mathsf{P}$. Then for any $\mathsf{t} \in \frac{1}{2}\mathsf{P}/\mathsf{P}$, we have $\mathsf{q}_{\mathsf{t}} + x  = y_{\mathsf{t}} \in \mathsf{V}$, and since $\pi (\mathsf{q}_{\mathsf{t}}) = \mathsf{t}$, we have 
$\mathsf{q}_{\mathsf{t}} \in 
\frac{1}{2}\mathsf{P}$.

To complete the proof, we are going to take an $\RR$-basis of $\RR^n$ out of $\{ \mathsf{q}_{\mathsf{t}} \mid \mathsf{t} \in \frac{1}{2}\mathsf{P} / \mathsf{P} \}$. To do that, it suffices to show that $\{ \mathsf{q}_{\mathsf{t}} \mid \mathsf{t} \in \frac{1}{2}\mathsf{P} / \mathsf{P} \}$ spans $\RR^n$. Let $L$ be the $\RR$-linear subspace spanned by this subset. Since $\mathsf{P} \otimes \QQ = \frac{1}{2}\mathsf{P} \otimes \QQ$ and $\{ \mathsf{q}_{\mathsf{t}} \mid \mathsf{t} \in\frac{1}{2}\mathsf{P} / \mathsf{P} \}$ is a subset of this $\QQ$-vector space, $L \cap \mathsf{P}$ is a lattice in $L$, and thus we have a subtorus
$L / L \cap \mathsf{P} \subseteq \RR^n / \mathsf{P}$. By the definition of $L$, we have $\frac{1}{2}\mathsf{P}/ \mathsf{P} \subseteq L / L \cap \mathsf{P}$, which means that the order of the subgroup of $2$-torsion points in $L / L \cap \mathsf{P}$ is not less than $2^n$. Since $L / L \cap \mathsf{P}$ a real torus of dimension $\dim (L)$, it follows that $\dim (L) \geq n$. Thus $L = \RR^n$. This completes the proof
\QED

For vectors $\mathsf{q}_1 , \ldots , \mathsf{q}_n$ in a vector space, let $\Delta_{\mathsf{q}_1 , \ldots , \mathsf{q}_n}$ denote the convex hull of $\{ 0 , \mathsf{q}_1 , \ldots , \mathsf{q}_n \}$, i.e.,
\[
\Delta_{\mathsf{q}_1 , \ldots , \mathsf{q}_n} = \left\{ t_1 \mathsf{q}_1 +  \cdots + t_n \mathsf{q}_n \;\left|\; 0 \leq t_1 , \ldots , 0 \leq t_n , \sum_{i=1}^n t_i \leq 1
\right.
\right\}.
\]

\begin{Lemma} \label{lemma:gole02}
Let $\mathsf{P}' \subseteq \RR^n$ be a lattice and let $\mathsf{q}_1 , \ldots , \mathsf{q}_n \in \mathsf{P}'$. Assume that  $\mathsf{q}_1 , \ldots , \mathsf{q}_n$ is $\RR$-linearly independent in $\RR^n$. Let $r \in \RR$ be $r \geq   (n-1)!$. Then there exists a $\ZZ$-basis $\mathsf{p}^{(1)} , \ldots , \mathsf{p}^{(n)}$ of the lattice $\frac{1}{r}\mathsf{P}'$ such that $\mathsf{p}^{(1)} , \ldots , \mathsf{p}^{(n)} \in \Delta_{\mathsf{q}_1 , \ldots , \mathsf{q}_n}$. 
\end{Lemma}

\Proof
First, we prove the lemma for $r= (n-1)!$.
We argue by induction on $n$. When $n=1$, it is trivial. Assume that the assertion holds when $n-1$ ($n \geq 2$). 
We set $L = \RR \mathsf{q}_1 + \cdots + \RR \mathsf{q}_{n-1}$, which is an $(n-1)$-dimensional $\RR$-vector space, and $L \cap \mathsf{P}'$ is a lattice in $L \cong \RR^{n-1}$ with $\mathsf{q}_1 , \ldots , \mathsf{q}_{n - 1} \in L \cap \mathsf{P}'$.  
By the induction hypothesis, there exists a $\ZZ$-basis $\mathsf{q}^{(1)} , \ldots , \mathsf{q}^{(n-1)}$ of $\frac{1}{(n-2)!} (L \cap \mathsf{P}') = L \cap \frac{1}{(n-2)!} \mathsf{P}'$ such that 
\[
\mathsf{q}^{(1)} , \ldots , \mathsf{q}^{(n-1)} \in \Delta_{\mathsf{q}_1 , \ldots , \mathsf{q}_{n-1}} \subseteq L. 
\]
Since the $\ZZ$-submodule $L \cap \frac{1}{(n-2)!} \mathsf{P}'$ of $\frac{1}{(n-2)!} \mathsf{P}'$ is saturated in $\frac{1}{(n-2)!} \mathsf{P}'$, there exists a $\mathsf{q^{(n)}} \in \frac{1}{(n-2)!} \mathsf{P}'$ such that $\mathsf{q}^{(1)} , \ldots , \mathsf{q}^{(n-1)} , \mathsf{q}^{(n)}$ is a basis of $\frac{1}{(n-2)!} \mathsf{P}'$. Then there exist $w_1 , \ldots , w_n \in \ZZ$ such that $\mathsf{q}_n = \sum_{j=1}^{n} w_j \mathsf{q}^{(j)}$.
Replacing $\mathsf{q}^{(n)}$ by $-\mathsf{q}^{(n)}$ if necessary, we may and do assume that $w_n >0$. 

We consider the case where $0 \leq w_1 , \ldots , w_{n-1} \leq  w_n$ and $(w_1 , \ldots , w_{n-1}) \neq (0 , \ldots , 0)$. 
Set $\mathsf{p}^{(j)} := \frac{\mathsf{q}^{(j)}}{n-1}$ for $j=1 , \ldots , n-1$ and $\mathsf{p}^{(n)} := \frac{1}{n-1}\sum_{j=1}^n \mathsf{q}^{(j)}$. Then, since $n$ vectors $\mathsf{q}^{(1)} , \ldots , \mathsf{q}^{(n-1)} , \sum_{j=1}^n \mathsf{q}^{(j)}$ form a basis of $ \frac{1}{(n-2)!} \mathsf{P}'$, we see that $\mathsf{p}^{(1)} , \ldots , \mathsf{p}^{(n)}$ is a basis of $\frac{1}{n-1} \left( \frac{1}{(n-2)!} \mathsf{P}'\right) = \frac{1}{(n-1)!}\mathsf{P}'$. Further, for $j = 1, \ldots , n-1$, we have
\[
\mathsf{p}^{(j)} = \frac{\mathsf{q}^{(j)}}{n-1} \in \Delta_{\mathsf{q}_1 , \ldots , \mathsf{q}_{n-1}} \subseteq \Delta_{\mathsf{q}_1 , \ldots , \mathsf{q}_n}.
\]

We show that $\mathsf{p}^{(n)} \in \Delta_{\mathsf{q}_1 , \ldots , \mathsf{q}_n}$. Note that
\begin{align} \label{align:p(n)expression}
\mathsf{p}^{(n)} =
\frac{1}{(n-1)w_n} \mathsf{q}_n + \left( \frac{1}{n-1} - \frac{w_1}{(n-1)w_n} \right) \mathsf{q}^{(1)} + \cdots + \left( \frac{1}{n-1} - \frac{w_{n-1}}{(n-1)w_n} \right) \mathsf{q}^{(n-1)}.
\end{align}
Here,
\begin{align} \label{align:coefficients[01]}
0 < \frac{1}{(n-1)w_n} \leq 1, \quad 
0 \leq \frac{1}{n-1} - \frac{w_i}{(n-1)w_n} \leq 1
\end{align}
for all $i=1 , \ldots , n-1$. 
Further, the sum of the coefficients in the left-hand side on (\ref{align:p(n)expression}) equals 
\[
\frac{1}{(n-1)w_n} + \sum_{i=1}^{n-1} \left( \frac{1}{n-1} - \frac{w_i}{(n-1)w_n} \right)
=
\frac{1- (w_1 + \cdots + w_{n-1})}{(n-1)w_n} + 1.
\]
This is greater than $0$ by (\ref{align:coefficients[01]}). Since $w_1 , \ldots , w_{n-1}$ are nonnegative integers and at least one of them is nonzero, it is not greater than $1$. It follows that
$
\mathsf{p}^{(n)}
\in \Delta_{\mathsf{q}^{(1)} , \ldots , \mathsf{q}^{(n-1)} , \mathsf{q}_n} \subseteq \Delta_{\mathsf{q}_1 , \ldots , \mathsf{q}_n}$.
This proves the assertion for $n$ in this case.

Next, we consider the general case for $r= (n-1)!$. 
For $1 \leq i \leq n-1$, there exist an $h_i \in \ZZ$ and a $w_i' \in \ZZ$ such that
$w_i + h_i w_n = w_i^\prime$ and $0 < w_i^\prime \leq w_n$. We set
\[
\mathsf{q}_n' := \mathsf{q}_n + \sum_{i=1}^{n-1} h_i w_n \mathsf{q}^{(i)} = \sum_{i=1}^{n-1} w_i' \mathsf{q}^{(i)} + w_n \mathsf{q}^{(n)}. 
\] 
Since $0 \leq w_1' , \ldots , w_{n-1}' \leq w_n$ and $(w_1' , \ldots , w_{n-1}') \neq (0, \ldots , 0)$, it follows from the above argument that there exists a basis $\mathsf{p}^{(1)} , \ldots , \mathsf{p}^{(n)}$ of $\frac{1}{(n-1)!} \mathsf{P'}$ such that 
\[
\mathsf{p}^{(1)} , \ldots , \mathsf{p}^{(n)} \in \Delta_{\mathsf{q}^{(1)} , \ldots , \mathsf{q}^{(n-1)} , \mathsf{q}_n'}.
\]

Since $\mathsf{q}^{(1)} , \ldots , \mathsf{q}^{(n)}$ is a basis of the lattice $\frac{1}{(n-2)!} \mathsf{P}'$ (and is an $\RR$-basis of $\RR^n$), the vectors $\mathsf{q}^{(1)} , \ldots , \mathsf{q}^{(n-1)} , \mathsf{q}^{(n)} + \sum_{i=1}^{n-1} h_i \mathsf{q}^{(i)}$ form a basis of $\frac{1}{(n-2)!} \mathsf{P}'$ (and an $\RR$-basis of $\RR^n$). It follows that there exists an automorphism $\psi : \RR^n \to \RR^n$ such that $\psi (\mathsf{q}^{(j)}) = \mathsf{q}^{(j)}$ for any $j=1 , \ldots , n-1$ and $\psi (\mathsf{q}^{(n)}) = \mathsf{q}^{(n)} + \sum_{i=1}^{n-1} h_i \mathsf{q}^{(i)}$. This automorphism restricts to an automorphism on $\frac{1}{(n-2)!} \mathsf{P}'$ and hence induces an automorphism on $\frac{1}{(n-1)!} \mathsf{P}'$. Since $\mathsf{p}^{(1)} , \ldots , \mathsf{p}^{(n)}$ is a basis of $\frac{1}{(n-1)!} \mathsf{P'}$, it follows that $\psi^{-1} (\mathsf{p}^{(1)}) , \ldots , \psi^{-1} (\mathsf{p}^{(n)})$ is a basis of $\frac{1}{(n-1)!} \mathsf{P}'$. Furthermore, we have
\[
\psi (\mathsf{q}_n) = \psi \left( \sum_{j=1}^n w_j \mathsf{q}^{(j)} \right) = \sum_{j=1}^{n-1} w_j \mathsf{q}^{(j)} + w_n \mathsf{q}^{(n)} + w_n \sum_{j=1}^{n-1} h_j \mathsf{q}^{(j)} = 
\mathsf{q}_n + \sum_{j=1}^{n-1} h_j w_n \mathsf{q}^{(j)} = \mathsf{q}_n'
,
\]
which shows that
\[
\psi^{-1} (\mathsf{p}^{(1)}) , \ldots , \psi^{-1} (\mathsf{p}^{(n)}) \in \psi^{-1} \left( \Delta_{\mathsf{q}^{(1)} , \ldots , \mathsf{q}^{(n-1)} , \mathsf{q}_n'} \right) = \Delta_{\mathsf{q}^{(1)} , \ldots , \mathsf{q}^{(n-1)} , \mathsf{q}_n}
\subseteq \Delta_{\mathsf{q}_1 , \ldots , \mathsf{q}_n}
.
\]
This completes the proof of the lemma for $r=(n-1)!$.

Finally, we consider an arbitrary $r \geq (n-1)!$. By what we have shown above, there exists a basis $\mathsf{p}^{(1)} , \ldots, \mathsf{p}^{(n)}$ of $\frac{1}{(n-1)!} \mathsf{P}$ such that each of them sits in $\Delta_{\mathsf{q}_1 , \ldots , \mathsf{q}_n}$. Then
$\frac{(n-1)!}{r} \mathsf{p}^{(1)} , \ldots, \frac{(n-1)!}{r} \mathsf{p}^{(n)}$
is a basis of $\frac{1}{r} \mathsf{P}'$. Further, since $0 \leq \frac{(n-1)!}{r} \leq 1$, they sit in $\Delta_{\mathsf{q}_1 , \ldots , \mathsf{q}_n}$. This completes the proof.
\QED

\begin{Proposition} \label{prop:gooddecomposition}
Let $d_1 , \ldots , d_n$ be positive integers such that $d_i \mid d_{i+1}$ for all $i=1, \ldots , n-1$. Suppose that $d_1 \geq 2 (n-1)!$. Then there exists a polyhedral tiling $\Sigma$ of $\mathsf{V}$ such that for any $\sigma \in \Sigma^{(n)}$, there exists a basis $\mathsf{p}^{(1)} , \ldots , \mathsf{p}^{(n)}$ of 
\[
\ZZ \frac{\mathsf{p}_1}{d_1} + \cdots + 
\ZZ  \frac{\mathsf{p}_n}{d_n}
\]
such that $\mathsf{p}^{(i)} + \sigma \subseteq \mathsf{V}$ for any $i = 1,\ldots , n$.
\end{Proposition}

\Proof
We set $\mathsf{S} := \{ \mathsf{p}' \in \frac{1}{2} \mathsf{P} \mid \| \mathsf{p}' \| \leq \mathrm{diam} (\mathsf{V}) \}$, where $\mathrm{diam} (\mathsf{V})$ is the diameter of $\mathsf{V}$. Then $\mathsf{S}$ is a finite set. 
Let $\Sigma$ be a polyhedral tiling of $\mathsf{V}$ as in Lemma~\ref{lemma:exist-polydecomp} for this $\mathsf{S}$.
We take a $\sigma \in \Sigma^{(n)}$. We fix an $x \in \relin (\sigma)$.
By Lemma~\ref{lemma:forcompatibledecomposition}, there exists a sequence $\mathsf{q}_1 , \ldots , \mathsf{q}_n \in \frac{1}{2} \mathsf{P}$ that is $\RR$-linearly independent and such that $\mathsf{q}_i + x \in \mathsf{V}$ for all $i=1,\ldots , n$. 

We fix any $i=1 , \ldots , n$. Since $x \in \mathsf{V}$ and $\mathsf{q}_i + x \in \mathsf{V}$, we note that $\| \mathsf{q}_i \| \leq \mathrm{diam} (\mathsf{V})$, and thus $\mathsf{q}_i \in \mathsf{S}$. By condition (i) in Lemma~\ref{lemma:exist-polydecomp}, there exists a unique $\mathsf{a}_i \in \mathsf{P}$ such that $\mathsf{q}_i + \relin (\sigma) \subseteq \mathsf{a}_i + \relin (\mathsf{V})$. Since $x \in \relin (\sigma)$ and $\mathsf{q}_i + x \in \mathsf{V} = \mathsf{0} + \mathsf{V}$, we have $(\mathsf{a}_i + \relin (\mathsf{V})) \cap (\mathsf{0} + \mathsf{V}) \neq \emptyset$. By condition (ii) in Lemma~\ref{lemma:exist-polydecomp}, it follows that $\mathsf{a}_i = \mathsf{0}$, and hence $\mathsf{q}_i + \relin (\sigma ) \subseteq \relin (V)$. Thus $\mathsf{q}_i + \sigma \subseteq \mathsf{V}$ for all $i=1 , \ldots , n$. Since $\mathsf{V}$ is convex and $\mathsf{0} \in \mathsf{V}$, it follows that $\Delta_{\mathsf{q}_1 , \ldots , \mathsf{q}_n} + \sigma \subseteq \mathsf{V}$.

For $i=1 , \ldots , n$, we set $\delta_i := \frac{d_i}{d_1}$ and
$\mathsf{P}' := \ZZ \frac{\mathsf{p}_1}{2\delta_1} + \cdots + \ZZ \frac{\mathsf{p}_n}{2 \delta_n}$. We put
$
\mathsf{P}'':=
\ZZ \frac{\mathsf{p}_1}{d_1} + \cdots + 
\ZZ  \frac{\mathsf{p}_n}{d_n}$. Then
$\mathsf{P}''=
 \frac{2}{d_1} \mathsf{P}'$.
By the assumption on $d_1 , \ldots , d_n$, we have $\delta_1 , \ldots , \delta_n \in \ZZ_{\geq 1}$. It follows that $\frac{1}{2} \mathsf{P}  \subseteq \mathsf{P}'$, and thus  $\mathsf{q}_1 , \ldots , \mathsf{q}_n  \in \mathsf{P}'$. Since this sequence is $\RR$-linearly independent in $\RR^n$ and since $\frac{d_1}{2} \geq (n-1)!$,
Lemma~\ref{lemma:gole02} gives us a basis $\mathsf{p}^{(1)} , \ldots , \mathsf{p}^{(n)}$ of 
$\frac{1}{d_1/2}\mathsf{P}' =\mathsf{P}''$
such that $\mathsf{p}^{(1)} , \ldots , \mathsf{p}^{(n)} \in \Delta_{\mathsf{q}_1 , \ldots , \mathsf{q}_n}$.
Then 
$
\mathsf{p}^{(i)} + \sigma \subseteq \Delta_{\mathsf{q}_1 , \ldots , \mathsf{q}_n} + \sigma \subseteq \mathsf{V}
$
for any $i=1 , \ldots , n$,
which completes the proof of the proposition.
\QED

\subsection{Proof of Proposition~\ref{cor:voronoi:variant:2}}

We prove Proposition~\ref{cor:voronoi:variant:2}. This proposition amounts to the following proposition. Let $\mathsf{p}_1 , \ldots , \mathsf{p}_n$ be the column vectors of the matrix $P$ defined by (\ref{eqn:def:pi}), and we regard them as elements in $\RR^n$.

\begin{Proposition}
\label{prop:voronoi:variant:new}
Let $\mathsf{P} := \ZZ \mathsf{p}_1 + \cdots + \ZZ \mathsf{p}_n$ and let $\mathsf{V}$ be the Voronoi cell of $\mathsf{P}$ around the origin. Assume that $d_1 \geq 2 (n-1)!$. Then there exists a polyhedral tiling $\Sigma$ of $\mathsf{V}$ such that for any $\mathsf{p} \in \mathsf{P}$ and $\sigma \in \Sigma$, there exist $\ell^{(0)}, \ell^{(1)}, \ldots, \ell^{(n)} \in \ZZ^n$ 
that satisfy the following conditions:
\begin{enumerate}
\item[(i)]
the vectors $\ell^{(1)} - \ell^{(0)}, \ldots, \ell^{(n)} - \ell^{(0)}$ 
form a basis of the free $\ZZ$-module $\ZZ^n$; 
\item[(ii)]
there exists an  $\widetilde{a} \in  \ZZ^n$ such that 
for any $(y_1 , \ldots , y_n) \in \RR^n$ with $\sum_{i=1}^n y_i \mathsf{p}_i \in \mathsf{p} + \sigma$ and for any $j = 0, 1, \ldots, n$, we have 
%for any $z \in \mathsf{p} + \sigma$ and for any $j = 0, 1, \ldots, n$, writing $z = \sum_{i=1}^n y_i \mathsf{p}_i$, we have
\[
\widetilde{a} \in 
\underset{a \in \ZZ^n} {\argmin} \ 
\left[
\sum_{i=1}^n \left(a_i + \frac{\ell_i^{(j)}}{d_i} + y_i\right) \mathsf{p}_i, \, 
\sum_{i=1}^n \left(a_i + \frac{\ell_i^{(j)}}{d_i} + y_i\right) \mathsf{p}_i
\right]
,
\]
where $\ell^{(j)} = (\ell_1^{(j)} , \ldots , \ell_n^{(j)})$ and $a = (a_1 , \ldots , a_n)$.
%for any $\widetilde{y} \in \mathsf{p} + \sigma$ and 
%for any $j = 0, 1, \ldots, n$. 
\end{enumerate}
\end{Proposition}

Indeed, assume that Proposition~\ref{prop:voronoi:variant:new} holds.
%Let $\mathsf{P} := \ZZ \mathsf{p}_1 + \cdots + \ZZ \mathsf{p}_n$ and let $\mathsf{V}$ be the Voronoi cell of $\mathsf{P}$ around the origin. 
Let $\mathsf{P}$ and $\mathsf{V}$ be as in this proposition.
Let $\Sigma$ be a polyhedral tiling of $\mathsf{V}$ as in Proposition~\ref{prop:voronoi:variant:new} and set $\widetilde{\Sigma} := \{ \mathsf{p} + \sigma \mid  \mathsf{p} \in \mathsf{P} , \sigma \in \Sigma \}$. Then 
%noting (\ref{eqn:argmin}), 
this choice of $\widetilde{\Sigma}$ suffices for Proposition~\ref{cor:voronoi:variant:2}. 

\medskip

\textsl{Proof of Proposition~\ref{prop:voronoi:variant:new}.}
%We prove Proposition~\ref{prop:voronoi:variant:new}. 
We take a polyhedral tiling $\Sigma$ for $\mathsf{V}$ as in Proposition~\ref{prop:gooddecomposition}. 
We take any $\mathsf{p} \in \mathsf{P}$ and $\sigma \in \Sigma$. 
Then, by the choice of $\Sigma$, there exists a basis $\mathsf{p}^{(1)} , \ldots , \mathsf{p}^{(n)}$ of 
$\mathsf{P}'' := 
\ZZ \frac{\mathsf{p}_1}{d_1} + \cdots + 
\ZZ  \frac{\mathsf{p}_n}{d_n}$
such that $\mathsf{p}^{(j)} + \sigma \subseteq \mathsf{V}$ for all $j=1,\ldots , n$. We set $\mathsf{p}^{(0)} := \mathsf{0}$.
%Set $\ell^{(0)} := \mathsf{0}$. 
For $j = 0, 1, \ldots , n$, we take $\ell^{(j)} = (\ell^{(j)}_1, \ldots , \ell^{(j)}_n) \in \ZZ^n$ in such a way that
 $\mathsf{p}^{(j)} = \sum_{i=1}^n  \ell_i^{(j)} \frac{\mathsf{p}_i}{d_i}$. We remark that $\ell^{(0)} = \mathsf{0}$ and that $\mathsf{p}^{(j)} + \sigma \subseteq \mathsf{V}$ holds also for $j=0$. 

We consider condition (i) in the proposition.
Since $\mathsf{p}^{(1)} , \ldots , \mathsf{p}^{(n)}$ is a basis of $\mathsf{P}''$, 
$\ell^{(1)} , \ldots , \ell^{(n)} $ is a basis of 
$\ZZ^n$. Since $\ell^{(0)} = \mathsf{0}$, 
This shows that 
those $\ell^{(0)}, \ell^{(1)}, \ldots, \ell^{(n)} \in \ZZ^n$ satisfy condition (i).

Finally, we consider condition (ii) in the proposition.
Since $\mathsf{p}_1 , \ldots , \mathsf{p}_n$ is a basis of $\mathsf{P}$, there exists a 
$c := (c_1, \ldots , c_n) \in \ZZ^n$ such that $\mathsf{p}
= \sum_{i=1}^n c_i  \mathsf{p}_i$. 
We 
%set $\widetilde{\sigma} = \mathsf{p} + \sigma$ and 
take any $(y_1 , \ldots , y_n) \in \RR^n$ with
$\sum_{i=1}^n y_i \mathsf{p}_i \in \mathsf{p} + \sigma$. We take any $j = 0,1, \ldots , n$. It suffices to show that
\[
- c \in
\underset{a \in \ZZ^n} {\argmin} \ 
\left[
\sum_{i=1}^n \left(a_i + \frac{\ell_i^{(j)}}{d_i} + y_i\right) \mathsf{p}_i, \, 
\sum_{i=1}^n \left(a_i + \frac{\ell_i^{(j)}}{d_i} + y_i\right) \mathsf{p}_i
\right]
,
\]
where $a = (a_1, \ldots ,a_n) \in \ZZ^n$. 
%Noting that $\mathsf{p}^{(j)} + \sigma \subseteq \mathsf{V}$, w
We compute that
\[
\sum_{i=1}^n
 \left( -c_i + \frac{\ell_i^{(j)}}{d_i} + y_i \right) \mathsf{p}_i = 
\sum_{i=1}^n y_i \mathsf{p}_i - \mathsf{p} + \mathsf{p}^{(j)} \in \sigma + \mathsf{p}^{(j)} \subseteq \mathsf{V}. 
\]
It follows that
\[
\mathsf{0} \in 
\underset{a \in \ZZ^n} {\argmin} \ 
\left[
\sum_{i=1}^n \left(a_i -c_i  + \frac{\ell_i^{(j)}}{d_i} + y_i\right) \mathsf{p}_i, \, 
\sum_{i=1}^n \left(a_i -c_i  + \frac{\ell_i^{(j)}}{d_i} + y_i \right) \mathsf{p}_i
\right]
,
\]
and thus 
\[
-c \in 
\underset{a \in \ZZ^n} {\argmin} \ 
\left[
\sum_{i=1}^n \left(a_i   + \frac{\ell_i^{(j)}}{d_i} + y_i\right) \mathsf{p}_i, \, 
\sum_{i=1}^n \left(a_i   + \frac{\ell_i^{(j)}}{d_i} + y_i \right) \mathsf{p}_i
\right]
.
\]
Thus the proof is complete. 
\QED

\setcounter{equation}{0}
\section{Uniformization theory of Raynaud and Bosch--L\"utkebohmert}
\label{section:uniformization}

We move on to the faithful tropicalization problem for the canonical skeleton of an abelian variety. The main ingredients will be the nonarchimedean theta functions and their Fourier expansions. We will discuss how they behave on the canonical skeleton for the proof of the main results. 

In this section, for the reader's convenience and to fix various sign conventions, we recall the nonarchimedean uniformization theory by Raynaud and Bosch--L\"utkebohmert, as the theory of nonarchimedean theta functions are developed over the uniformization. The basic reference is \cite{BoschLutke-DAV}. 

\subsection{Notation and conventions}
We fix the notation and conventions in the sequel. Let $k$ be a complete and algebraically closed nonarchimedean valued field. Assume that the absolute value $|\ndot|_k$ of $k$ is nontrivial. Let $k^{\circ}$ denote the valuation ring of $k$ and let $k^{\circ \circ}$ denote the maximal ideal. We set $|k^{\times}| := \{ |\alpha|_k \mid \alpha \in k^{\times} \}$.

For an algebraic scheme $X$ over $k$, let $X^{\an}$ denote the Berkovich analytic space associated to $X$. Recall that the topological space $X^{\an}$ can be described as follows.
The underlying set is given by
\[
X^{\an} =
\{ (p ,|\ndot | ) \mid \text{$p \in X$, $|\ndot|$ is an absolute value on $\kappa (p)$ that extends $|\ndot|_k$} \}
,
\]
where $\kappa (p)$ denotes the residue field of $X$ at $p$. For any $x = (p,|\ndot|) \in X^{\an}$ and  for any regular function $f$ (in scheme theoretic sense) on a neighborhood of $p$, we write $|f(x)| := |f (p)|$, where $f(p)$ is the class in $\kappa (p)$. The topology of $X^{\an}$ is the weakest topology such that for any Zariski open subset $U \subseteq X$ and for any $f \in \OO_{X}(U)$, $U^{\an}$ is an open subset of $X^{\an}$ and the map $|f| : U^{\an} \to \RR$ given by $x \mapsto |f(x)|$ is continuous.

For any $x = (p , |\ndot|) \in X^{\an}$, let $\mathscr{H} (x)$ denote the completed residue field of $X^{\an}$ at $x$, which is the completion of $\kappa (p)$ with respect to $|\ndot|$. 

There exists a natural injection $X(k) \hookrightarrow X^{\an}$. Via this, we regard $X(k) \subseteq X^{\an}$, and a point in $X(k)$ is called a \emph{classical points}. If we write $X^{\an} (k)$ for the set of points in $X^{\an}$ with completed residue field $k$, then $X(k) = X^{\an}(k)$.

Let $f : X \to Y$ be a morphism of algebraic schemes over $k$. Then we have a natural morphism $f^{\an} : X^{\an} \to Y^{\an}$. Indeed, $f^{\an} (p,|\ndot|) = (f(p) , |f^{\ast} (\ndot)|)$. We call $f^{\an}$ the \emph{analytification} of $f$.

Let $L$ be a line bundle on $X$. Then $L^{\an}$ is a line bundle on $X^{\an}$. For any $x \in X^{\an}$, let 
$\rest{L^{\an}}{x}$ denote the restriction of $L^{\an}$ over the single point $\{ x \}$, or in other words, the fiber of the line bundle $L^{\an}$ over $x$. If $x \in X(k)$ and if $\rest{L}{x}$ denote the fiber of $L$ over $x$, then $\rest{L}{x}^{\an} = \rest{L^{\an}}{x}$ naturally. In particular, when we consider the classical points of the fiber of $L^{\an}$ over a classical point $x$, we do not have to distinguish $L$ and $L^{\an}$, as $\rest{L}{x} (k) = \rest{L^{\an}}{x}(k)$. 

Let $s$ be a section of a line bundle $L$ over a Zariski open subset $U \subseteq X$. Then a section $s^{\an}$ of $L^{\an}$ over $U^{\an}$, as an analytification of $s$, canonically arises. When there is no danger of confusion, let the same symbol $s$ denote the analytification of $s$.

We identify a line bundle with the invertible sheaf in a \emph{covariant} way, while in \cite{BoschLutke-DAV, FRSS}, they do in a contravariant way. For a locally ringed space, the structure sheaf is regarded as a trivial line bundle in an obvious way. 

A \emph{formal model} of $X^{\an}$ means an admissible formal scheme $\mathscr{X}$ over $k^{\circ}$ equipped with an isomorphism $\mathscr{X}^{\an} \cong X^{\an}$, where $\mathscr{X}^{\an}$ denotes the Raynaud generic fiber of $\mathscr{X}$. A formal model $\mathscr{X}$ of $X^{\an}$ is said to be \emph{proper} if $\mathscr{X}$ is proper over $k^{\circ}$. In addition, let $L$ be a line bundle on $X$. A \emph{formal model} of $(X^{\an},L^{\an})$ is a pair $(\mathscr{X} , \mathscr{L})$ of a formal model $\mathscr{X}$ of $X^{\an}$ and a line bundle $\mathscr{L}$ on $\mathscr{X}$ equipped with an isomorphism $\rest{\mathscr{L}}{X^{\an}} \cong L^{\an}$ via the specified isomorphism $\mathscr{X}^{\an} \cong X^{\an}$. A formal model $(\mathscr{X} , \mathscr{L})$ of $(X^{\an} , L^{\an})$ is said to be \emph{proper} if the formal model $\mathscr{X}$ of $X$ is proper.

Let $\mathcal{B} \to S$ is an abelian scheme (or a formal abelian scheme) with zero-section $0$. Let $\mathcal{L}$ be a line bundle on $\mathcal{B}$. A \emph{rigidification} for $\mathcal{L}$ means a isomorphism $0^{\ast} (\mathcal{L}) \cong \OO_S$. A line bundle with a rigidification is called a \emph{rigidified line bundle}. An isomorphism between rigidified line bundles is an isomorphism between line bundle that respects the rigidifications. Note that an isomorphism between rigidified line bundles is unique.

\subsection{Raynaud extensions}

We recall the notion of Raynaud extensions for abelian varieties. The basic reference is \cite{BoschLutke-DAV}; we also refer to \cite{BoschLutke-SR2}.

\subsubsection{Formal analytic Raynaud extensions}

An analytic group $G$ is said to be \emph{formal} if it is the Raynaud generic fiber of an admissible formal group scheme $\mathscr{G}$ over $R$ with reduced special fiber. We say that a formal analytic group $G$ has \emph{semiabelian reduction} if the above $\mathscr{G}$ can be take to have semiabelian reduction.

\begin{Remark}
There is a category of formal analytic varieties, and it is known that this category is equivalent to the category of admissible formal schemes with reduced fiber; see \cite[Section~1]{BoschLutke-NM} and \cite[Section~1]{Gubler-Heights}. Thus by \cite[Korollar~1.5]{Bosch}, for $G$ above, such an admissible formal group scheme $\mathscr{G}$ as above  is unique. 
\end{Remark}

Let $A$ be an abelian variety over $k$. Then by \cite[Theorem~8.2]{BoschLutke-SR2}, there exists a unique analytic subdomain $A^{\an}_1$ of $A^{\an}$ that is a analytic subgroup of $A^{\an}$ and is a formal analytic group with semi-abelian reduction;  further, there exists an exact sequence
\begin{align} \label{align:formalanalyticRE}
\begin{CD}
1 @>>> T^{\an}_1 @>>> A^{\an}_1 @>>> \mathscr{B}^{\an} @>>> 0
\end{CD}
\end{align}
of analytic groups, where $T^{\an}_1$ is an affinoid torus and $\mathscr{B}$ is an admissible formal abelian scheme over $k^{\circ}$. This exact sequence is unique up to canonical isomorphism, and it is called the \emph{formal analytic Raynaud extension} for $A$. Furthermore, by \cite[Theorem~1.1]{BoschLutke-DAV}, there exists an exact sequence
\begin{align} \label{align:formalRE3}
\begin{CD}
1 @>>> \mathscr{T}_1 @>>> \mathscr{A}_1 @>>> \mathscr{B} @>>> 0
\end{CD}
\end{align}
of admissible formal groups schemes such that the sequence of analytic groups obtained by taking the Raynaud generic fibers of (\ref{align:formalRE3}) agrees with (\ref{align:formalanalyticRE}). We call (\ref{align:formalRE3}) the \emph{formal Raynaud extension} for $A$.
It is known that the formal Raynaud extension is a $\mathscr{T}_1$-torsor; see the second paragraph after Theorem~1.1 in \cite{BoschLutke-DAV}.

Let 
\[
\begin{CD}
1 @>>> T^{\an}_1 @>>> A^{\an}_1 @>{q_1}>> \mathscr{B}^{\an} @>>> 0 \\
@. @VVV @VVV @V{\mathrm{id}}VV \\
1 @>>> T^{\an} @>>> \mathbf{E} @>{q}>> \mathscr{B}^{\an} @>>> 0.
\end{CD}
\]
be the pushout of by the canonical inclusion $T^{\an}_1 \hookrightarrow T^{\an}$. The second exact sequence in the above commutative diagram of analytic spaces is called the \emph{Raynaud extension} of $A$.

Since $q_1 : A^{\an}_1 \to \mathscr{B}^{\an}$ is a $T^{\an}_1$-torsor, $q : \mathbf{E} \to \mathscr{B}^{\an}$ is a $T^{\an}$-torsor. Let 
\[
\{ (V_{\lambda})_{\lambda} , (\rho_{\lambda} : T^{\an} \times V \to q^{-1} (V_{\lambda}))_{\lambda} , (g_{\mu \lambda} : V_{\lambda} \cap V_{\mu} \to T^{\an} )_{\lambda, \mu} \}
\]
be a trivialization of the $T^{\an}$-torsor $q^{\an}$. We call this a \emph{formal $T^{\an}_1$-trivialization} of the Raynaud extension for $A$ if $V_{\lambda}$ is the Raynaud generic fiber of an open subset of the formal model $\mathscr{B}$ and the transition morphisms of the trivialization value in $T^{\an}_1$.
Since the Raynaud extension is the pushout of the formal analytic Raynaud extension and the formal analytic Raynaud extension is the Raynaud generic fiber of the formal Raynaud extension, it has a formal $T^{\an}_1$-trivialization.

\subsubsection{Valuation map}

We keep the above notation. We recall the valuation maps arising from the Raynaud extension. Since $\mathscr{T}_1$ is a formal torus over $k^{\circ}$, we write $\mathscr{T}_1 = \mathrm{Spf} (k^{\circ} [\chi^M])$ for some free $\ZZ$-module $M$ of finite rank, where $\chi^{M}$ is the character group of $M$. Then $T = \Spec (k [\chi^M])$, and we call $M$ the \emph{character lattice} of $T$. We set $N: = \mathrm{Hom} (M , \ZZ)$ and $N_{\RR} :=N \otimes \RR = \mathrm{Hom} (M , \RR)$.

First, we recall the valuation map $\val : T^{\an} \to N_{\RR}$. We take any $x \in T^{\an}$. Then for any $u \in M$, we define a map $u \mapsto - \log |\chi^u (x)|$, which is a group homomorphism from $M$ to $\RR$. This assignment defines a map $\val : T^{\an} \to N_{\RR} = \mathrm{Hom} (M,\RR)$. 

Using a formal $T^{\an}_1$-trivialization, we construct a valuation map $\val_{\mathbf{E}} : \mathbf{E} \to N_{\RR}$. We take a formal $T^{\an}_1$-trivialization 
\[
\{ (V_{\lambda})_{\lambda} , (\rho_{\lambda} : T^{\an} \times V_{\lambda} \to q^{-1} (V_{\lambda}))_{\lambda} , (g_{\mu \lambda} : V_{\lambda} \cap V_{\mu} \to T^{\an} )_{\lambda, \mu}\}
\] 
of the Raynaud extension. For each $\lambda$, we define a map $\val_{\lambda} : q^{-1} (V_{\lambda}) \to N_{\RR}$ to be the composite 
\[
\begin{CD}
q^{-1} (V_{\lambda}) @>{\rho_{\lambda}^{-1}}>> T^{\an} \times V_{\lambda} @>{\mathrm{pr}_1}>> T^{\an} @>{\val}>> N_{\RR} ,
\end{CD}
\]
where $\mathrm{pr}_1$ is the projection to the first factor.
Then, since $g_{\mu \lambda}$ values in $T^{\an}_1$, those maps for all $\lambda$ patch together to be a continuous map 
\[
\val_{\mathbf{E}} : \mathbf{E} \to N_{\RR}, 
\]
which is called the \emph{valuation map} for $\mathbf{E}$. This is well-defined for the Raynaud extension and does not depend on the choice of the formal $T_1^{\an}$-trivialization.

By \cite[Theorem~8.8]{BoschLutke-SR2} and \cite[Theorem~1.2]{BoschLutke-DAV}, the canonical homomorphism $A^{\an}_1 \hookrightarrow A^{\an}$ extends to unique homomorphism $p : \mathbf{E} \to A^{\an}$. This is surjective, and $\Ker (p)$ is a lattice in $E^{\an}$ in the sense that it is a discrete subgroup of $\mathbf{E} (k)$, $\val_{\mathbf{E}} (\Ker (p))$ is a (full) lattice in $N_{\RR}$, and the homomorphism $\rest{\val}{\Ker (p)} : \Ker (p) \to \val_{\mathbf{E}} (\Ker (p))$ is an isomorphism.
We set $M' := \val (\Ker (p)) \subseteq N_{\RR}$ and let $\widetilde{\Phi} : M' \to \mathbf{E}$ be the injective map given by the inverse of the isomorphism $\val_{\mathbf{E}} : \Ker (p) \to M'$. We often identify $M'$ with a $\Ker (p)$ via this isomorphism.
The valuation map $\val_{\mathbf{E}} : \mathbf{E} \to N_{\RR}$ descends to map $\val_{A^{\an}} : A^{\an} \to N_{\RR}/M'$, and thus we have a commutative diagram
\[
\begin{CD}
\mathbf{E} @>{\val_{\mathbf{E}}}>> N_{\RR} \\
@VVV @VVV \\
A^{\an} @>{\val_{A^{\an}}}>> N_{\RR} / M' .
\end{CD}
\]
Note that $N_{\RR} / M'$ is a real torus with an integral structure. We call $\val_{A^{\an}}$ the \emph{valuation map} for $A^{\an}$.
We call $N_{\RR}/M'$ the \emph{canonical tropicalization} of $A$.

\begin{Remark} \label{remark:algebraizable-extension}
Since $A^{\an} = \mathbf{E} / M'$ is algebraizable, it follows from \cite[Theorem~6.13]{BoschLutke-DAV} that $\mathscr{B}^{\an}$ is algebraizable. Thus there exists an abelian variety (unique up to isomorphism) such that $\mathscr{B}^{\an} = B^{\an}$. As is mentioned in the last paragraph in \cite[Page~655]{BoschLutke-DAV}, it turns out that the Raynaud extension is algebraizable, i.e., there exists an exact sequence
\begin{align} \label{align:alg:RE}
\begin{CD}
1 @>>> T @>>> E @>>> B @>>> 0 
\end{CD}
\end{align}
of which analytification is the Raynaud extension of $A$. 
\end{Remark}

We call the exact sequence of algebraic group in (\ref{align:alg:RE}) the \emph{algebraized Raynaud extension} of $A$.

In summary, we have the following diagram, which we call the \emph{Raynaud cross} of $A$:
\begin{align} \label{align:RaynaudCross}
\begin{CD}
@. M' \\
@. @VV{\widetilde{\Phi}}V \\
 T^{\an} @>>> E^{\an} @>{q^{\an}}>> B^{\an}  \\
@. @V{p}VV \\
@. A^{\an}
\end{CD}
\end{align}
where $T = \Spec ( k[\chi^M])$ is an algebraic torus with character lattice $M$ of rank $n = \dim (T)$, $B$ is an abelian variety with good reduction, $M'$ is a lattice in $N_{\RR}$ with $\rank (M') = \rank (M)$, and $\widetilde{\Phi}$ is the natural injective homomorphism with $\widetilde{\Phi} (M') \subseteq E (k)$. We call $n = \dim (T)$ the \emph{torus rank} of $A$. The horizontal sequence is algebraizable, but the vertical one is not unless $A$ has torus rank $0$.

\subsubsection{Poincar\'e bundle valued bilinear form}
We define a bilinear form $t$ on $M' \times M$ with values in the Poincar\'e bundle of the abelian part $B$ of $A$. We refer to \cite[Section~3]{BoschLutke-DAV} for the detail.

Let $B'$ be the dual abelian variety of $B$ and let $P$ be the rigidified Poincar\'e line bundle on $B \times B'$. Here, the rigidification means the rigidification as a line bundle on the abelian variety $B \times B'$ over $k$. We always assume that the Poincar\'e line bundle is normalized so that $\rest{P}{\{  0 \} \times B'} \cong \OO_{\{ 0 \} \times B'}$. Then the rigidification on $P$ induces a rigidification $\phi: \rest{P}{\{ 0 \} \times B'} \cong \OO_{\{ 0 \} \times B'}$ on $P$ as a line bundle on the abelian scheme $B \times B'$ over $B'$ such that the restriction of this rigidification $\phi$ over the point $(0,0) \in B \times B'$ coincides with the original rigidification of $P$ as a line bundle on the abelian variety $B \times B'$ over $k$.

\begin{Remark} \label{remark:hom-algebraicextension}
As is explained in \cite[Page~665]{BoschLutke-DAV}, there is a bijective correspondence between $\mathrm{Hom} (M , B'(k))$ and the isomorphism classes of extensions of $B^{\an}$ by $T^{\an}$ in the category of analytic group. Since an extension of $B^{\an}$ by $T^{\an}$ is algebraizable, as is noted in Remark~\ref{remark:algebraizable-extension}, it follows that there is a bijective correspondence between $\mathrm{Hom} (M , B'(k))$ and the isomorphism classes of extensions of $B$ by $T$ in the category of algebraic groups.
\end{Remark}

We take $\Phi' \in \mathrm{Hom} (M , B'(k))$ corresponding, in the sense of Remark~\ref{remark:hom-algebraicextension}, to the algebraized Raynaud extension
\[
\begin{CD}
1 @>>> T @>>> E @>>> B @>>> 0
\end{CD}
\]
of $A$. For each $u \in M$, take the pushout 
\[
\begin{CD}
1 @>>> T @>>> E @>>> B @>>> 0 \\
@. @V{\chi^u}VV @V{e_u}VV @VV{\id_B}V \\
1 @>>> \mathbb{G}_{m} @>>> E_u @>>> B @>>> 0
\end{CD}
\]
by $\chi^u$. Then $E_u$ is an algebraic group and the homomorphism $E_u \to B$ is an $\GG_m$-torsor. Note that a $\GG_m$-torsor over $B$ is naturally identified with a line bundle on $B$. Then, as is explained at the beginning of \cite[Page~666]{BoschLutke-DAV}, we can regard $E_u \to B$ as a homologically trivial rigidified line bundle, and hence we have an identification 
\begin{align} \label{align:E-poincare}
E_u = \rest{P}{B \times \{ \Phi' (u) \}}
\end{align} of rigidified line bundles; in other words, $\Phi' (u) \in B'(k)$ corresponds to the isomorphism class of (the line bundle corresponding to) $E_u$. Note that, via the identification $E_u^{\an} = \rest{P}{B \times \{ \Phi' (u) \}}^{\an}$, we have $e_u^{\an} (x) \in \rest{P}{(q^{\an}(x) , \Phi' (u) )}^{\an}$ for any $x \in B^{\an}$. Thus we have a map $E^{\an} \times M \to P^{\an}$ given by $(x,u) \mapsto e_u^{\an} (x)$. This map is bilinear in the sense that for any $u \in M$, $e_u^{\an} : E^{\an} \to E_u^{\an}$ is a homomorphism of analytic groups and that
$e_{u_1 + u_2}^{\an} (x) =e_{u_1}^{\an} (x) \otimes e_{u_2}^{\an} (x)$ via the canonical isomorphism $E_{u_1 + u_2}^{\an} = E_{u_1}^{\an} \otimes E_{u_2}^{\an}$; see the paragraph just before Proposition~3.1 in \cite{BoschLutke-DAV}. 

Note that if  $x \in B(k)$, then $e_u (x) = e^{\an} (x) \in E_u(k) = \rest{P}{(q(x) ,\Phi (u))} (k)$. Then we define a map $t : M' \times M \to P(k)$ by 
\begin{align} \label{align:def:t}
t (u' , u) := e_{u} (\widetilde{\Phi} (u')) \in \rest{P}{(q(x) ,\Phi (u))} (k) \subseteq P(k)
\end{align} 
for any $(u',u) \in M' \times M$.
It follows from what we mentioned above that $t$ is bilinear. We remark that if $u =0$ or $u'=0$, then $t(u' , 0) = t (0,u) = 1 \in k = \rest{P}{(0,0)} (k)$, where the last equality is the rigidification of $P$.

\subsection{Canonical metrics} \label{subsection:canonicalmetric}

In this subsection, we briefly recall the canonical metrics on rigidified line bundles on abelian varieties with good reduction. These metrics are mentioned as ``the canonical absolute value on the Poincar\'e bundle'' in \cite{BoschLutke-DAV} and as ``the canonical model metrics'' in \cite[Subsection~3.1]{FRSS}.

Let $X$ be an algebraic scheme over $k$. Let $L$ be a line bundle on $X$. A \emph{metric} on $L^{\an}$, or on $L$ by abuse of words, is a map $\| \ndot \| : L^{\an} \to \RR$ such that for any $x = (p, | \ndot |)  \in X^{\an}$ and for any local sections $s$ of $L$ and $f$ of $\OO_X$ over a Zariski open neighborhood of  $p \in X$, 
$\| s^{\an} (x) \| = 0$ implies $s^{\an} (x) = 0$ and
$\| f^{\an} (x) s^{\an} (x)\| = |f^{\an} (x)| \| s^{\an} (x) \|$ hold.

Assume that $X$ is proper over $k$. Let $(\mathscr{X} , \mathscr{L})$ be a proper formal model of $(X,L)$. Then we have a formal model metric $\| \ndot \|$ induced from $(\mathscr{X} , \mathscr{L})$, which is characterized as follows. Take any $x = (p , |\ndot|) \in X^{\an}$. Let $\mathscr{H} (x)^{\circ}$ be the valuation ring of the completed residue field $\mathscr{H}(x)$ at $x$. Since $\mathscr{X}$ is proper over $k^{\circ}$, we have a morphism $\widetilde{x} : \mathrm{Spf} (\mathscr{H} (x)^{\circ}) \to \mathscr{X}$ over $k^{\circ}$ whose restriction to the Raynaud generic fiber equals the point $x \in X^{\an}$. Here, we regard the line bundles on $\mathscr{X}$ as invertible sheaves. Then there exists an isomorphism $\phi : \widetilde{x}^{\ast} (\mathscr{L}) \to \mathscr{H} (x)^{\circ}$. Let $\phi_{\mathscr{H} (x)} : \widetilde{x}^{\ast} (\mathscr{L}) \otimes_{\mathscr{H}(x)^{\circ}} \mathscr{H}(x) \to \mathscr{H}(x)$ be the base change of $\phi$ to $\mathscr{H} (x)$. For a local section $s$ of $L$ at $p$, regarding $s^{\an} (x) \in \widetilde{x}^{\ast} (\mathscr{L}) \otimes_{\mathscr{H} (x)^{\circ}} \mathscr{H} (x)$, we have $\phi_{\mathscr{H} (x)} (s^{\an}(x)) \in \mathscr{H} (x)$. Then $\| s^{\an} (x) \| = |\phi_{\mathscr{H} (x)} (s^{\an}(x))|$. One checks that a formal model metric is continuous in the sense that for any Zariski open subset $U \subseteq X$ and for any section $s$ of $L$ over $U$, the map $U^{\an} \to \RR$ given by $x \mapsto \| s^{\an} (x) \|$ is continuous.

\begin{Remark} \label{remark:value:norm}
With the above notation, if $x \in X(k)$, then $\| s^{\an} (x) \| = |\phi_{\mathscr{H} (x)} (s^{\an}(x))| \in |k^{\times}|$. Further, if $s$ is a global section of $L$, then it follows from the maximal modulus principle that $\| s \|_{\sup} \in |k^{\times}|$.
\end{Remark}

Let $B$ be an abelian variety over $k$ and assume that it has good reduction. Let $L$ be a rigidified line bundle on $B$. Then there exists a formal model $(\mathscr{B} , \mathscr{L})$ of $(B,L)$ such that $\mathscr{L}$ has a rigidification whose restriction to $B^{\an}$ equals the given rigidification on $L$. Then we have a formal model metric $\| \ndot \|$ on $L$. Since such a formal model is unique up to canonical isomorphism of rigidified line bundles, this formal model metric is uniquely determined. We call this metric on $L^{\an}$ the \emph{canonical metric}.

In general, there is a notion of canonical metric on any rigidified line bundle on any abelian variety over $k$, which is characterized by a certain compatibility condition with $m$-times endomorphism for an integer $m \geq 2$, and the canonical metric in this sense is unique for a rigidified line bundle. For a rigidified line bundle on an abelian variety with good reduction, the canonical metric defined above by using the formal model equals the canonical metric in this sense.

\begin{Remark} \label{remark:canonicalmetric}
One checks that the following hold for the canonical metrics:
\begin{enumerate}
\item
the pullback of the canonical metric by a homomorphism is the canonical metric;
\item
the tensor product of the canonical metrics is the canonical metric;
\item
if $L_1$ and $L_2$ are rigidified line bundles with canonical metric and if $\phi : L_1 \to L_2$ is an isomorphism of rigidified line bundles, then $\phi$ is an isometry.
\end{enumerate}
\end{Remark}

\subsection{Nonarchimedean Appell--Humbert theorem}

In this subsection, we recall the nonarchimdean Appell--Humbert theory for line bundles on abelian varieties. We refer to \cite[Sections~3--6]{BoschLutke-DAV} for the detail.

For a rigidified line bundle $L$ on $B$, we set
\begin{align*} 
\mathscr{D}_2^- (L) := m^{\ast} (L^{\otimes -1}) \otimes p_1^{\ast} (L) \otimes p_2^{\ast} (L)
.
\end{align*}
Then by the universal property of the rigidified Poincar\'e bundle, there exists a unique homomorphism $\varphi_L^- : B \to B'$ such that
\begin{align} \label{align:translate-canonicalisom00}
(\id_B \times \varphi_L^-)^{\ast} (P) = \mathscr{D}_2^- (L).
\end{align}

We remark that in \cite{FRSS}, they use the notation $\mathscr{D}_2 L$ and $\varphi_L$. Note that our $\mathscr{D}_2^- (L)$ is the dual of $\mathscr{D}_2 L$, and $\varphi_L^-$ is $[-1] \circ \varphi_L$. We put ``$-$'' to distinguish our notation it form that in \cite{FRSS}.

Let $L$ be a rigidified line bundle on $B$ and put the canonical metric on it. Then we put a metric on $\mathscr{D}_2^- (L)$ by the pullbacks and the tensor products, and it turns out that this metric is the canonical metric by Remark~\ref{remark:canonicalmetric}~(1) and (2).  By pulling back the canonical metric on $P$, $(\id_B \times \varphi_L^-)^{\ast} (P) $ is equipped with the canonical metric. By Remark~\ref{remark:canonicalmetric}~(3), the isomorphism (\ref{align:translate-canonicalisom00}) is an isometry with respect to those metrics. Furthermore, restricting the isomorphism on (\ref{align:translate-canonicalisom00}) to $B \times \{ \varphi_L^- (x) \}$ for any $x \in B(k)$, we obtain a canonical isomorphism
\begin{align} \label{align:translate-canonicalisom0}
\rest{P}{B \times \{ \varphi_L^- (x) \}} = T_x^{\ast} (L^{\otimes -1}) \otimes L \otimes L (x)
,
\end{align}
where $T_x : B \to B$ is the translation by $x$ and $L(x)$ is the trivial line bundle on $B$ with fiber $\rest{L}{x}$ over $0 \in B(k)$. Since (\ref{align:translate-canonicalisom00}) is an isometry with respect to the canonical metrics, the isomorphism on (\ref{align:translate-canonicalisom0}) is also isometry with respect to the metric given by restriction.

We define $\Phi : M' \to B(k)$ by $\Phi := q \circ \widetilde{\Phi}$. Recall that we have $\Phi' \in \mathrm{Hom} (M , B' (k))$ corresponding to the Raynaud extension; see Remark~\ref{remark:hom-algebraicextension}. Let $\lambda : M' \to M$ be a homomorphism and assume that it satisfies the compatibility condition 
\begin{align} \label{align:compatibility}
\Phi' \circ \lambda = \varphi_L^- \circ \Phi.
\end{align} 
Then we have
\begin{align} \label{align:valueoft}
t (u_1' , \lambda (u_2')) \in \rest{P}{(\Phi (u_1') , \Phi' (\lambda (u_2'))} (k) =
\left( \rest{L^{\otimes -1}}{\Phi (u_1') + \Phi (u_2')} \otimes \rest{L}{\Phi (u_1')} \otimes \rest{L}{\Phi (u_2')} \right)(k)
,
\end{align}
where $t : M' \times M \to P(k)$ is the map defined in (\ref{align:def:t}). Indeed,  for any $(u_1' , u_2') \in M' \times M'$, restricting the identification on (\ref{align:translate-canonicalisom00}) to the fibers over $(\Phi (u_1') , \Phi (u_2'))$, we have an identity
\[
\rest{(\id_B \times \varphi_L^-)^{\ast} (P)}{( \Phi (u_1') , \Phi (u_2'))} = \rest{L^{\otimes -1}}{\Phi (u_1') + \Phi (u_2')} \otimes \rest{L}{\Phi (u_1')} \otimes \rest{L}{\Phi (u_2')}
,
\]
and since $\varphi_L^- \circ \Phi = \Phi' \circ \lambda$, we obtain the equality in  (\ref{align:valueoft}).
Note that $\tau^{\ast} (\mathscr{D}_2^- (L)) = \mathscr{D}_2^- (L)$,
where $\tau : B \times B \to B \times B$ is the switching of the first and the second factors.
Then we have the canonical identification
\[
\rest{\tau^{\ast} (\id_B \times \varphi_L^-)^{\ast} (P)}{( \Phi (u_1') , \Phi (u_2'))} = \rest{L^{\otimes -1}}{\Phi (u_1') + \Phi (u_2')} \otimes \rest{L}{\Phi (u_1')} \otimes \rest{L}{\Phi (u_2')}
,
\]
and $t (u_2' , \lambda (u_1'))$ is also regarded as a $k$-point of the fiber of this line bundle.

\begin{Remark} \label{remark:tissymmetric}
We keep the assumption that $\lambda$ satisfies (\ref{align:compatibility}). By \cite[Proposition~3.7]{BoschLutke-DAV}, $t (u_1' , \lambda (u_2')) = t (u_2' , \lambda (u_1'))$ as classical points of $\rest{L^{\otimes -1}}{\Phi (u_1') + \Phi (u_2')} \otimes \rest{L}{\Phi (u_1')} \otimes \rest{L}{\Phi (u_2')}$. 
\end{Remark}

Let $c : M' \to \Phi^{\ast} (L^{\otimes -1})$ be a trivialization as a rigidified line bundle on $M'$. Then
\[
c(u'_1 + u'_2) \otimes c(u'_1)^{\otimes -1} \otimes c(u'_2)^{\otimes -1}
\in
\left( \rest{L^{\otimes -1}}{\Phi (u_1') + \Phi (u_2')} \otimes \rest{L}{\Phi (u_1')} \otimes \rest{L}{\Phi (u_2')} \right) (k)
.
\]
Now, we come to define the notion of descent data over $B$.

\begin{Definition} \label{def:NADD}
We call a triple $(L, \lambda , c)$ as follows a \emph{descent datum} over $B$: $L$ is a rigidified line bundle on $B$, $\lambda : M' \to M$ is a homomorphism, and $c : M' \to \Phi^{\ast} (L^{\otimes -1})$ is a trivialization as a rigidified line bundle on $M'$ such that
$\Phi' \circ \lambda = \varphi_L^- \circ \Phi$ and
\begin{align} \label{align:lambda-c}
c(u'_1 + u'_2) \otimes c(u'_1)^{\otimes -1} \otimes c(u'_2)^{\otimes -1} = t (u'_1 , \lambda(u'_2))
\end{align}
for any $u'_1 , u'_2 \in M'$. 
\end{Definition}

\begin{Remark} \label{remark:c0}
Since $t (u'_1 , \lambda(u'_2)) = 1$, it follows from (\ref{align:lambda-c}) that $c(0) = 1 \in k = \rest{L^{\otimes -1}}{0}$, where the equality is the rigidification of $L^{\otimes -1}$.
\end{Remark}

Let $(L,\lambda,c)$ be a descent datum over $B$. From $(L,\lambda , c)$, one constructs an $M'$-linearlization of $q^{\ast} (L)$. Since $\Phi' \circ \lambda = \varphi_L^- \circ \Phi$, equality (\ref{align:E-poincare}) implies $\rest{P}{B \times \{ \varphi_L^- (\Phi (u')) \}} = E_{\lambda (u')}$ for any $u' \in M$, and thus 
\begin{align} \label{align:translate-canonicalisom}
T_{\widetilde{\Phi} (u')}^{\ast} (q^{\ast} (L)) = L(\Phi (u')) \otimes q^{\ast} (E_{\lambda (u')}^{\otimes -1}) \otimes  q^{\ast} (L)
.
\end{align}
Since this isomorphism comes from the pullback of the isometry on (\ref{align:translate-canonicalisom0}), it is also isometry with respect to the metrics arising from the canonical metrics. Then the $M'$-linearlization is characterized by the following data: for any $u' \in M'$, an isomorphism 
\[
q^{\ast} (L) \to T_{\widetilde{\Phi}(u')} (q^{\ast}(L)) ;
\quad
s \mapsto c(u')^{\otimes -1} \otimes e_{\lambda (u')}^{\otimes -1} \otimes s 
\]
under the canonical identification (\ref{align:translate-canonicalisom}). Here, we identify a line bundle with the corresponding invertible sheaf in a covariant way and $s$ is a local section of this sheaf. By the descent theory (\cite[Theorem~6.7]{BoschLutke-DAV}), such an $M'$-linearlized rigidified line bundle descends to a line bundle on $A^{\an} = E^{\an} / M'$. By GAGA, there is a unique rigidified line bundle
\[
\tilde{L}_A(\lambda ,c)
\] 
on $A$ whose analytification equals this line bundle on $A^{\an}$.

Note that, since $B$ has good reduction, so do $B'$ and $B \times B'$. Let $\| \ndot \|$ be the canonical metric on the rigidified Poincar\'e bundle $P$.
We define a map $t_{\trop}(\ndot, \ndot) : M' \times M \to \RR$ by
\begin{align} \label{label:fundamentalbilinearform}
t_{\trop}(\ndot, \ndot) := - \log \| t (\cdot , \cdot)\|
.
\end{align}
This is a bilinear form.

\begin{Lemma} \label{lemma:bilinearforms}
With the above notation, for any $(u' , u) \in M' \times M$, we have $t_{\trop} (u' , u) = \langle u , u' \rangle$.
\end{Lemma}

\Proof
We have
$
t_{\trop} (u' , u) = - \log \| e_{u} (\widetilde{\Phi} (u'))\| 
$,
where 
$\| \ndot \|$ is the canonical metric on the rigidified $P^{\an}$. Recall that 
$\val (\tilde{\Phi} (u')) = u'$ for any $u' \in M'$. Then we see from  \cite[Remark~3.12]{FRSS} and the definition of the valuation map that $- \log \| e_{u} (\widetilde{\Phi} (u'))\| = \langle u , u' \rangle$. Thus we have the desired equality. 
\QED

Let $Q_{\lambda}$ denote the $\RR$-bilinear form on $N_{^\RR} \times N_{\RR}$ corresponding to $\lambda$ by Lemma~\ref{lem:polarization:equiv}. 
For a tropical descent datum $(L,\lambda ,c)$, we call $\lambda$ a \emph{polarization on the torus part (of $A$)} if $\lambda$ is a polarization on the real torus $N_{\RR}/M'$, i.e., $Q_{\lambda}$ is positive-definite. We define the type of $\lambda$ to be the type of $\lambda$ as a polarization on $N_{\RR}/M'$.

By \cite[Theorems~6.13]{BoschLutke-DAV}, we have the following nonarchimedean Appell--Humbert theorem.

\begin{Theorem} [Theorems~6.13 in \cite{BoschLutke-DAV}] \label{thm:NAH}
The analytification of any line bundle $\tilde{L}$ on $A$ is given by the above procedure for some descent datum $(L,\lambda, c)$, i.e., $\tilde{L} \cong \tilde{L}_A (\lambda , c)$. For descent data $(L_1 , \lambda_1 , c_1)$ and $(L_2 , \lambda_2 , c_2)$ over $B$, $\tilde{L_1}_A (\lambda_1 , c_1) \cong \tilde{L_2}_A (\lambda_2 , c_2)$ if and only if $\lambda_1 = \lambda_2$ and there exists a $u \in M$ such that $L_2 = L_1 \otimes E_{u}$ and $c_2^{\otimes -1} = c_1^{\otimes -1} \otimes \Phi^{\ast} (q^{\ast} (e_u))$. Further, $\tilde{L} \cong \tilde{L}_A (\lambda , c)$ is ample if and only if $L$ is ample and $\lambda$ is a polarization on the torus part.
\end{Theorem}

Let $\tilde{L}$ be an ample line bundle. By Theorem~\ref{thm:NAH}, there exists a descent datum $(L, \lambda , c)$ over $B$ such that $\tilde{L} \cong \tilde{L}_A (\lambda , c)$. Furthermore, $\lambda$ is uniquely determined from $\tilde{L}$ and is a polarization on the torus part. We call the type of $\lambda$ the \emph{polarization type of $\tilde{L}$ on the torus part}.

\begin{Remark} \label{remark:equiv:two:polarizations}
The canonical tropicalization of $A$ is a tropical abelian variety. Indeed, since $A$ is projective, there exists an ample line bundle, and by Theorem~\ref{thm:NAH}, it follows that there exists a polarization $\lambda : M' \to M$ on the torus part. By the definitions, this is also a polarization on the tropicalization $N_{\RR}/M'$ of $A$.
\end{Remark}

\setcounter{equation}{0}
\section{Faithful tropicalization of the canonical skeleton}
\label{sec:FT:skeleton}

In this section, we state the results on the faithful tropicalization of the canonical skeleton, which is one of the main results of this paper.

\subsection{Canonical skeleton}

Let $A$ be an abelian variety over $k$. We keep the notation in the previous section. 
We  recall the canonical skeleton in $E^{\an}$ and that in $A^{\an}$. In \cite[Example~7.2]{Gubler-NACM}, Gubler describes the canonical skeleton $\Sigma_{E^{\an}}$ of $E^{\an}$. He constructs a continuous section $\sigma_{E^{\an}} : N_{\RR} \to E^{\an}$ of  $\val_{E_{\an}} : E^{\an} \to N_{\RR}$, called the canonical section, and defined $\Sigma_{E^{\an}}$ to be the image of $\sigma_{E^{\an}}$. We describe this section $\sigma_{E^{\an}}$ by using a formal $T^{\an}_1$-trivialization. 
Let $\{ (V_{\lambda})_{\lambda} , (\rho_{\lambda} : T^{\an} \times V_{\lambda} \to q^{-1} (V_{\lambda}))_{\lambda} \}$ be a formal $T^{\an}_1$-trivialization of the Raynaud extension. We fix any $\lambda$. Then it follows from Gubler's construction that $\Sigma_{E^{\an}} \subseteq q^{-1} (V_{\lambda})$. Let $g$ be an analytic function on $q^{-1} (V_{\lambda})$. We identify $q^{-1} (V_{\lambda})$ with $T^{\an} \times V_{\lambda}$ via $\rho_{\lambda}$. If we regard $g$ as a function on $T^{\an} \times V_{\lambda}$ via this identification, we have a Laurent series expansion
\[
g = \sum_{u \in M} g_u \chi^u
,
\]
where $g_u \in \OO_{V_{\lambda}} (V_{\lambda})$. Let $\xi \in  B^{\an}$ be the Shilov point associated to the model $\mathscr{B}^{\an}$. 
Since $V_{\lambda}$ is the analytification of an affine formal open subset of $\mathscr{B}$, we have $\xi \in V_{\lambda}$. We take any $v \in N_{\RR}$. Then, $\sigma_{E^{\an}} (v) \in \Sigma_{E^{\an}}$ is characterized by the condition that, for any $g$,
\[
|g (\sigma_{E^{\an}} (v))| =\left| \left( \sum_{u \in M} g_u \chi^{u} \right) (\sigma_{E^{\an}} (v)) \right| = \max_{u \in M} |g_u (\xi)| e^{- \langle u , v \rangle}
.
\] 
Since $\xi$ is also the Shilov point of $V_{\lambda}$ corresponding the affine formal open subset of $\mathscr{B}$, it corresponds to the supremum norm $|\ndot|_{V_{\lambda} , \sup}$ on $V_{\lambda}$. Thus
\begin{align} \label{align:valueskeleton1}
|g (\sigma_{E^{\an}} (v))| =\max_{u \in M} |g_u|_{V_{\lambda} , \sup} e^{- \langle u , v \rangle}
.
\end{align}

Now, we construct the canonical skeleton of $A^{\an}$. Recall that we have a  homomorphism $p : E^{\an} \to A^{\an}$ with kernel $M' \cong \Ker (p) \subseteq E(k)$. By \cite[Lemma~4.2]{FRSS}, 
$M'$ acts on the canonical skeleton of $E^{\an}$, and we have a closed subset $\Sigma = \Sigma_{E^{\an}} / M'$ in $A^{\an}$. We call $\Sigma$ the \emph{canonical skeleton} of $A^{\an}$. Note that the canonical section $\sigma_{E^{\an}} :  N_{\RR} \to E^{\an}$ of $\val_{E^{\an}} : E^{\an} \to N_{\RR}$ descends to a section $\sigma : N_{\RR}/M' \to A^{\an}$ of the valuation map $\val_{A^{\an}} : A^{\an} \to N_{\RR}/M'$ of $A^{\an}$, which we call the \emph{canonical section} of $\val_{A^{\an}}$.

\begin{Remark} \label{remark:skeleton-birationalpoint}
Note that for any proper closed subvariety $Z \subseteq A$, $\Sigma_{A^{\an}} \cap Z^{\an} = \emptyset$. Indeed, it follows from the description of $\Sigma_{A^{\an}}$ in \cite[Example~7.2]{Gubler-NACM} that for any $v \in N_{\RR}$, any nonzero rational function on $A$ that is regular at $\sigma (v)$ does not vanish at this point, which means that this point cannot be contained in $Z^{\an}$ for any proper closed subset $Z \subseteq A$. 
\end{Remark}

\subsection{Main results for arbitrary abelian varieties}
We are ready to state our main results precisely.
Let $A$ be an abelian variety over $k$ with canonical skeleton $\Sigma$. Let $\varphi : A \dasharrow \PP^m_k$ be a rational map. Let $\varphi^{\trop} : A^{\an} \dasharrow \TT\PP^m$ denote the composite of $\varphi^{\an}$ with the tropicalization map $(\PP^m_k)^{\an} \to \TT\PP^m$. Although $\varphi^{\trop}$ is not necessarily defined on the whole $A^{\an}$, this map gives rise to a continuous map $\rest{\varphi^{\trop}}{\Sigma} : \Sigma \to \TT\PP^m$ by Remark~\ref{remark:skeleton-birationalpoint}.

\begin{Definition}
We keep the above notation.
\begin{enumerate}
\item
We say that $\varphi$ \emph{faithfully (resp. unimodularly) tropicalizes} the canonical skeleton if the map $\rest{\varphi^{\trop}}{\Sigma} : \Sigma \to \TT\PP^m$ is a faithful embedding (resp. is unimodular).
\item
Let $\tilde{L}$ be a line bundle on $A$. We say that $\tilde{L}$ admits a \emph{rational faithful (resp. rational unimodular) tropicalization} of the canonical skeleton if there exist global sections $s_0 , \ldots , s_m \in H^0 (A,\tilde{L}) \setminus \{ 0 \}$ such that the rational map $\varphi_{s_0 , \ldots , s_m} : A \dashrightarrow \PP^m_k$ faithfully (resp. unimodularly) tropicalizes the canonical skeleton. 
\item
We say that \emph{$\tilde{L}$ admits a faithful (resp. unimodular) tropicalization} of the canonical skeleton if there exist global sections $s_0 , \ldots , s_m \in H^0 (A,\tilde{L}) \setminus \{ 0 \}$ such that $\bigcap_{i=0}^m \Supp (\zero (s_i)) = \emptyset$, where %$\mathrm{Bs} (s_i)$ is the zero-locus of $s_i$
$\Supp (\zero (s_i))$ is the support of the zero divisor $\zero (s_i)$ of $s_i$, and such that the morphism $\varphi_{s_0 , \ldots , s_m} : A \dashrightarrow \PP^m_k$ faithfully (resp. unimodularly) tropicalizes the canonical skeleton. 
\end{enumerate}
\end{Definition}

For the notion of faithful tropicalization of arbitrary varieties by rational maps given  by rational functions, see \cite[Definition 9.2]{GRW}. The formulation using line bundles as above is proposed in \cite{KY1,KY2}.

\begin{Remark} \label{remark:RFT+BPF}
With the above notation, suppose that $\tilde{L}$ admits a rational faithful tropicalization and is basepoint free. Then it admits a faithful (resp. unimodular) tropicalization. Indeed, in this setting, we can take $s_0 , \ldots , s_l$ that gives a rational faithful tropicalization, and we can take global sections $s_{l+1} , \ldots , s_m$ such that $\bigcap_{i=0}^m \Supp (\zero  (s_i)) = \emptyset$. Since a part of $s_0 , \ldots  s_m$ has already given rise to a faithful (resp. unimodular) tropicalization of the canonical skeleton, those sections also give a faithful (resp. unimodular) tropicalization; see also \cite[Lemma 9.3]{GRW}.
\end{Remark}

The following is our main theorem on the faithful tropicalization of the canonical skeleton of an abelian variety.

\begin{Theorem} \label{thm:main:FT1}
Let $A$ be an abelian variety over $k$ and let $\tilde{L}$ be an ample line bundle on $A$. Let $(d_1 , \ldots , d_n)$ be the polarization type of $\tilde{L}$ on the torus part, where $n$ is the torus rank of~$A$. 
\begin{enumerate}
\item
Suppose that $d_1 \geq 3$.
Then there exist $s_0 , \ldots , s_m \in H^0 (A, \tilde{L})$ such that the rational map $\varphi_{s_0 , \ldots , s_m} : A \dashrightarrow \PP^m_k$ induces a homeomorphism $\varphi_{s_0 , \ldots , s_m}^{\trop} : \Sigma \to \TT\PP^m$ onto its image, where $\Sigma$ is the canonical skeleton of $A^{\an}$.
\item
Suppose that $d_1 \geq 2 (n-1)!$. Then $\tilde{L}$ admits a rational unimodular tropicalization of the canonical skeleton. Further, if $\tilde{L}$ is basepoint free, then $\tilde{L}$ admits a unimodular tropicalization of the canonical skeleton.
\item
Suppose that $d_1 \geq \max \{ 3 , 2 (n-1)! \}$. Then $\tilde{L}$ admits a rational faithful embedding of the canonical skeleton. Further, if $\tilde{L}$ is basepoint free, then $\tilde{L}$ admits a faithful tropicalization of the canonical skeleton.
\end{enumerate}
\end{Theorem}

The basepoint freeness assumption in (2) and (3) is necessary, because we impose the condition $\bigcap_i {\rm Bs}(s_i) = \emptyset$ on the definition of a faithful tropicalization but the ample line bundle $\tilde{L}$ is not basepoint free in general.

Note that in the theorem, (3) is an immediate consequence of (1) and (2) with Remark~\ref{remark:RFT+BPF}. We will give a proof of (1) and (2) in the next section.

In the rest of this section, we discuss consequences of Theorem~\ref{thm:main:FT1}. First, note that Theorem~\ref{thm:intro:main:FT} is a corollary of Theorem~\ref{thm:main:FT1}; see (2) in the following corollary.

\begin{Corollary} \label{cor:main:FT1}
Let $A$ be an abelian variety over $k$ and let $\tilde{L}$ be an ample line bundle on $A$. Let $d$ be an integer.
\begin{enumerate}
\item
Suppose that $d \geq 3$.
Then there exist $s_0 , \ldots , s_m \in H^0 (A, \tilde{L}^{\otimes d}) \setminus \{ 0 \}$ such that $\bigcap_{i=0}^m \Supp (\zero (s_i)) = \emptyset$ and the morphism $\varphi_{s_0 , \ldots , s_m} : A \to \PP^m_k$ induces a homeomorphism $\varphi_{s_0 , \ldots , s_m}^{\trop} : \Sigma \to \TT\PP^m$ onto its image.
\item
Suppose that $d \geq \max \{ 3 , 2 (n-1)! \}$ (resp. $d \geq 2 (n-1)!$). Then $\tilde{L}$ admits a faithful (resp. unimodular) tropicalization of the canonical skeleton. 
\end{enumerate}
\end{Corollary}

\Proof
Let $(d_1 , \ldots , d_n)$ be the polarization type of $\tilde{L}$ on the torus part. Then $\tilde{L}^{\otimes d}$ has polarization type $(dd_1 , \ldots , dd_n)$. Further, since $\tilde{L}$ is ample, $\tilde{L}^{\otimes d}$ is basepoint free for $d \geq 2$. Then we get the corollary from Theorem~\ref{thm:main:FT1}.
\QED

\subsection{Totally degenerate case}

In this subsection, we consider the case where $A$ is totally generate. In this case, the abelian part of the Raynaud extension of $A$ is trivial, so that the descent datum is a pair $(\lambda ,c)$ consists of a homomorphism $\lambda : M' \to M$ and $c : M' \to k^{\times}$ satisfying the equality in Definition~\ref{def:NADD}. Note that the compatibility (\ref{align:compatibility}) is trivially satisfied. We call the polarization type on the torus part the polarization type, because the abelian part of $A$ is trivial.

\begin{Lemma} \label{lemma:lambda-c:TD}
Let $\lambda : M' \to M$ be a homomorphism. Let $s : M' \times M' \to k^{\times}$ be a symmetric bilinear form (as $\ZZ$-modules) such that $s(u_1' , u_2')^2 = t (u_1' , \lambda (u_2'))$. Then the following hold.
\begin{enumerate}
\item
Let $c : M' \to k^{\times}$ be a map. 
Then the pair $(\lambda , c)$ is a descent datum if and only if the map
$h : M' \to k^{\times}$ given by $h : u' \mapsto c (u')/s(u',u')$ is a group homomorphism.
\item
Let $l : M' \to k^{\times}$ be a homomorphism. Let $c_{l} : M' \to k^{\times}$ be the map defined by $c_{l} (u') := l (u') s (u',u')$ for any $u' \in M'$. Then $(\lambda , c_{l})$ is a descent datum.
\end{enumerate}
\end{Lemma}

\Proof
We prove (1). 
By (\ref{align:compatibility}), the condition that $(\lambda , c)$ is a descent datum is equivalent to saying that
\[
c (u_1' + u_2') c (u_1')^{-1} c (u_2)^{-1} = s (u_1' , u_2')^2 
\]
for any $u_1' , u_2' \in M'$. Since $s$ is a symmetric bilinear form, we have $s(u_1' + u_2' , u_1' + u_2') = s(u_1',u_1') s(u_2',u_2') s (u_1' ,u_2')^2$, which means the above equality is equivalent to saying that $h$ is a homomorphism. Thus (1) holds.

Assertion (2) immediately follows from (1). Thus the proof of the lemma completes.
\QED

\begin{Proposition} \label{prop:divisible:linebundle:TD}
Let $A$ be an abelian variety over $k$ of dimension $n$ and assume that it is totally degenerate. Let $\tilde{L}$ be an ample line bundle on $A$ of polarization type $(d_1 , \ldots , d_n)$. Then there exists an ample line bundle $\tilde{L}_1$ on $A$ such that $\tilde{L} \cong \tilde{L}_1 ^{\otimes d_1}$.
\end{Proposition}

\Proof
By Theorem~\ref{thm:NAH}, there exists a descent datum $(\lambda , c)$ corresponding to $\tilde{L}$. Since $\frac{1}{d_1} \lambda (M') \subseteq M$, we define a homomorphism $\lambda_1 : M' \to M$ to be $\lambda / d_1$. Since the compatibility (\ref{align:compatibility}) is trivially satisfied for $\lambda_1$, we have $t (u_1' , \lambda_1 (u_2')) = t (u_2 , \lambda_1 (u_1'))$ for all $u_1' , u_2' \in M'$, as is noted in Remark~\ref{remark:tissymmetric}.

We construct a map $c_1 : M' \to k^{\times}$ such that $(\lambda_1 , c_1)$ is a descent datum and $c_1^{d_1} = c$. We fix a basis $m_1' , \ldots , m_n'$ of $M'$. For each $i=1 , \ldots , n$, we take a $c (m_i')^{1/d_1} \in k^{\times}$ and set $c_1 (m_i') := c (m_i')^{1/d_1}$. For any $i,j = 1 , \ldots , n$, we take $\alpha_{ij} \in k^{\times}$ in such a way that $\alpha_{ij}^2 = t (m_i' , \lambda_1 (m_j'))$ and $\alpha_{ij} = \alpha_{ji}$; this is possible because the bilinear form $t (\ndot , \lambda_1 (\ndot))$ on $M' \times M'$ is symmetric. Then, since $m_1' , \ldots , m_n'$ is a basis of $M'$, there exists a symmetric bilinear form $\sqrt{t (\ndot , \lambda (\ndot))} : M' \times M' \to k^{\times}$ such that $\sqrt{t (m_i' , \lambda_1 (m_j'))} = \alpha_{ij}$ for all $i,j$ and such that $\sqrt{t (u_1' , \lambda_1 (u_2'))}^2 = t (u_1' , \lambda_1 (u_2'))$ for all $u_1' , u_2' \in M'$. 
The bilinear form $\sqrt{t (\ndot , \lambda_1 (\ndot))}^{d_1}$ is symmetric. Further, Since $d_1 \lambda_1 = \lambda$, we note that $\left( \sqrt{t (\ndot , \lambda_1 (\ndot))}^{d_1} \right)^2 = t (\ndot , \lambda (\ndot))$. Then, since $(\lambda , c)$ is a descent datum, the map
$h : M' \to k^{\times}$ defined by $h : u' \mapsto c(u') / \sqrt{t (u' , \lambda_1 (u'))}^{d_1} $ is a homomorphism by Lemma~\ref{lemma:lambda-c:TD}~(1).

For any $i=1 , \ldots , n$, we fix a $c(m_i')^{1/d_1} \in k$. Then, since $m_1' , \ldots  , m_n'$ is a basis of the free $\ZZ$-module $M'$, there exists a unique homomorphism $l : M' \to k^{\times}$ such that $l (m_i') =  c(m_i')^{1/d_1}/ \sqrt{t (m_i' , \lambda_1 (m_i'))}$. Define $c_1 : M' \to k$ by $c_1(u') := l (u' )\sqrt{t (u' , \lambda_1 (u'))}$. Then, since the bilinear form $\sqrt{t (\ndot , \lambda_1 (\ndot))}$ on $M' \times M'$ is symmetric and $\sqrt{t (\ndot , \lambda_1 (\ndot))}^{2} = t (\ndot , \lambda_1 (\ndot))$, $(\lambda_1 , c_1)$ is a descent datum by Lemma~\ref{lemma:lambda-c:TD}~(2). Further, for any $i=1, \ldots , n$, we have
\[
 l(m_i')^{d_1} = (c(m_i')^{1/d_1})^{d_1} / \sqrt{t (m_i' , \lambda_1 (m_i'))}^{d_1} = c(m_i') / \sqrt{t (m_i' , \lambda_1 (m_i'))}^{d_1} = h (m_i').
\]
Since $l^{d_1}$ and $h$ are homomorphism, it follows that $l^{d_1} = h$ as maps. Thus, for any $u' \in M'$, we have
\[
c_1 (u')^{d_1} = l (u')^{d_1}  \sqrt{t (u' , \lambda_1 (u'))}^{d_1} = h(u') \sqrt{t (u' , \lambda_1 (u'))}^{d_1} = c (u')
,
\]
which proves that $(d_1 \lambda , c_1^{d_1}) = (\lambda, c)$. 

Let $\tilde{L}_1$ be the line bundle on $A$ corresponding to the descent datum $(d_1 \lambda , c_1^{d_1})$. Then, since $(d_1 \lambda , c_1^{d_1}) = (\lambda, c)$, we have $\tilde{L}_1^{\otimes d_1} = \tilde{L}$. Further, since $\tilde{L}$ is ample, so is $\tilde{L}_1$. This completes the proof. 
\QED

The following is a consequence of Corollary~\ref{cor:main:FT1}~(2) and Proposition~\ref{prop:divisible:linebundle:TD}.

\begin{Corollary} \label{cor:main:FT3}
Let $A$ be an abelian variety of dimension $n$ and assume that it is totally degenerate. Let $\tilde{L}$ be an ample line bundle on $A$ with polarization type $(d_1 , \ldots , d_n)$. Then, if $d_1 \geq \max \{ 3 , 2 (n-1)! \}$ (resp. $d_1 \geq 2 (n-1)!$), $\tilde{L}$ admits a faithful (resp. unimodular) tropicalization of the canonical skeleton.
\end{Corollary}

\Proof
By Proposition~\ref{prop:divisible:linebundle:TD}, there exists an ample line bundle $\tilde{L}_1$ on $A$ such that $\tilde{L} \cong \tilde{L}_1^{\otimes d_1}$. If $d_1 \geq 2$, then $\tilde{L}$ is basepoint free by the Lefschetz theorem. Thus the corollary follows from Corollary~\ref{cor:main:FT1}~(2)
\QED

\setcounter{equation}{0}
\section{Construction of nonarchimedean theta functions}
\label{section:construction:NA:theta}

In this section, we construct nonarchimedean theta functions that are lifts of the tropical theta functions as in (\ref{eqn:def:theta:b0}). Using such nonarchimedean theta functions together with the theorems of faithful embeddings and unimodular embeddings for the tropical abelian varieties, respectively, we show the faithful and unimodular tropicalization theorems, respectively.

We fix the notation. Let $A$ be an abelian variety over $k$. We keep the notation of Raynaud cross in (\ref{align:RaynaudCross}); $q : E \to B$, $T = \Spec (k [\chi^M])$, $p : E^{\an} \to A^{\an}$, and $\widetilde{\Phi} : M' \to E(k)$ are as there. Let $\Phi : M' \to B(k)$ be the map given by $\Phi := q \circ \widetilde{\Phi}$. Let $P \to B \times B'$ be the rigidified Poincar\'e line bundle and let $t : M' \times M \to P(k)$ be as in (\ref{align:def:t}). 
Let $\Phi' : M \to B'(k)$ be the homomorphism corresponding to the Raynaud extension of $A$ in Remark~\ref{remark:hom-algebraicextension}.  Let $\langle \ndot , \ndot \rangle : M \times N_{\RR} \to \RR$ be the canonical pairing.

\subsection{Tropicalization of nonarchimedean theta functions}
\label{subsec:trop:na:theta:functions}

Let $(L , \lambda , c)$ be a nonarchimedean descent datum over $B$.
Let $\tilde{L} = \tilde{L}_A (\lambda, c)$ be the rigidified line bundle on $A$ arising from $(L , \lambda , c)$. Then there is a canonical isomorphism between the pullback of $p^{\ast} (\tilde{L}^{\an})$ by $p : E^{\an} \to A^{\an}$ and $q^{\ast}(L)^{\an}$. Thus we have an injection $p^{\ast} : H^0 (A , \tilde{L}) \hookrightarrow H^0 (E^{\an}, q^{\ast} (L)^{\an})$.

By \cite[Proposition~4.9]{BoschLutke-DAV}, we know that for any $f \in H^0 (E^{\an}, q^{\ast} (L)^{\an})$, we have $f \in p^{\ast} ( H^0 (A , \tilde{L}))$ if and only if $T_{\Phi (u')}^{\ast} (f) = c(u')^{\otimes -1}  \otimes e_{\lambda (u')}^{\otimes -1} \otimes  f$
for any $u' \in M'$, where we use the identification (\ref{align:translate-canonicalisom}).
An element of $p^{\ast} (H^0 (A , \tilde{L})) \subseteq H^0 (E^{\an}, q^{\ast} (L)^{\an})$ that satisfies the above equivalent conditions is called a \emph{nonarchimedean theta function} with respect to $(L, \lambda , c)$.

We explain the tropicalization of a nonarchimedean theta function given in \cite{FRSS}. First, we explain the tropicalizations of descent data. We define $c_{\trop} : M' \to \RR$ by
\[
c_{\trop} (u') = - \log \| c (u') \|
,
\]
where $\| \ndot \|$ is the canonical metric on $L^{\otimes -1}$. By condition (\ref{align:lambda-c}) for $c$ together with Lemma~\ref{lemma:bilinearforms}, we have
\[
c_{\trop} (u'_1 + u'_2) - c_{\trop} (u'_1) - c_{\trop} (u'_2)
= \langle  \lambda (u'_2) , u_1'\rangle
.
\]
Thus $(\lambda , c_{\trop})$ is a tropical descent datum. We call this the \emph{tropicalization of the descent datum $(L, \lambda , c)$}. 

Let $\tilde{\vartheta}$ be a nontrivial nonarchimedean theta function with respect to $(L, \lambda ,c)$. In \cite[Definition~4.6]{FRSS}, they define a function $\tilde{\vartheta}_{\trop} : N_{\RR} \to \RR$ by
\[
\tilde{\vartheta}_{\trop}(v) := - \log \| \tilde{\vartheta} (\sigma_{E^{\an}} (v))\|
,
\]
where $\sigma_{E^{\an}} : N_{\RR} \to E^{\an}$ is the canonical section and $\| \ndot \|$ is the pullback by $q : E \to B$ of the canonical metric on $L$. Further, they show in \cite[Theorem~4.7]{FRSS} that $\tilde{\vartheta}_{\trop}$ is a tropical theta function with respect to $(\lambda , c_{\trop})$.

\subsection{Lifts of the tropical theta functions}
\label{subsection:construct:lift}

In this subsection, let $(L , \lambda ,c)$ be a nonarchimedean descent datum over $B$.
For a tropical theta function $\vartheta$ with respect to the tropicalization $(\lambda , c_{\trop})$ of $(L,\lambda,c)$, we call a nonarchimedean theta function $\tilde{\vartheta}$ with respect to $(L,\lambda,c)$ such that $\tilde{\vartheta}_{\trop} = \vartheta$ a \emph{lift} of $\vartheta$. 

The goal of this subsection is to construct a lift of $\vartheta_b^{(\lambda ,c_{\trop})}$ defined in (\ref{eqn:def:theta:b0}) for $(\lambda , c_{\trop})$.
By Theorem~\ref{thm:generators0}, (\ref{eqn:def:theta:b2}), and (\ref{eqn:def:theta:b0}), any tropical theta function can be written as 
\[
\min_{u \in M} ( \langle u , x \rangle + c_u ),
\] 
where $c_u \in \TT$. We call this expression the \emph{(tropical) Fourier expansion} of the tropical theta function.
Since the tropical theta function $\vartheta_b^{(\lambda ,c_{\trop})}$ in (\ref{eqn:def:theta:b0}) is given by tropical Fourier expansion, we will use the nonarchimedean Fourier expansion to construct a lift of it, taking into account of  \cite[Theorem~4.7~(1)]{FRSS}, which is cited as Proposition~\ref{prop:tropicalization_FE} below. 

\subsubsection{Fourier expansions}

For any $u \in M$, we regard 
$H^0 (B , L \otimes E_u^{\otimes -1} ) \otimes e_u \subseteq H^0 (E^{\an} , q^{\ast} (L^{\an}))$
via the pullback by $q: E \to B$. We fix a basis $v_1 , \ldots , v_n$ of $N$, where $n:=\rank (N)$. For any $u$, we set $|u| := \sum_{i=1}^n |\langle u, v_i \rangle|$.
By the first paragraph of Section~5 in \cite{BoschLutke-DAV} (see also \cite[Remark~3.12]{FRSS}),  the following Fourier decomposition holds:
\[
H^0 (E^{\an} , q^{\ast} (L^{\an})) = \hat{\bigoplus}_{u \in M} H^0 ( B , L \otimes E_{u}^{\otimes -1}) \otimes e_u
.
\]
Here, the right-hand side means 
\[
\left\{\left.
\sum_{u \in M} g_u  \otimes e_u \;\right|\; \text{$g_u \in H^0 (  B , L \otimes E_u^{\otimes -1} )$, $\lim_{|u| \to \infty} - \log \|g_u\|_{\sup} + \langle u ,v \rangle = + \infty$ for any $v \in N_{\RR}$}
\right\}
,
\]
where $\|g_u\|_{\sup}$ is the supremum norm of $g_u$ with respect to the canonical metric on $L \otimes E_u^{\otimes -1}$. This does not depend on the choice of a basis $v_1 , \ldots , v_n$.

When we argue the convergence of formal series of Fourier expansions, the following lemma will be useful.

\begin{Lemma} \label{lemma:criterion-convergence}
The following are equivalent to each other:
\begin{enumerate}
\item[(i)]
for any $v \in N_{\RR}$, $\lim_{|u| \to  \infty} - \log \|g_u\|_{\sup} + \langle u ,v \rangle = + \infty$;
\item[(ii)]
for any $v \in N_{\RR}$, the subset $\{ - \log \|g_u\|_{\sup} + \langle u , v \rangle \mid u \in M \}$ of $\RR$ has a minimum.
\end{enumerate}
\end{Lemma}

\Proof
It is obvious that (i) implies (ii). To prove that (ii) implies (i), we show the contraposition. Assume that (i) does not hold. Then there exist a $v \in N_{\RR}$, a sequence $(u_l)_{l}$ on $M$, and a $C \in \RR \cup \{ - \infty \}$ such that 
\begin{align} \label{align:lemma:criterion-convergence}
\begin{split} 
&\lim_{l \to \infty} |u_l| = + \infty , \\
&\lim_{l \to \infty} - \log \|g_{u_l} \|_{\sup} + \langle u_l ,v \rangle = C.
\end{split}
\end{align}
Recall that we have fixed a basis $v_1 , \ldots , v_n$ of $N_{\RR}$ and that $|u| = \sum_{i=1}^n |\langle u , v_i\rangle|$ for any $u \in M$. Since $\lim_{l \to \infty} |u_l| = + \infty$, there exists an $i = 1 , \ldots , n$ such that, replacing $(u_l)_l$ by a subsequence if necessary, we have $\lim_{l \to \infty} \langle u_l ,v_i \rangle = + \infty$ or $\lim_{l \to \infty} \langle u_l ,v_i \rangle = - \infty$. We write $\lim_{l \to \infty} \langle u_l , v_i \rangle = \epsilon \infty$, where $\epsilon = \pm 1$. 
By (\ref{align:lemma:criterion-convergence}), we have
\[
\lim_{l \to \infty} - \log \|g_{u_l}\|_{\sup} + \langle u_l ,\lambda v - \epsilon v_i \rangle  
= \lim_{l \to \infty} - \log \|g_{u_l}\|_{\sup} + \langle u_l , v \rangle  - \epsilon \langle u_l , v_i \rangle =
- \infty.
\]
This shows that (ii) does not hold. Thus the proof completes. 
\QED

The following proposition shows how the Fourier expansion of a nonarchimedean theta function is tropicalized to a tropical Fourier expansion. 

\begin{Proposition} [Theorem~4.7~(1) in \cite{FRSS}] \label{prop:tropicalization_FE}
Let $\tilde{\vartheta} = \sum_{u \in M} g_u \otimes e_{u}$ be the Fourier expansion of a nonarchimedean theta function $\tilde{\vartheta}$. Then,
\[
\tilde{\vartheta}_{\trop} (v) = \min_{u \in M} ( \langle u , v \rangle - \log \| g_u \|_{\sup} )
\]
for any $v \in N_{\RR}$, where $\| \ndot \|_{\sup}$ is the supremum norm arising from the canonical metric on $L \otimes E_u^{\otimes -1}$.
\end{Proposition}

\begin{Remark} \label{remark:tropicalization-inequality}
For any nonarchimedean theta functions $\tilde{\vartheta}$ and $\tilde{\vartheta}'$ with respect to $(L, \lambda , c)$, we have $(\tilde{\vartheta} + \tilde{\vartheta}')_{\trop} \geq \min \{ \tilde{\vartheta}_{\trop} , \tilde{\vartheta}'_{\trop} \}$. Indeed, we write $\tilde{\vartheta} = \sum_{u \in M} g_u \otimes e_{u}$ and $\tilde{\vartheta}' = \sum_{u \in M} g_u' \otimes e_{u}$ as in Proposition~\ref{prop:tropicalization_FE}. Then, since $\| g_u + g_u' \|_{\sup} \leq \max \{ \| g_u \|_{\sup} , \| g_u'\|_{\sup} \}$, we have the desired inequality. 
\end{Remark}

\subsubsection{Construction of lifts}

Recall from the beginning of this section that $\Phi' : M \to B'(k)$ is the homomorphism corresponding to the Raynaud extension of $A$ in Remark~\ref{remark:hom-algebraicextension} and that  $\Phi : M' \to B(k)$ is given by $\Phi := q \circ \widetilde{\Phi}$. Let $\varphi_L^- : B \to B'$ be the homomorphism given by (\ref{align:translate-canonicalisom00}). 

\begin{Lemma} \label{lemma:identitybdagger}
Assume that $b \in M$ and $\tilde{b} \in B(k)$ satisfy $\varphi_L^- (\tilde{b}) = \Phi' (b)$. Then for any $u' \in M'$, we have a canonical identity
\[
\rest{L}{ \tilde{b} + \Phi (u')} = \rest{L} { \tilde{b}} \otimes \rest{L}{\Phi (u')} \otimes \rest{E_{b}^{\otimes -1}}{\Phi (u')}
,
\]
where we identify $E_b$ with the corresponding rigidified line bundle.
\end{Lemma}

\Proof
By (\ref{align:translate-canonicalisom0}), we have a canonical identity
$T_{\tilde{b}}^{\ast} (L)^{\otimes -1} \otimes L  \otimes L(\tilde{b})
= E_{b}$. Taking the fibers of both sides over $\Phi (u')$, we have
$\rest{T_{\tilde{b}}^{\ast} (L)^{\otimes -1}}{\Phi (u')} \otimes \rest{L}{\Phi (u')}  \otimes \rest{L}{\tilde{b}}
= \rest{E_{b}}{\Phi (u')}$.
Since $\rest{T_{\tilde{b}}^{\ast} (L)}{\Phi (u')} = \rest{L}{\tilde{b} + \Phi (u')}$, this proves the desired identity. 
\QED

Assume that $b \in M$ and $\tilde{b} \in B(k)$ satisfy $\varphi_L^- (\tilde{b}) = \Phi' (b)$. We fix a $c_{\tilde{b},0} \in \rest{L^{\otimes -1}}{\tilde{b}} (k)$ with $c_{\tilde{b},0} \neq 0$. For any $u' \in M'$, we set
\[
c_{\tilde{b}} (u') := c_{\tilde{b},0} \otimes c(u') \otimes t(u' , b) \in \left(
 \rest{L^{\otimes -1}}{\tilde{b}} \otimes \rest{L^{\otimes -1}}{\Phi (u')} \otimes \rest{E_{b}}{\Phi (u')} \right) (k)
.
\]
By Lemma~\ref{lemma:identitybdagger}, $c_{\tilde{b}} (u')  \in \rest{L^{\otimes -1} }{\tilde{b} + \Phi (u')} (k)$ naturally. Since $c(0) = 1$ as in Remark~\ref{remark:c0}, we note that $c_{\tilde{b}} (0) = c_{\tilde{b},0}$. Further, it is straightforward to see that
\begin{align} \label{align:quasiperiodcitydagger}
c_{\tilde{b}} (u'_1 + u'_2) = c_{\tilde{b}} (u'_1) \otimes c (u_2') \otimes t (u_2' , b + \lambda (u_1'))
.
\end{align}

Fix an $s \in H^0 (B , L)$. For any $u' \in M'$, set
\begin{align} \label{align:def:s..}
s_{\tilde{b} ,u'} :=c_{\tilde{b}} (u') \otimes T_{\tilde{b} + \Phi (u')}^{\ast} (s)
,
\end{align}
which is a global section of $L^{\otimes -1} (\tilde{b} + \Phi (u')) \otimes T_{ \tilde{b} + \Phi (u')}^{\ast} (L)$.
By (\ref{align:translate-canonicalisom0}) and (\ref{align:E-poincare}), we have
\begin{align} \label{align:translate2}
T_{\tilde{b} + \Phi (u')}^{\ast} (L)^{\otimes -1} \otimes L  \otimes L(\tilde{b} + \Phi (u'))
= E_{b + \lambda(u')}
.
\end{align}
Thus $s_{\tilde{b} ,u'} \in H^0 (B, L \otimes E_{-b - \lambda (u')})$, and $e_{b + \lambda (u')} \otimes s_{\tilde{b} ,u'} \in H^0 (B,L)$.

\begin{Proposition} \label{prop:comparison-metric}
We keep the above notation and assume that $\varphi_L^- (\tilde{b}) = \Phi' (b)$.
Then the following hold.
\begin{enumerate}
\item
We have
\[
- \log \| s_{\tilde{b} ,u'} \|_{L \otimes  E_{-b - \lambda (u')}} =
 c_{\trop} (u') + \langle b ,u' \rangle - \log \| s \|_{L} - \log \| c_{\tilde{b} ,0} \|_{L^{\otimes -1}},
\]
where $\| \ndot \|_{L \otimes  E_{-b - \lambda (u')}}$, $\| \ndot \|_L$, and $\| \ndot \|_{L^{\otimes -1}}$ are the canonical metrics on $L \otimes  E_{-b - \lambda (u')}$, $L$, and $L^{\otimes -1}$, respectively.
\item
Assume in addition that $\lambda : M' \to M$ is a polarization on the torus part of $A$. Then the infinite sum $\sum_{u' \in M'}
e_{b + \lambda (u')} \otimes s_{\tilde{b} ,u'}$ converges in $H^0 (E^{\an} , q^{\ast} (L^{\an}))$. 
\end{enumerate}
\end{Proposition}

\Proof
Let $\| \ndot \|'$ be the metric on $T_{\tilde{b} + \Phi (u')}^{\ast} (L)$ that is the pullback by $T_{\tilde{b} + \Phi (u')}$ of the canonical  metric $\| \ndot \|_L$ on $L$. Since  (\ref{align:translate-canonicalisom0}) holds in the level of formal models and the metrics of the line bundles here are the formal metrics, (\ref{align:translate2}) is an isometry. From that fact, we compute
\[
\| s_{\tilde{b} ,u'} \|_{L \otimes  E_{b + \lambda (u')}} = 
\| c_{\tilde{b}}(u') \| \| T_{\tilde{b} + \Phi (u')}^{\ast} (s) \|'
.
\]
Since $\| \ndot \|'$ is the pullback of $\| \ndot \|_L$ by $T_{\tilde{b} + \Phi (u')}$, we have $\| T_{\tilde{b} + \Phi (u')}^{\ast} (s) \|' = \| s \|_L$, and thus $\| s_{\tilde{b} ,u'} \|_{L \otimes  E_{\tilde{b} + \lambda (u')}} = \| c_{\tilde{b}} (u') \| \| s \|_{L}$. By the definition of $c_{\tilde{b}} (u')$, we have
\[
\| c_{\tilde{b}} (u') \| = \| c_{\tilde{b} ,0}\| \|c(u')\| \| t (u' , b)\| ,
\]
where the metrics are the canonical metrics. Noting Lemma~\ref{lemma:criterion-convergence} and taking $- \log$ prove (1).

To prove (2), it suffices to show that for any $ v \in N_{\RR}$,
\[
\{ - \log \| s_{\tilde{b}, u'} \|_{\sup} + \langle b + \lambda (u') , v \rangle \mid u' \in M' \}
\]
has a minimum. By (1),
\begin{multline} \label{align:theta-computation}
- \log \| s_{\tilde{b} ,u'} \|_{\sup} + \langle b + \lambda (u') , v \rangle \\
= \langle b + \lambda (u') , v \rangle+ c_{\trop} (u') + \langle b ,u'  \rangle - \log \| s \|_{\sup} - \log \| c_{\tilde{b} ,0} \|_{L^{\otimes -1}}
.
\end{multline}
Note that $(\lambda ,c_{\trop})$ is a tropical descent datum on $N_{\RR}$.
Since $\lambda$ is a polarization on the canonical tropicalization $N_{\RR}/M'$ of $A$ (cf, Remark~\ref{remark:equiv:two:polarizations}), it follows from Lemma~\ref{lemma:towardSumiTheta_new1} that
the right-hand side in (\ref{align:theta-computation}) has a minimum. This completes the proof.
\QED

We construct an element in $H^0 (E^{\an} , q^{\ast} (L^{\an}))$ when $\tilde{L}_A(\lambda , c)$ is ample. Assume that $\tilde{L}_A(\lambda , c)$ is ample. Fix any $b \in M$. Since the line bundle $L$ on $B$ is ample by Theorem~\ref{thm:NAH}, $\varphi_L^- : B(k) \to B'(k)$ is surjective, and thus there exists a $\tilde{b} \in B(k)$ such that $\varphi_L^-(\tilde{b}) =\Phi'(b)$. Further, since $\lambda : M' \to M$ is a polarization on the torus part by Theorem~\ref{thm:NAH}, the following definition makes sense by Proposition~\ref{prop:comparison-metric}~(2).

\begin{Definition} \label{def:nonarchimedean:lift}
Assume that the line bundle $\tilde{L}_A(\lambda , c)$ is ample. For $b \in M$ and $\tilde{b} \in B(k)$ as above, we fix a nonzero $c_{\tilde{b},0} \in \rest{L^{\otimes -1}}{\tilde{b}} (k)$ and fix a nonzero $s \in H^0 (B,L)$. Let $s_{\tilde{b},u'}$ is the one given in (\ref{align:def:s..}). Then we set
\[
\tilde{\vartheta}_{\tilde{b},s}
:=
\sum_{u' \in M'}
e_{b + \lambda (u')} \otimes s_{\tilde{b} ,u'}
.
\]
\end{Definition}

We remark that $\tilde{\vartheta}_{\tilde{b},s}$ depends also on the choice of a nonzero $c_{\tilde{b} , 0} \in \rest{L^{\otimes -1}}{\tilde{b}}$, but this plays only a minor role.

If $\tilde{L}_A (\lambda ,c)$ is ample, then $\lambda$ is a polarization on the canonical tropicalization $N_{\RR} /M'$ (cf. Remark~\ref{remark:equiv:two:polarizations}), and thus we can consider $\vartheta_b^{(\lambda , c_{\trop})}$.
The following theorem shows that $\tilde{\vartheta}_{\tilde{b},s}$ is a nonarchimedean theta function and is a lift of $\vartheta_b^{(\lambda , c_{\trop})}$ up to an additive constant.

\begin{Theorem} \label{thm:liftingtheorem}
With the above notation, assume that the line bundle $\tilde{L}(\lambda , c)$ is ample. Then the following hold.
\begin{enumerate}
\item
For any $\widetilde{u}' \in M'$, we have
$T_{\widetilde{\Phi} (\widetilde{u}')}^{\ast} (\tilde{\vartheta}_{\tilde{b},s}) =
\tilde{\vartheta}_{\tilde{b},s} \otimes c(\widetilde{u}')^{\otimes -1} \otimes e_{\lambda (\widetilde{u}')}^{\otimes -1} $.
\item
We have
\[
(\tilde{\vartheta}_{\tilde{b} , s})_{\trop}
=
 \vartheta_b^{(\lambda , c_{\trop})}  
- \log \| s \|_{\sup} - \log \| c_{\tilde{b} ,0} \|_{L^{\otimes -1}}
.
\]
\end{enumerate}
\end{Theorem}

\Proof
We take any $z \in E^{\an}$.
We have
\begin{align*}
T_{\widetilde{\Phi} (\widetilde{u}')}^{\ast} (e_{b + \lambda (u')}) (z) &=
e_{b + \lambda (u')} (z + \widetilde{\Phi} (\widetilde{u}')) \\
&= e_{b + \lambda (u' + \widetilde{u}')} (z) \otimes e_{-\lambda (\widetilde{u}')} (z) \otimes e_{b + \lambda (u')} (\widetilde{\Phi} (\widetilde{u}')) \\
&= 
e_{b + \lambda (u' + \widetilde{u}')} (x) \otimes e_{-\lambda (\widetilde{u}')} (x) \otimes t ( \widetilde{u}' , b + \lambda (u'))
.
\end{align*}
Further, using (\ref{align:quasiperiodcitydagger}), we compute
\begin{align*}
T_{\Phi (\widetilde{u}')}^{\ast} (s_{\tilde{b} ,u'})
&= c_{\tilde{b}} ( u') \otimes T_{\tilde{b} + \Phi (\lambda (u'+\widetilde{u}'))}^{\ast} (s)
\\
&= c_{\tilde{b}} (u' + \widetilde{u}') \otimes c (\widetilde{u}')^{\otimes -1} \otimes t (\widetilde{u}' , b + \lambda (u'))^{\otimes -1} \otimes T_{\tilde{b} + \Phi (\lambda (u'+\widetilde{u}'))}^{\ast} (s)
\\
&=
c (\widetilde{u}')^{\otimes -1} \otimes t (\widetilde{u}' , b + \lambda (u'))^{-1} \otimes s_{\tilde{b} ,u' + \widetilde{u}'}
.
\end{align*}
Thus the equality in (1) holds.

We take any $v \in N_{\RR}$. Note that $\| s \|_{\sup} = \| s (\xi )\|$, where $\xi \in B^{\an}$ is the Shilov point corresponding the formal model $\mathscr{B}$ of $B^{\an}$. Let $\sigma_{E^{\an}} : N_{\RR} \to E^{\an}$ be the canonical section of $\val_{E^{\an}} : E^{\an} \to N_{\RR}$. Then by Proposition~\ref{prop:tropicalization_FE} and Proposition~\ref{prop:comparison-metric}~(1),
\begin{align*}
(\tilde{\vartheta}_{\tilde{b} , s})_{\trop} (v) 
&= 
- \log \| \tilde{\vartheta}_{b,s} (\sigma_{E^{\an}} (v)) \| \\
&= \min_{u' \in M'}
( - \log \|  s_{-b - \lambda (u') }  \|_{\sup}  + \langle b + \lambda (u') , v \rangle )
\\
&= 
\min_{u' \in M'}
( \langle b + \lambda (u') , v \rangle + c_{\trop} (u') + \langle b ,u' \rangle ) - \log \| s \|_{\sup} - \log \| c_{\tilde{b} ,0} \|_{L^{\otimes -1}}
\\
&=
\vartheta_b^{(\lambda , c_{\trop})} (v) - \log \| s \|_{\sup} - \log \| c_{\tilde{b} ,0} \|_{L^{\otimes -1}}
.
\end{align*}
Thus the (2) holds. 
\QED

\begin{Corollary} \label{cor:preciselift}
Assume that $\tilde{L}_A (\lambda , c)$ is ample. Let $b \in M$. Then there exists a $\tilde{b} \in B(k)$ such that $\varphi_L^- (\tilde{b}) = \Phi' (b)$, and there exist a $c_{\tilde{b} ,0} \in \rest{L^{\otimes -1}}{\tilde{b}} (k) \setminus \{ 0 \}$ and an $s \in H^0 (B,L) \setminus \{ 0 \}$ such that $\left( \tilde{\vartheta}_{\tilde{b},s} \right)_{\trop} = \vartheta_b^{(\lambda , c_{\trop})}$, where $\tilde{\vartheta}_{\tilde{b},s}$ is defined in Definition~\ref{def:nonarchimedean:lift}.  In particular, for any $b \in M$, there exists a lift of $\vartheta_b^{(\lambda , c_{\trop})}$.
\end{Corollary}

\Proof
Since $\tilde{L}_A (\lambda ,c)$ is ample, $L$ on $B$ is also ample by Theorem~\ref{thm:NAH}. It follows that $\varphi_L^- : B(k) \to B'(k)$ is surjective, and hence there exists a $\tilde{b} \in B(k)$ such that $\varphi_L^- (\tilde{b}) = \Phi' (b)$. Since the metric on $L$ is a formal model metric and $\tilde{b} \in B(k)$, Remark~\ref{remark:value:norm} gives us $c_{\tilde{b} , 0} \in \rest{L^{\otimes -1}}{\tilde{b}} (k)$ such that $\| c_{\tilde{b} , 0} \| = 1$. Since $L$ is ample, it follows that $L$ has a nontrivial global section $s$; this follows from the two theorems at the beginning of Section~16 in \cite{Mumford-AV}. 
By Remark~\ref{remark:value:norm}, $ \| s \|_{\sup} \in |k^{\times}|$. By multiplying a suitable nonzero element of $k$ to $s$, we may assume that $\| s \|_{\sup} = 1$. 
Then the corollary follows from Theorem~\ref{thm:liftingtheorem}~(2).
\QED

\subsection{Proof of Theorem~\ref{thm:main:FT1}}

We keep the notation in Theorem~\ref{thm:main:FT1} and that of the beginning of this section. 
There exists a descent datum $(L , \lambda , c)$ over $B$ such that $\tilde{L} \cong \tilde{L}_A (\lambda , c)$. Let $(\lambda , c_{\trop})$ be the tropicalization of this descent datum. 
Recall that $(d_1 , \ldots , d_n)$ is the type of $\lambda$. Let $\{ b_1 , \ldots , b_D \}$ be a complete system of representatives of $\mathrm{Coker} (\lambda)$, where $D = |\mathrm{Coker} (\lambda)| = d_1 \cdots d_n$. 

Assume that $d_1 \geq 3$. Since $\tilde{L}\cong \tilde{L}_A (\lambda , c)$ is ample, Corollary~\ref{cor:preciselift} gives us  a lift $\tilde{\vartheta}_j$ of $\vartheta_{b_j}^{(\lambda ,c_{\trop})}$ for each $j=1 , \ldots , D$.
Then we have a rational map
\[
\tilde{\varphi}: A \dashrightarrow \PP^{D-1}_k ; \quad
x \mapsto (\tilde{\vartheta}_{1} (x) : \cdots : \tilde{\vartheta}_{D} (x)),
\]
and this induces a 
continuous map $\tilde{\varphi}^{\trop} : \Sigma \to \TT\PP^{D-1}$.

Let $\val_{A^{\an}} : A^{\an} \to N_{\RR} / M'$ be the valuation map for $A^{\an}$. Let $\sigma : N_{\RR} /M' \to A^{\an}$ be the canonical section. 
We define $\varphi : N_{\RR} /M' \to \TT\PP^{D-1}$ by 
$\varphi := (\bar{\vartheta}_{b_1}^{(\lambda,c_{\trop})} : \cdots : \bar{\vartheta}_{b_D}^{(\lambda,c_{\trop})})$, 
where $\bar{\vartheta}_{b_j}^{(\lambda , c_{\trop})}$ indicates that $\vartheta_{b_j}^{(\lambda , c_{\trop})}$ is regarded as a regular global section of the corresponding tropical line bundle on $N_{\RR}/M'$. 
Since
$\sigma^{\ast} \left( - \log \| \tilde{\vartheta}_{j} \| \right) = (\tilde{\vartheta}_{j})_{\trop} = \vartheta_{b_j}^{(\lambda, c_{\trop})}$, 
we have
\begin{align*} 
\rest{\tilde{\varphi}^{\trop}}{\Sigma} \circ \sigma = \varphi
.
\end{align*} 
Then, since $\sigma$ give an isomorphism $N_{\RR}/M' \to \Sigma$ of polyhedral sets, the theorem is an  immediate consequence of Theorem~\ref{thm:main:1-1}. Thus the proof is complete.
\QED

\setcounter{equation}{0}
\section{Surjectivity of the tropicalization for the theta functions}
\label{section:surjective:trop:theta}
 
Let $A$ be an abelian variety over $k$ and  let
\[
\begin{CD}
1 @>>> T @>>> E @>>> B @>>> 0
\end{CD}
\]
be the Raynaud extension for $A$. Let $N$ be the dual of the character lattice of $T$.
Let $(L , \lambda , c)$ be a nonarchimedean descent datum over $B$. Let $(\lambda , c_{\trop})$ be its tropicalization, which is a tropical descent datum on the canonical tropicalization of $A$. Then we have a tropicalization map from the set of nonarchimedean theta functions on $E^{\an}$ with respect $(L , \lambda , c)$ to the set of tropical theta functions on $N_{\RR}$ with respect to the tropicalization $(\lambda , c_{\trop})$ of $(L , \lambda , c)$. In general, this map is not surjective because the value group $\Gamma$ of $k$ is only a subset of $\RR$ and does not equal $\RR$. This indicates that to discuss the surjectivity, we need to restrict ourselves to some class of tropical theta functions, taking into account the value group. 

In this section, we prove that if we restrict the target of the tropicalization map to the set of ``$\Gamma$-rational'' tropical theta functions, then it will turn out that the map is surjective.

\subsection{$\Gamma$-rational tropical theta functions}

In this section, we define the notion of $\Gamma$-rational tropical theta functions for a subgroup $\Gamma$ of $\RR$. 
We fix our setting. Let $M$, $N$, and $M' \hookrightarrow N_{\RR}$ be as in Subsection~\ref{subsection:notation1}. Recall that $n= \dim (N_{\RR})$. Let $\Gamma \subseteq \RR$ be a subgroup. We set $N_{\Gamma} := N \otimes \Gamma$. Note that $N_{\Gamma} \subseteq N_{\RR}$. 
Assume that $M' \subseteq N_{\Gamma}$. 

Recall that in our convention, a regular function on $N_{\RR}$ means a constant function $+ \infty$ or  concave piecewise $\ZZ$-affine function (cf. Remark~\ref{remark:regular-locallypiecewiseaffine}). 
We say a $\ZZ$-affine function $F$ on $N_{\RR}$ is \emph{$\Gamma$-rational} if there exist a $u \in M$ and a $b \in \Gamma$ such that $F = u + b$. Note that $u+b$ with $u \in M$ and $b \in \RR \setminus \Gamma$ can never be $\Gamma$-rational. We say that regular function on $N_{\RR}$ is \emph{$\Gamma$-rational} if all the $\ZZ$-affine pieces of it are $\Gamma$-rational. Note that if $u+ b$ is a $\ZZ$-affine function that equals a $\Gamma$-rational regular function on an $n$-dimensional  polyhedron, then $b \in \Gamma$, because such a $u+ b$ is unique.

Let $\gamma : M' \to \RR$ be a map. We say that $\gamma$ is \emph{$\Gamma$-rational} if $\gamma (M') \subseteq \Gamma$. For any $u \in M$, since $\langle u,u' \rangle \in \Gamma$ for any $u' \in M' \subseteq N_{\Gamma}$, $\gamma + \rest{u}{M'}$ is also $\Gamma$-rational

Let $(\lambda, \gamma)$ be a tropical descent datum on $X:=N_{\RR}/M'$. Throughout  this subsection, we assume that $\gamma$ is $\Gamma$-rational. A tropical theta function with respect to $(\lambda,\gamma)$ is said to be \emph{$\Gamma$-rational} if it is a $\Gamma$-rational regular function. Recall that $L(\lambda ,\gamma)$ denotes the (rigidified) tropical line bundle on $X$ associated to $(\lambda , \gamma)$. In this section, we identify $H^0 (X,L(\lambda , \gamma))$ with the $\TT$-semimodule of tropical theta functions with respect to $(\lambda , \gamma)$. We set $\overline{\Gamma} := \Gamma \cup \{ + \infty \}$. Then the subset of $H^0 (X , L(\lambda , \gamma))$ consisting of $\Gamma$-rational tropical theta functions is a $\overline{\Gamma}$-semimodule, for which we write $H^0 (X,L(\lambda , \gamma))_{\Gamma}$.

\begin{Proposition} \label{prop:Gammarationalbasis}
Under the above setting,  assume that $\lambda$ is a polarization on $X$. Let $\mathfrak{B}$ be a complete system of representatives of $M / \lambda (M')$. For each $b \in \mathfrak{B}$, let $\vartheta_{b}^{(\lambda , \gamma)}$ be the tropical theta function given by (\ref{eqn:def:theta:b0}). Then, as a $\overline{\Gamma}$-semimodule, $\left\{ \vartheta_{b}^{(\lambda , \gamma)} \right\}_{b \in \mathfrak{B}}$ generates $H^0 (X,L(\lambda , \gamma))_{\Gamma}$.
\end{Proposition}

\Proof
We take any $\vartheta \in H^0 (X,L(\lambda , \gamma))_{\Gamma} \setminus  \{ + \infty \}$.
By Theorem~\ref{thm:generators0} and (\ref{eqn:def:theta:b2}), $\left\{ \vartheta_b^{(\lambda, \gamma)} \right\}_{b \in \mathfrak{B}}$ generates $H^0 (X,L(\lambda , \gamma))$ as a $\TT$-semimodule. It follows that for any $b \in \mathfrak{B}$, there exists a $c_b \in \TT$  such that $\vartheta = \min_{b \in \mathfrak{B}} \left( c_b + \vartheta_b^{(\lambda , \gamma)}\right)$. For $b \in \mathfrak{B}$ such that there is no $n$-dimensional polyhedron on $N_{\RR}$ on which $\vartheta$ coincides with $c_b + \vartheta_b^{(\lambda , \gamma)}$, then we may replace $c_b$ by $+ \infty$. Suppose that $\vartheta$ coincides with $c_b + \vartheta_b^{(\lambda , \gamma)}$ on some $n$-dimensional polyhedron. Then, since $c_b + \vartheta_b^{(\lambda , \gamma)}$ is $\Gamma$-rational, we have $c_b \in \Gamma$. This proves the proposition. 
\QED

\subsection{Surjectivity of the tropicalization}

In this subsection, assume that $\Gamma$ is the value group of $k$, i.e., $\Gamma := \{ - \log |\lambda |_k \in \RR \mid  \lambda \in k^{\times} \}$. Let $A$ be an abelian variety over $k$. We consider the Raynaud cross of $A$ and use the notation in (\ref{align:RaynaudCross}). We will prove Theorem~\ref{thm:surjective:tropicalization:theta} below, which is the goal of this section.

We begin with a lemma.

\begin{Lemma} \label{lemma:tropicalization:Gamma-rational}
Let $B$ be the abelian part of an abelian variety $A$ over $k$, and let $(L,\lambda,c)$ be a nonarchimedean descent datum over $B$. Then the tropicalization $(\lambda ,c_{\trop})$ is $\Gamma$-rational. Further, for any $f \in H^0 (A , \tilde{L}_A (\lambda , c)) \setminus \{ 0 \}$, $f_{\trop}$ is a $\Gamma$-rational tropical theta function with respect to $(\lambda , c_{\trop})$. 
\end{Lemma}

\Proof
Let $\| \ndot \|$ be the canonical metric on the line bundle $(L^{\otimes -1})^{\an}$ on $B^{\an}$. Since $\| \ndot \|$ is a formal model metric and $c(u') \in L^{\otimes -1} (k) \setminus \{0 \}$ for any $u' \in M'$, we have $\| c (u') \| \in |k_{\times}|$ and thus $c_{\trop} (u') = - \log \| c (u') \| \in \Gamma$ for any $u' \in M'$. This proves the first assertion of the lemma.

Let $\| \ndot \|'$ be the canonical metric on $L^{\an}$, which is a formal model metric. For a nontrivial global section $g$ of $L$, we have $\| g \|'_{\sup} \in \Gamma$ by the maximal modulus principle. By Proposition~\ref{prop:tropicalization_FE}, we see that the tropicalization of a nonarchimedean theta function is $\Gamma$-rational. This completes the proof.
\QED

Let $X$ be the canonical tropicalization of $A$. Then
Lemma~\ref{lemma:tropicalization:Gamma-rational} says, under the setting of the lemma, that the tropicalization map $H^0 (A , \tilde{L}_A (\lambda ,c)) \to H^0 (X,L (\lambda , c_{\trop}))$ restricts to a map $H^0 (A , \tilde{L}_A (\lambda ,c)) \to H^0 (X,L (\lambda , c_{\trop}))_{\Gamma}$. The following theorem asserts that this map is surjective unless $H^0 (A , \tilde{L}_A (\lambda ,c))$ is trivial.

\begin{Theorem} \label{thm:surjective:tropicalization:theta}
Let $A$ be an abelian variety over $k$ with abelian part $B$. Let $X = N_{\RR}/M'$ be the canonical tropicalization of $A$. Let $(L,\lambda,c)$ be a nonarchimedean descent datum over $B$. Then, if $\tilde{L}_A (\lambda , c)$ has a nontrivial global section, the map $H^0 (A , \tilde{L}_A (\lambda , c)) \to H^0 (X,L(\lambda , c_{\trop}))_{\Gamma}$ given by $f \mapsto f_{\trop}$ is surjective.
\end{Theorem}

In the proof of Theorem~\ref{thm:surjective:tropicalization:theta}, we use the following proposition to reduce the surjectivity in the theorem to that when $\tilde{L}_A (\lambda , c)$ is ample.

\begin{Proposition} \label{prop:Nori's_remark}
Let $\tilde{L}$ be a rigidified line bundle on an abelian variety $A$. Suppose that $H^0 (A,\tilde{L}) \neq 0$. Then there exist a surjective homomorphism of abelian varieties $\phi : A \to A_0$ and a rigidified ample line bundle $\tilde{L}_0$ on $A_0$ such that $\phi$ has reduced connected fiber and $\tilde{L} = \phi^{\ast} (\tilde{L}_0)$. Further, $\phi$ induces an isomorphism $\phi^{\ast} : H^0 (A_0 , \tilde{L}_0 ) \to H^0 (A,\tilde{L})$.
\end{Proposition}

\Proof
First, we assume that there exists an $s \in H^0 (A,\tilde{L}) \setminus \{ 0 \}$ such that $s(0) \neq 0$. Set $D := \zero (s)$. Then this divisor does not pass through $0 \in A$.
There exists a reduced subgroupscheme $H(D)$ of $A$ such that $H(D)(k) = \{ x \in A(k) \mid T_{x}^{\ast}(D) = D  \}$; see the theorem in Nori's remark on page 84 in \cite{Mumford-AV}. Let $G$ be the connected component of $H(D)$ with $0 \in G$. Then $G$ naturally acts on $A$, and let $\mu : G \times A \to A$ denote this natural action. This is the restriction of the addition $A \times A \to A$ to $G \times A$.

We define a $G$-linearlization on $\tilde{L}$. First, we note that $G \cap \Supp (D) = \emptyset$. Indeed, if $x \in G(k) \cap \Supp (D)$, then since $D -x = T_x^{\ast}(D)  = D$, we have $0 \in \Supp (D)$, which contradicts the assumption. 

Let $A'$ be the dual abelian variety of $A$ and let $\tilde{P}$ be the rigidified Poincar\'e line bundle on $A' \times A$. We have a unique homomorphism $\varphi_{\tilde{L}}^{-} : A \to A'$ such that
\[
m^{\ast} (\tilde{L}) = p_1^{\ast} (\tilde{L}) \otimes p_2^{\ast} (\tilde{L}) \otimes (\varphi_{\tilde{L}}^{-} \times \id_A)^{\ast} (\tilde{P})^{\otimes -1}
.
\]
Note here that this is an isomorphism of rigidified line bundles. Note also that $\varphi_{\tilde{L}}^{-} (G) = \{ 0 \}$ by the definition of $G$. We consider the restriction of this isomorphism over $G \times A$. Since $\varphi_{\tilde{L}}^{-} (G) = \{ 0 \}$, we have $(\varphi_{\tilde{L}}^{-} \times \id_A)^{\ast} (P)^{\otimes -1} = \OO_{G \times A}$. Since $G \cap \Supp (D) = \emptyset$, the rigidified line bundle $\rest{\tilde{L}}{G}$ is trivial, and hence $\rest{p_1^{\ast} (\tilde{L})}{G \times A} = \OO_{G \times A}$. It follows that
\[
\mu^{\ast} ( \tilde{L} ) = q_2^{\ast} (\tilde{L})
,
\]
where $q_2 : G \times A \to A$ is the second projection.
Since the above identity is an isomorphism of rigidified line bundles and since the morphisms $G \times G \times A \to G \times A$ arising from the addition on $G$, the action of $G$ on $A$, and the projection of the last two factors, respectively, are homomorphism of abelian varieties, we see that the above identity satisfies so-called the cocycle condition. Thus it is indeed a $G$-linearlization on $\tilde{L}$.

Let $\phi : A \to A_0$ be the quotient of $A$ by $G$. By the descent theory, $\tilde{L}$ descends to a unique (up to canonical isomorphims) rigidified line bundle $\tilde{L}_0$ on $A_0$; in particular, $\tilde{L} = \phi^{\ast} (\tilde{L}_0)$. Since $\phi$ has connected reduced fibers, the last assertion of the theorem holds. It remains to show that $\tilde{L}_0$ is ample. Note that $\tilde{L}^{\otimes 2}$ is basepoint free, and we take a homomorphism $h : A \to \PP^m$ arising from the complete linear system associated to $\tilde{L}^{\otimes 2}$. Then by the theorem in Nori's remark on page 84 in \cite{Mumford-AV}, $h(G)$ is a singleton. Further, since $G$ is reduced, $h(G)$ is a reduced. It follows that there exists a morphism $h_0 : A_0 \to \PP^m$ such that $h_0 \circ \phi = h$. By \cite[Theorem on page 84]{Mumford-AV} again, $h_0$ is finite. Thus the line bundle $\tilde{M}:= h_{0}^{\ast} (\OO (1))$ on $A_0$ is ample. Since $\phi$ has reduced connected fibers and $A$ is normal, it follows from the identities $\phi^{\ast} (\tilde{M}) = \tilde{L}^{\otimes 2} = \phi^{\ast} \left(\tilde{L}_0^{\otimes 2}\right)$ that $\tilde{M} = \tilde{L}_0^{\otimes 2}$. This shows that $\tilde{L}_0$ is ample. This completes the proof of the proposition when $\tilde{L}$ has a global section that does not vanish at $0 \in A$.

Next, we assume that for any $s \in H^0 (A , \tilde{L})$, $s(0) = 0$. We ignore the given rigidification on $\tilde{L}$ for a while. Since $H^0 (A ,\tilde{L}) \neq 0$, there exists an $s \in H^0 (A ,\tilde{L}) \setminus \{ 0 \}$. We set $D = \zero (s)$ again. There exists a point $z \in A(k)$ such that $T_z^{\ast} (s)$ does not vanish at $0$ as a global section of the line bundle $T_z^{\ast} (L)$. We put a rigidification on $T_z^{\ast} (\tilde{L})$. Since this line bundle has global section that does not vanish at $0$, it follows from what we have proved above that there exist an abelian subvariety $G$ of $A$ and a rigidified ample line bundle $\tilde{L}_0'$ on $A_0$, where $\phi : A \to A_0$ is the quotient of $A$ by G, such that $\phi^{\ast} (\tilde{L}_0') = T_z^{\ast} (\tilde{L})$. This isomorphism of line bundles, induces an isomorphism $\phi^{\ast} (T_{\phi(-z)}^{\ast}( \tilde{L}_0')) \cong \tilde{L}$. We set $\tilde{L}_0 := T_{\phi(-z)}^{\ast}( \tilde{L}_0')$ and put a rigidification on it in such a way that the isomorphism $\phi^{\ast} (\tilde{L}_0) \cong \tilde{L}$ is an isomorphism of rigidified line bundles. Then this $\tilde{L}_0$ suffices.
\QED

\textsl{Proof of Theorem~\ref{thm:surjective:tropicalization:theta}.}
First, assume that $\tilde{L}_A (\lambda , c)$ is ample. Then $L$ is ample and $\lambda : M' \to M$ is a polarization by Theorem~\ref{thm:NAH}, where $M$ is the character lattice of the torus part $T$ of $A$. Let $\mathfrak{B}$ be a complete system of representatives of $M / \lambda (M')$. 

By Lemma~\ref{lemma:tropicalization:Gamma-rational}, the tropical descent datum $(\lambda , c_{\trop})$ is $\Gamma$-rational. We take any $\vartheta \in H^0 (X,L(\lambda , c_{\trop}))_{\Gamma} \setminus  \{ + \infty \}$. By Proposition~\ref{prop:Gammarationalbasis}, there exists a $(c_{b})_{b \in \mathfrak{B}} \in \overline{\Gamma}^{\mathfrak{B}}$ such that
\[
\vartheta = \min_{b \in \mathfrak{B}} \left( c_b + \vartheta_{b}^{(\lambda , c_{\trop})} \right)
.
\]
Since $N_{\RR} / M'$ is compact, there exists a polyhedron $\varPi \subseteq N_{\RR}$ such that $\varPi + M' = N_{\RR}$. By the quasi-periodicity, it suffices to show that there exists a nonarchimedean theta function $\tilde{\vartheta}$ such that $\tilde{\vartheta}_{\trop}$ coincides with $\vartheta$ on $\varPi$. For any $b \in \mathfrak{B}$, Corollary~\ref{cor:preciselift} gives us a nonarchimedean theta function $\vartheta_{b}^{(\lambda , c)}$ such that $\left( \tilde{\vartheta}_{b}^{(\lambda ,c)} \right)_{\trop} = \vartheta_{b}^{(\lambda , c_{\trop})}$. We fix a $\tilde{c}_b^{(0)} \in k$ such that $- \log |\tilde{c}_b^{(0)}| = c_b$. For any $b \in \mathfrak{B}$, we have 
\begin{align} 
\label{align:tropicalization:basis2}
c_b + \vartheta_b^{(\lambda , c_{\trop})} = \left( \tilde{c}_b^{(0)} \tilde{\vartheta}_b^{(\lambda , c)} \right)_{\trop}
.
\end{align}

We take a finite polyhedral decomposition $\Sigma$ of $\varPi$ such that for each $\Delta \in \Sigma$, $\vartheta$ is $\ZZ$-affine over $\Delta$ and there exists a $b_{\Delta} \in \mathfrak{B}$ such that $c_{b_{\Delta}} + \vartheta_{b_{\Delta}}^{(\lambda , c_{\trop})} = \vartheta$ on $\Delta$. Note that 
\begin{align} \label{align:ineq:bDelta}
c_{b} + \vartheta_{b}^{(\lambda , c_{\trop})} \geq c_{b_{\Delta}} + \vartheta_{b_{\Delta}}^{(\lambda , c_{\trop})}.
\end{align}

Let $\Sigma^{(n)}$ be the subset of $\Sigma$ consisting of the $n$-dimensional faces, where $n= \dim (N_{\RR})$. For each $\Delta \in \Sigma^{(n)}$, we fix a $v_{\Delta} \in \relin (\Delta)$, where $\relin (\Delta)$ denotes the relative interior of $\Delta$. Let $\sigma_{E^{\an}} : N_{\RR} \to E^{\an}$ be the canonical section of the valuation map for $E^{\an}$. 
For $b \in \mathfrak{B}$ and $\Delta \in \Sigma^{(n)}$, we consider 
\[
\alpha_{b, \Delta} := \frac{\tilde{c}_b^{(0)} \tilde{\vartheta}_{b}^{(\lambda, c)} (\sigma_{E^{\an}} (v_{\Delta}))}{ \tilde{c}_{b_{\Delta}}^{(0)}\tilde{\vartheta}_{b_{\Delta}}^{(\lambda, c)} (\sigma_{E^{\an}} (v_{\Delta}))}
\in \mathscr{H} (\sigma_{E^{\an}} (v_{\Delta}))
,
\]
where $\mathscr{H} (\sigma_{E^{\an}} (v_{\Delta}))$ is the completed residue field of $E^{\an}$ at the point $\sigma_{E^{\an}} (v_{\Delta}) \in E^{\an}$. Since we have (\ref{align:tropicalization:basis2}) and (\ref{align:ineq:bDelta}), we see that 
\[
|\alpha_{b,\Delta}| =
\left| \frac{\tilde{c}_b^{(0)} \tilde{\vartheta}_{b}^{(\lambda, c)} (\sigma_{E^{\an}} (v_{\Delta}))}{ \tilde{c}_{b_{\Delta}}^{(0)}\tilde{\vartheta}_{b_{\Delta}}^{(\lambda, c)} (\sigma_{E^{\an}} (v_{\Delta}))} \right| \leq 1
\] 
for any $b \in \mathfrak{B}$ and $\Delta \in \Sigma^{(n)}$. Let $\widetilde{\alpha_{b, \Delta}}$ be the residue class in the residue field $\widetilde{\mathscr{H} (\sigma_{E^{\an}} (v_{\Delta}))}$ of $\mathscr{H} (\sigma_{E^{\an}} (v_{\Delta}))$. Note that $\widetilde{\alpha_{b_{\Delta} , \Delta}} = 1$.

Let $(z_b)_{b \in \mathfrak{B}}$ be a system of indeterminates indexed by $\mathfrak{B}$. Consider a finite system of linear forms
\[
\left( \sum_{b \in \mathfrak{B}} z_b \widetilde{\alpha_{b, \Delta}}
\right)_{\Delta \in \Sigma^{(n)}}
\]
on $(z_b)_{b \in \mathfrak{B}}$ indexed by $\Sigma^{(n)}$. Since each linear form is nontrivial and the residue field $\widetilde{k} := k^{\circ} / k^{\circ \circ}$ is algebraically closed and hence infinite, there exists a $(\lambda_b)_{b \in \mathfrak{B}} \in (k^{\circ})^{\mathfrak{B}}$ such that for any $\Delta \in \Sigma^{(n)}$, 
$\sum_{b \in \mathfrak{B}} \widetilde{\lambda_b} \widetilde{\alpha_{b, \Delta}}
\neq 0$, where $\widetilde{\lambda_b} \in \widetilde{k}$ is the residue class of $\lambda_b$. Thus for any $\Delta \in \Sigma^{(n)}$, we have 
\begin{align} \label{align:solution:linear:eq}
\left| \sum_{b \in \mathfrak{B}} \lambda_b \alpha_{b, \Delta} \right| = 1
.
\end{align}

We take any $\Delta \in \Sigma^{(n)}$ and consider 
\[
\eta := \rest{\left( \sum_{b \in \mathfrak{B}} \lambda_i \tilde{c}_b^{(0)} \tilde{\vartheta}_{b}^{(\lambda,c)} \right)_{\trop}}{\Delta}
-
\rest{\vartheta}{\Delta}
.
\] 
Since $\left( \sum_{b \in \mathfrak{B}} \lambda_i \tilde{c}_b^{(0)} \tilde{\vartheta}_{b}^{(\lambda,c)} \right)_{\trop}$ is concave and $\rest{\vartheta}{\Delta}$ is linear, $\eta$ is concave. Noting (\ref{align:tropicalization:basis2}) and Remark~\ref{remark:tropicalization-inequality}, we have
\begin{align*} 
\left( \sum_{\mathfrak{B}} \lambda_b \tilde{c}_b^{(0)} \tilde{\vartheta}_{b}^{(\lambda,c)} \right)_{\trop}
\geq 
\min_{b \in \mathfrak{B}} \left( c_b + \vartheta_{b}^{(\lambda,c_{\trop})} \right) = \vartheta
,
\end{align*}
and thus $\eta \geq 0$.
Further, by (\ref{align:solution:linear:eq}), we have
\[
\left| \sum_{b \in \mathfrak{B}} \lambda_b 
\tilde{c}_b^{(0)}  \tilde{\vartheta}_{b}^{(\lambda ,c)} (\sigma_{E^{\an}} (v_{\Delta}))  \right|  = \left| \tilde{c}_{b_{\Delta}}^{(0)}\tilde{\vartheta}_{b_{\Delta}}^{(\lambda, c)} (\sigma_{E^{\an}} (v_{\Delta})) \right|
,
\]
and hence
\begin{multline} \label{align:trop:equal:v}
\left( \sum_{b \in \mathfrak{B}} \lambda_b \tilde{c}_b^{(0)} \tilde{\vartheta}_{b}^{(\lambda,c)} \right)_{\trop} (v_{\Delta}) 
=
- \log \left| \sum_{b \in \mathfrak{B}} \lambda_b 
\tilde{c}_b^{(0)}  \tilde{\vartheta}_{b}^{(\lambda,c)} (\sigma_{E^{\an}} (v_{\Delta}))  \right| 
\\
=
- \log \left| \tilde{c}_{b_{\Delta}}^{(0)}\tilde{\vartheta}_{b_{\Delta}}^{(\lambda,c)} (\sigma_{E^{\an}} (v_{\Delta})) \right| = c_{b_{\Delta}} + \vartheta_{b_{\Delta}}^{(\lambda,c_{\trop})} (v_{\Delta}) = \vartheta (v_{\Delta}).
\end{multline}
Thus $\eta$ is a nonnegative concave function on $\Delta$ such that $\eta (v_{\Delta}) = 0$ with $v_{\Delta} \in \relin (\Delta)$. It follows that $\eta$ is a constant function $0$, i.e., 
\[
\rest{\left( \sum_{b \in \mathfrak{B}} \lambda_i \tilde{c}_b^{(0)} \tilde{\vartheta}_{b}^{(\lambda,c)} \right)_{\trop}}{\Delta}
=
\rest{\vartheta}{\Delta}
.
\]
Since we have taken $\Delta \in \Sigma^{(0)}$ arbitrarily and since $\varPi = \bigcup_{\Delta \in \Sigma^{(n)}} \Delta$, this shows that 
\[
\rest{\left( \sum_{b \in \mathfrak{B}} \lambda_i \tilde{c}_b^{(0)} \tilde{\vartheta}_{b}^{(\lambda,c)} \right)_{\trop}}{\varPi}
=
\rest{\vartheta}{\varPi}
.
\]
Thus the proof when
$\tilde{L}_A (\lambda , c)$ is ample is complete.

We consider the general case. Since $\tilde{L}_A (\lambda,c)$ has a nontrivial global section, so does $L$ on $B$ by \cite[Remark~5.4]{BoschLutke-DAV}. By Proposition~\ref{prop:Nori's_remark}, there exist an abelian variety $A_0$ over $k$, a rigidified ample line bundle $\tilde{L}_0$ on $A_0$ and a surjective homomorphism $\phi : A \to A_0$ by a such that $\phi^{\ast} (\tilde{L}_0) = \tilde{L}_A (\lambda , c)$ and the induced homomorphism $\phi^{\ast} : H^0 (A_0 , \tilde{L}_0) \to H^0 (A , \tilde{L}_A (\lambda , c))$ is an isomorphism. Let 
\begin{align*} 
\begin{CD}
@. M'_0 \\
@. @VV{\widetilde{\Phi}_0}V \\
 T_0^{\an} @>>> E_0^{\an} @>{q_0}>> B_0^{\an}  \\
@. @V{p_0}VV \\
@. A_0^{\an}
\end{CD}
\end{align*}
be the Raynaud cross for $A_0$. Let $M_0$ be the character lattice of $T_0$, and recall that $M'$ is a lattice in $N_{\RR} := \mathrm{Hom} (M,\RR)$. By Theorem~\ref{thm:NAH}, there exists a nonarchimedean descent datum $(L_0 , \lambda_0 , c_0)$ over $B_0$ such that $\tilde{L}_0 = \tilde{L_0}_{A_0} (\lambda_0 , c_0)$. By \cite[Proposition~6.10]{BoschLutke-DAV}, $\phi$ lifts to a homomorphism $\psi : E \to E_0$, and $\psi$ corresponds to a pair $(\mu , \psi_{\mathrm{ab}})$, where $\mu : M_0 \to M$ is a homomorphism between the character lattices that induces the homomorphism $T \to T_0$ that is the restriction of $\psi$ and $\psi_{\mathrm{ab}} : B \to B_0$ is a homomorphism such that $\psi_{\mathrm{ab}} \circ q = q_0 \circ \psi$. Further, there exists a homomorphism $\mu' : M' \to M_0'$ such that $\psi \circ \widetilde{\Phi} = \widetilde{\Phi}_0 \circ \mu'$. Since $\phi^{\ast} (\tilde{L}_0) = \tilde{L}$, it follows from \cite[Theorem~6.7]{BoschLutke-DAV} that $\mu \circ \lambda_0 \circ \mu'$ and that there exists a $u \in M$ such that $L = \psi_{\mathrm{ab}}^{\ast} (\tilde{L}_0)  \otimes E_u$ and $c^{\otimes -1} = (\mu')^{\ast} (c_0)^{\otimes -1}  \otimes \widetilde{\Phi}^{\ast} (q^{\ast} (e_u))$. Moreover, for any nonarchimedean theta function $\tilde{\vartheta}$ with respect to $(L_0 , \lambda_0, c_0)$, we see that $\phi^{\ast} (\vartheta) \otimes q^{\ast} (e_u)$ is a nonarchimedean theta function with respect to 
\[
(\psi_{\mathrm{ab}}^{\ast} (L_0) \otimes q^{\ast} (E_u), \mu \circ \lambda_0 \circ \mu' ,(\mu')^{\ast} (c_0) \otimes \widetilde{\Phi} (q^{\ast} (e_u))^{\otimes -1}) = (L,\lambda,c).
\]
From those descriptions, we see that the isomorphism $H^0 (A_0 , \tilde{L_0} (\lambda_0 , c_0)) \to H^0 (A, \tilde{L} (\lambda , c))$ 
is given by $\tilde{\vartheta_0} \mapsto \phi^{\ast}(\tilde{\vartheta_0}) \otimes q^{\ast} (e_u)$. 

Set $X_0 := (N_{0})_{\RR} / M_0'$, which is the canonical tropicalization of $A_0$. The homomorphism $\mu : M_0 \to M$ induces an $\RR$-linear map $N_{\RR} \to (N_0)_{\RR}$, which restricts to $\mu' : M' \to M_0'$. Thus $\phi$ induces a homomorphism $\phi_{\trop} : X \to X_0$. Noting the above description of the isomorphism $H^0 (A_0 , \tilde{L_0} (\lambda_0 , c_0)) \to H^0 (A, \tilde{L} (\lambda , c))$, the diagram 
\[
\begin{CD}
H^0 (A_0 , \tilde{L_0}_{A_0} (\lambda_0 , c_0)) @>{\phi^{\ast} (\ndot) \otimes \tilde{\Phi}^{\ast} (q^{\ast} (e_u))}>> H^0 (A , \tilde{L}_A (\lambda , c)) \\
@VVV @VVV \\
H^0 (X_0 , L (\lambda_0 , (c_{0})_{\trop}) ) @>>{\phi_{\trop}^{\ast} (\ndot) + u}> H^0 (X , L(\lambda , c_{\trop}))
\end{CD}
\]
commutes, where $L(\lambda , c_{\trop})$ and $L (\lambda_0 , (c_{0})_{\trop})$ are the line bundles on $X$ and $X_0$ associated to the tropical descent data $(\lambda , c_{\trop})$ and $ (\lambda_0 , (c_{0})_{\trop})$, respectively, and the vertical arrows are the tropicalization maps. Since $\tilde{L_0}_{A_0} (\lambda_0 , c_0)$ is ample, it follows from what we have just proved above that the vertical arrow on the left-hand side is surjective. Since the horizontal arrows are isomorphisms (see Remark~\ref{remark:translate:u} for the bottom horizontal arrow), the vertical arrows on the right-hand side is also surjective. This completes the proof.
\QED

\end{document}